\documentclass{article}

\usepackage[letterpaper, margin=1in]{geometry}
\usepackage{amsmath,amssymb}
\usepackage{mathrsfs}
\usepackage{graphicx}
\usepackage{booktabs}
\usepackage[amsmath,thmmarks,framed]{ntheorem}
\usepackage{authblk}
\usepackage{xcolor}
\usepackage{comment}
\definecolor{darkgreen}{rgb}{0.0,0.5,0.0}
\usepackage[colorlinks,linkcolor=red,citecolor=darkgreen,urlcolor=blue]{hyperref}
\hypersetup{
  pdftitle={Asymptotic Variance and Quantitative Central Limit Theorem for Underdamped Langevin Dynamics on Riemannian Manifolds}
}
\usepackage[numbers,square,comma,sort&compress]{natbib}

\newtheorem{theorem}{Theorem}
\newtheorem{lemma}{Lemma}
\newtheorem{corollary}[theorem]{Corollary}
\newtheorem{assumption}{Assumption}
\newtheorem{definition}{Definition}
\newtheorem{remark}{Remark}
\newenvironment{proof}[1][Proof]{\par\medskip\noindent\textit{#1.}\ }{\hfill$\square$\par\medskip}

\newcommand{\Le}{\leqslant}
\newcommand{\Ge}{\geqslant}
\newcommand{\D}{\mathrm{d}}
\newcommand{\norm}[1]{\lVert #1 \rVert}
\newcommand{\Ric}{\mathrm{Ric}}

\title{Asymptotic Variance and Quantitative Central Limit Theorem for Underdamped Langevin Dynamics on Riemannian Manifolds}
\author[1]{Xuda Ye}
\author[2]{Yuliang Wang}
\author[3]{Zhihan Zhang}
\affil[1]{Department of Mathematics, Purdue University, West Lafayette, IN 47907, USA}
\affil[2]{Department of Mathematics, Duke University, Durham, NC 27708, USA}
\affil[3]{Department of Computer Science and Operations Research, Université de Montréal, Montréal, QC H3T 1N8, Canada}
\date{}

\numberwithin{equation}{section}

\begin{document}
\maketitle

\begin{abstract}
For the underdamped Langevin dynamics with the damping rate $\gamma$ on a complete Riemannian manifold, we study the convergence of the time averages of test functions of the position. We assume a Poincaré inequality with the constant $m$ and a lower bound $-K$ of the sum of the Ricci curvature and the Hessian of the potential, and we prove the following results. First, the semigroup decays in $L^2$ with a rate $\lambda$, which satisfies
\begin{equation*}
    \lambda \Ge \frac{m\gamma}{\gamma^2 + 2 \sqrt{\frac{\sqrt5 - 1}{2}} \gamma\sqrt m + \frac{3 + \sqrt5}{2} m + K} .
\end{equation*}
Second, we define the normalized variance $\mathbb V(\gamma)$ as the largest ratio of the asymptotic variance to the variance of the test function. The normalized variance bounds the mean-square error of the time averages, and we bound $\mathbb V(\gamma)$ from below and from above. Third, the mixing time $t_{\mathrm{mix}}(\gamma)$ of the dynamics is equivalent to $\mathbb V(\gamma) + \gamma^{-1}$ up to universal constants. Finally, we establish a quantitative central limit theorem: the Wasserstein-1 distance between the law of the normalized time average and its Gaussian limit decays at the rate $t^{-\frac12}$ in the time $t$.


\end{abstract}


\section{Introduction}
\label{sec: intro}

Let $\mathcal X$ be a complete, connected and oriented Riemannian manifold of dimension $d$ without boundary and with the tangent bundle $T\mathcal X$. Let $U(\cdot)$ be a smooth potential function on $\mathcal X$, let $\gamma > 0$ be the damping rate and consider the underdamped Langevin dynamics \citep{jorgensen1978construction, grothaus2020hypocoercivity, eberle2026space}
\begin{equation}
    \D X_t = V_t \,\D t , \qquad
    \D V_t = -\bigl( \nabla U(X_t) + \gamma V_t \bigr) \D t + \sqrt{2\gamma}\,\D B_t ,
    \label{eq: uld}
\end{equation}
where $X_t \in \mathcal X$ is the position, $V_t \in T_{X_t}\mathcal X$ is the velocity and $(B_t)_{t \Ge 0}$ is a Brownian motion in the tangent spaces $T_{X_t}\mathcal X$. For $\mathcal X = \mathbb R^d$, \eqref{eq: uld} is the underdamped Langevin dynamics in $\mathbb R^{2d}$. The dynamics \eqref{eq: uld} leaves invariant the Boltzmann distribution
\begin{equation}
    \mu(\D x\D v) \propto \exp\Bigl( -U(x) - \frac12 |v|^2 \Bigr) \D x \D v,\qquad x \in \mathcal X, \quad v \in T_x\mathcal X.
    \label{eq: boltzmann}
\end{equation}
Here $\D x$ is the Riemannian volume of $\mathcal X$ and $\D v$ is the Lebesgue measure of $T_x\mathcal X$. The marginal distribution of the position under $\mu$ is denoted by $\bar\mu(\D x) \propto \exp(-U(x))\D x$.

Let $C^1(\mathcal X)$ be the space of the continuously differentiable functions on $\mathcal X$.
For a function $\psi$ on $T\mathcal X$, denote its average under the Boltzmann distribution $\mu$ by
\begin{equation}
	\mu(\psi):=\int_{T\mathcal X} \psi(x,v)\, \mu(\D x\D v).
\end{equation}
Given a test function $f \in C^1(\mathcal X)$, we aim to quantify how fast the time average $t^{-1}\int_0^t f(X_s)\D s$ converges to $\mu(f)$, with explicit dependence on the damping rate $\gamma>0$.

The major approach for quantifying the time average convergence is the Poisson equation \citep{bhattacharya1982functional, mattingly2010convergence, komorowski2012fluctuations}. Let $\mathcal L$ be the generator of the dynamics \eqref{eq: uld} and $(e^{t\mathcal L})_{t \Ge 0}$ be the associated semigroup. Suppose the test function $f$ has a bounded gradient, $\norm{\nabla f}_\infty := \sup_{x \in \mathcal X} |\nabla f(x)| < \infty$. The (strong) Poisson equation is defined as
\begin{equation}
	\mathcal L \phi = v \cdot \nabla_x\phi - \nabla U \cdot \nabla_v\phi - \gamma v \cdot \nabla_v\phi + \gamma \Delta_v\phi = \mu(f) - f,
	\label{eq: poisson equation}
\end{equation}
where $\nabla_x$ and $\nabla_v$ are the gradients in the position and in the velocity and $\Delta_v$ is the Laplacian in the velocity. With the decay of the semigroup $(e^{t\mathcal L})_{t\Ge0}$ in $L^2(\mu)$ \citep{cao2023explicit, eberle2026space}, the Poisson solution can be represented as
\begin{equation}
    \phi = \int_0^\infty \bigl( e^{t\mathcal L}f - \mu(f) \bigr) \D t .
    \label{eq: poisson solution}
\end{equation}
The Poisson solution allows a convenient It\^o formula \citep{bhattacharya1982functional, kipnis1986central, komorowski2012fluctuations}
\begin{equation}
	\int_0^t \bigl( f(X_s) - \mu(f) \bigr) \D s = \phi(X_0,V_0) - \phi(X_t,V_t) + \sqrt{2\gamma} \int_0^t \nabla_v\phi(X_s,V_s) \cdot \D B_s ,
	\label{eq: poisson decomposition}
\end{equation}
where the stochastic integral is a martingale with the asymptotic variance $\sigma_f^2 := 2\gamma\, \mu\bigl( |\nabla_v\phi|^2 \bigr)$. If the initial state $(X_0,V_0) \sim \mu$ is independent of $(B_t)_{t \Ge 0}$, then we have formally the central limit theorem
\begin{equation}
    \frac{1}{\sqrt t} \int_0^t \bigl( f(X_s) - \mu(f) \bigr) \D s \xrightarrow{d} \mathcal N(0, \sigma_f^2) ,
    \label{eq: clt convergence}
\end{equation}
where $\xrightarrow{d}$ denotes the convergence in distribution as $t \to \infty$, and $\mathcal N(0,\sigma_f^2)$ denotes the point mass at $0$ if $\sigma_f^2 = 0$. Classical results \citep{bhattacharya1982functional, cattiaux2012central, komorowski2012central} prove convergences of the form \eqref{eq: clt convergence}.

The classical central limit theorem \eqref{eq: clt convergence} does not indicate how fast the time average converges in $t$, nor the dependence of the convergence on parameters such as the damping rate $\gamma$. To state our results, we summarize the assumptions on the manifold $\mathcal X$ and the potential function $U(\cdot)$:
\begin{itemize}
	\setlength{\itemsep}{0pt}
	\item For some $K\Ge0$, the Ricci curvature $\Ric(\cdot)$ of $\mathcal X$ and the Hessian $\nabla^2 U$ of the potential function satisfy
	\begin{equation*}
		w \cdot \bigl( \Ric + \nabla^2 U \bigr)(x)\, w \Ge -K |w|^2 , \qquad \forall x \in \mathcal X ,\quad w \in T_x\mathcal X .
	\end{equation*}
	\item For some $m>0$, the distribution $\bar\mu$ satisfies a Poincaré inequality with the smallest constant $m^{-1}$. That is, $m^{-1}$ is the smallest constant such that, for any $g\in C^1(\mathcal X)$ with $\norm{\nabla g}_\infty < \infty$,
	\begin{equation*}
		\mu\bigl( (g - \mu(g))^2 \bigr) \Le \frac1m\, \mu\bigl( |\nabla g|^2 \bigr) .
	\end{equation*}
\end{itemize}

Our contributions concern four aspects.
First, Theorem~\ref{thm: decay} proves the decay of the semigroup in $L^2(\mu)$ with an explicit lower bound of a new rate $\lambda$: for every $\psi \in L^2(\mu)$ with $\mu(\psi) = 0$ and every $t \Ge 0$,
\begin{equation*}
    \mu\bigl( ( e^{t\mathcal L}\psi )^2 \bigr) \Le e^{1 - 2\lambda t}\, \mu(\psi^2) , \qquad
    \lambda \Ge \frac{m\gamma}{\gamma^2 + 2 \sqrt{\frac{\sqrt5 - 1}{2}} \gamma\sqrt m + \frac{3 + \sqrt5}{2} m + K} .
\end{equation*}
Theorem~\ref{thm: decay} achieves a better rate than the rates in \citep{fan2026sharp, eberle2026space}.
Second, we define the normalized variance
\begin{equation*}
    \mathbb V(\gamma) := \sup_{f \in L^2(\bar\mu)} \frac{\sigma_f^2}{\mu\bigl( ( f - \mu(f) )^2 \bigr)} ,
\end{equation*}
and Theorem~\ref{thm: poisson} and Corollary~\ref{cor: second moment} bound the mean-square error with $(X_0,V_0) \sim \mu$: for every $t > 0$,
\begin{equation*}
	\sup_{\mu( ( f - \mu(f) )^2 ) \Le 1} \mathbb E\biggl[ \Bigl( \frac1t \int_0^t f(X_s) \,\D s - \mu(f) \Bigr)^2 \biggr] \Le \Bigl( \frac1t + \frac{2}{\gamma t^2} \Bigr) \mathbb V(\gamma) \Le \Bigl( \frac1t + \frac{2}{\gamma t^2} \Bigr) \frac{2 ( \gamma^2 + m + K )}{m\gamma} .
\end{equation*}
Third, Theorem~\ref{thm: mixing} bounds the mixing time $t_{\mathrm{mix}}(\gamma)$, the smallest $t \Ge 0$ with $\norm{e^{t\mathcal L}\psi}_{L^2(\mu)} \Le e^{-1} \norm{\psi}_{L^2(\mu)}$ for every $\psi \in L^2(\mu)$ with $\mu(\psi) = 0$, from below and from above by the normalized variance,
\begin{equation*}
    \frac{1 - e^{-1}}{2} \max\Bigl( \mathbb V(\gamma) , \frac1\gamma \Bigr) \Le t_{\mathrm{mix}}(\gamma) \Le \Bigl( 1 + \frac\pi2 \Bigr) \Bigl( \mathbb V(\gamma) + \frac3\gamma \Bigr) .
\end{equation*}
Finally, with a universal constant $C$, Theorem~\ref{thm: clt} bounds the Wasserstein-1 distance between the law of the left side of \eqref{eq: clt convergence} and $\mathcal N(0,\sigma_f^2)$,
\begin{equation*}
    \mathcal W_1\biggl( \mathrm{Law}\Bigl( \frac{1}{\sqrt t} \int_0^t \bigl( f(X_s) - \mu(f) \bigr) \D s \Bigr) , \mathcal N(0, \sigma_f^2) \biggr)
    \Le \frac{C ( \gamma + \sqrt m + \sqrt K )^2}{m^{\frac32}\gamma} \Bigl( 1 + \log\frac{\gamma + \sqrt m + \sqrt K}{\gamma} \Bigr) \frac{\norm{\nabla f}_\infty}{\sqrt t}.
\end{equation*}
To the best of our knowledge, Theorem~\ref{thm: clt} is the first quantitative central limit theorem obtained for the underdamped Langevin dynamics with the explicit rate $t^{-\frac12}$. Related quantitative central limit theorems are known for Markov chains in discrete time \citep{srikant2024rates, zhang2026wasserstein}, for martingales \citep{haeusler1988rate, guo2024rate} and for a kinetic diffusion with a periodic force \citep{hairer2004periodic}.

The paper is organized as follows. Section~\ref{sec: related} reviews related works. Section~\ref{sec: preliminaries} states the background on the manifold, the dynamics and the Poisson equation. Section~\ref{sec: notation} fixes the remaining notations and states the assumptions. Section~\ref{sec: results} states the results. Section~\ref{sec: conclusion} concludes. The appendices prove the results of Section~\ref{sec: results}, in the order of its subsections.

\section{Related works}
\label{sec: related}

This section reviews two lines of work: the accelerated convergence of the underdamped Langevin dynamics, and the convergence of time averages by the Poisson equation and Stein's method.

\subsection{Accelerated convergence of the underdamped Langevin dynamics}
\label{sec: related convergence}

The theory of hypocoercivity gives exponential convergence of the underdamped Langevin dynamics in a modified $H^1$ norm \citep{villani2009hypocoercivity} and in a weighted $L^2$ norm \citep{dolbeault2015hypocoercivity}. Coupling methods give contraction in Wasserstein distances for strongly convex potentials \citep{cheng2018underdamped, dalalyan2020sampling} and for nonconvex potentials \citep{eberle2019couplings, schuh2024global}. In $\mathbb R^d$, under the lower bound $-K$ of the Hessian of $U$ and the Poincaré inequality with the constant $m^{-1}$, \citet{cao2023explicit} prove that the semigroup decays in $L^2(\mu)$ at the rate
\begin{equation}
    \frac{m\gamma}{C ( \gamma + \sqrt m + \sqrt K )^2}
    \label{eq: euclidean rate}
\end{equation}
with a universal constant $C$, by the variational approach of \citet{albritton2024variational}. This is the first result to show accelerated convergence for convex potentials: for $K = 0$ and $\gamma = \sqrt m$, the rate is a multiple of $\sqrt m$, while the semigroup of the overdamped Langevin dynamics decays at the rate $m$. \citet{fan2026sharp} prove the rate \eqref{eq: euclidean rate} with $C = 6$. For a class of lifts of reversible diffusions that contains the underdamped Langevin dynamics, \citet{brigati2026hypocoercivity} recover the order $\sqrt m$ of the rate for convex potentials. For convex potentials under a logarithmic Sobolev inequality, \citet{lu2026sharp} proves the decay of the relative entropy and \citet{li2026space} prove the hypercontractivity of the semigroup, both at a rate of the order of the square root of the logarithmic Sobolev constant. Under a logarithmic Sobolev inequality and a lower bound of the Hessian of $U$, \citet{lu2026renyi} prove the decay of the Rényi divergence of order $2$ for the underdamped Langevin dynamics and for piecewise deterministic samplers, at a rate that depends on the initial divergence. \citet{du2026sharp} prove a decay rate of the order $\sqrt m$ for convex potentials in bounded convex domains of $\mathbb R^d$ with specular reflection, and \citet{eberle2026space} prove a decay rate of the form \eqref{eq: euclidean rate} on Riemannian manifolds with a locally convex boundary and a reflection of the velocity at the boundary. \citet{grothaus2020hypocoercivity} prove exponential convergence in $L^2(\mu)$ for the underdamped Langevin dynamics on the tangent bundle of a complete Riemannian manifold. \citet{villani2026kinetic} proves exponential convergence in $L^2$ for the kinetic Fokker--Planck equation without a potential on the tangent bundle of a compact Riemannian manifold, and \citet{kong2024convergence} prove convergence in the Wasserstein-2 distance for the underdamped Langevin dynamics on compact Lie groups. \citet{lelievre2012langevin} study the Langevin dynamics constrained to a submanifold, and \citet{angst2015kinetic} study the kinetic Brownian motion on the unit tangent bundle of a Riemannian manifold.

Lower bounds of the relaxation time show that the order $\sqrt m$ of the decay rate cannot be improved. The relaxation time is the smallest time after which the semigroup contracts every function of mean $0$ by the factor $e^{-1}$ in $L^2$ of the invariant distribution. For finite Markov chains, \citet{diaconis2000analysis} show that a nonreversible sampler can improve on the diffusive behavior of reversible samplers, and \citet{chen1999lifting} bound the mixing time of a lift of a chain from below. In continuous time, \citet{eberle2024non} show that the relaxation time of the underdamped Langevin dynamics is at least a multiple of the square root of the relaxation time of the overdamped Langevin dynamics. In $\mathbb R^d$, \citet{eberle2026relaxation} bound the relaxation time of the underdamped Langevin dynamics from below by a multiple of $\gamma^{-1}$, and by a multiple of $\gamma$ divided by a constant that is at least $m$ and equals $m$ for quadratic potentials. Theorem~\ref{thm: mixing} bounds the relaxation time, which this paper calls the mixing time, from below and from above on Riemannian manifolds.

\subsection{Convergence of time averages}
\label{sec: related averages}

The convergence \eqref{eq: clt convergence} is known without a rate for stationary ergodic Markov processes and test functions in the range of the generator \citep[Theorem~2.1]{bhattacharya1982functional}, for reversible Markov processes \citep{kipnis1986central}, for ergodic diffusions including kinetic models \citep{cattiaux2012central}, and for Lipschitz test functions when the semigroup contracts exponentially in the Wasserstein-1 distance \citep[Theorem~2.1]{komorowski2012central}. \citet[Theorem~2.1]{bhattacharya1982functional} also expresses the asymptotic variance by the solution of the Poisson equation, and the monograph of \citet{komorowski2012fluctuations} presents the martingale approximation behind these results. For the underdamped Langevin dynamics, $\sigma_f^2$ is at most $2 \mu\bigl( ( f - \mu(f) )^2 \bigr)$ times the norm of the inverse of $\mathcal L$ on the functions $\psi \in L^2(\mu)$ with $\mu(\psi) = 0$. Hypocoercive estimates bound this norm by a multiple of $\max( \gamma, \gamma^{-1} )$ \citep{roussel2018spectral}, and \citet{bernard2022hypocoercivity} make the bound explicit in $\gamma$ and in the Poincaré constant. \citet{chak2023optimal} study the asymptotic variance as a function of the friction matrix, and compute the optimal friction for Gaussian distributions and quadratic test functions. Theorem~\ref{thm: poisson} and Corollary~\ref{cor: variance} bound $\sigma_f^2$ from above and from below with explicit constants.

The Poisson equation also underlies the error analysis of time averages and invariant distributions of numerical integrators \citep{mattingly2010convergence, leimkuhler2016computation}. Discrete Poisson equations bound the mean-square error of time averages of integrators of the overdamped Langevin dynamics \citep{ye2026mse, ye2026sglm} and of the underdamped Langevin dynamics \citep{lu2025mse}. For Markov chains in discrete time, \citet{latuszynski2013nonasymptotic} bound the mean-square error after $n$ steps, with the leading term given by the asymptotic variance of the chain divided by $n$, and \citet{joulin2010curvature} prove nonasymptotic bounds of the convergence of empirical means under a positive curvature condition. For nonreversible Markov processes, \citet{xu2026singular} bounds the variance of the time averages at a finite time from above and from below by the singular value gap of the generator, uniformly over the functions in $L^2$. Corollary~\ref{cor: second moment} bounds the mean-square error of the time average of the underdamped Langevin dynamics itself.

Rates of convergence in the central limit theorem are known for independent samples and for Markov chains. For $n$ independent samples, Stein's method gives the rate $n^{-\frac12}$ in the Wasserstein-1 distance, which is the best possible \citep[Theorem~3.2]{ross2011fundamentals}. For Markov chains in discrete time with a positive definite asymptotic covariance, \citet{srikant2024rates} solves the Poisson equation of the chain, applies Stein's method to the resulting martingale and obtains the rate $(\log n)^2 / \sqrt n$ after $n$ steps in the Wasserstein-1 distance, and \citet{zhang2026wasserstein} improve it to $n^{-\frac12}$. For the martingale central limit theorem, \citet{haeusler1988rate} proves rates of convergence in discrete and continuous time, and \citet{guo2024rate} proves rates in Wasserstein distances, where the rate in the Wasserstein-1 distance is optimal up to a constant. For ergodic Markov processes with a spectral gap, \citet{lezaud2001chernoff} proves Berry--Esseen inequalities for the empirical means. In Stein's method for the invariant distribution of an elliptic diffusion, the first and second derivatives of the Poisson solution are bounded by the Lipschitz constant of the test function \citep{gorham2019measuring, fang2019multivariate}. For the degenerate noise of the underdamped Langevin dynamics, Theorem~\ref{thm: poisson} and Corollary~\ref{cor: variance} bound the Poisson solution in $L^2(\mu)$ and $L^4(\mu)$, and Theorem~\ref{thm: clt} proves the rate $t^{-\frac12}$ in continuous time and assumes no lower bound of $\sigma_f^2$.

\section{Preliminaries}
\label{sec: preliminaries}

This section recalls the Riemannian geometry of $\mathcal X$ and of the tangent bundle $T\mathcal X$, formulates the underdamped Langevin dynamics \eqref{eq: uld}, its generator and the decay of its semigroup, and introduces the Poisson equation and the asymptotic variance. Throughout this section, the functions are assumed to be smooth and integrable enough for the identities to hold.

\subsection{Riemannian manifold and its tangent bundle}
\label{sec: manifold}

For $x \in \mathcal X$, $T_x\mathcal X$ is the tangent space of $\mathcal X$ at $x$, $a \cdot b$ denotes the inner product of $a, b \in T_x\mathcal X$ and $|a|$ the norm of $a$. For a function $g$ on $\mathcal X$, $\nabla g$ is the gradient, $\nabla^2 g$ is the Hessian with respect to the Levi-Civita connection, and $\Delta g$ is the trace of the Hessian. At each point $x$, the Hessian $\nabla^2 g$ and the Ricci curvature $\Ric(\cdot)$ of $\mathcal X$ are symmetric linear maps of $T_x\mathcal X$. For a symmetric linear map $A$ of $T_x\mathcal X$, $|A|^2$ is the trace of $A^2$. For $K \Ge 0$, the sum of the Ricci curvature and the Hessian of the potential has the lower bound $-K$ if
\begin{equation}
    w \cdot \bigl( \Ric + \nabla^2 U \bigr)(x)\, w \Ge -K |w|^2 , \qquad \forall x \in \mathcal X ,\quad w \in T_x\mathcal X .
    \label{eq: curvature}
\end{equation}
The lower bound \eqref{eq: curvature} is the curvature condition of \citet{bakry1985diffusions} with the constant $-K$.

A state is $z = (x,v) \in T\mathcal X$ with $x \in \mathcal X$ and $v \in T_x\mathcal X$, and $\D z := \D x \D v$ is the measure of \eqref{eq: boltzmann}, the product of the Riemannian volume of $\mathcal X$ in the position and the Lebesgue measure of $T_x\mathcal X$ in the velocity. The distribution of the position under $\mu$ is $\bar\mu(x) \propto e^{-U(x)}$ with respect to the Riemannian volume. A function of the position is a function $g$ on $\mathcal X$, also regarded as the function $(x,v) \mapsto g(x)$ on $T\mathcal X$, with the average $\mu(g)$. With $m > 0$, the Poincaré inequality of $\bar\mu$ with the constant $m^{-1}$ is, for a function $g \in C^1(\mathcal X)$,
\begin{equation}
    \mu\bigl( (g - \mu(g))^2 \bigr) \Le \frac1m\, \mu\bigl( |\nabla g|^2 \bigr) ,
    \label{eq: poincare}
\end{equation}
and $m^{-1}$ is the smallest such constant.
The (elliptic) diffusion operator on $\mathcal X$ is
\begin{equation}
    \bar{\mathcal L} := \Delta - \nabla U \cdot \nabla ,
    \label{eq: overdamped}
\end{equation}
which is self-adjoint in $L^2(\bar\mu)$, and $\bar\mu$ is the invariant distribution, that is, $\mu(\bar{\mathcal L} g) = 0$ for a function $g$ of the position. Since $m^{-1}$ is the smallest constant of \eqref{eq: poincare}, $m$ is the spectral gap of $\bar{\mathcal L}$ in $L^2(\bar\mu)$ \citep[Section~2]{eberle2026space}. Integration by parts on $\mathcal X$ gives, for functions $g_1$ and $g_2$ of the position,
\begin{equation}
    \mu\bigl( (-\bar{\mathcal L} g_1)\, g_2 \bigr) = \mu( \nabla g_1 \cdot \nabla g_2 ) .
    \label{eq: position form}
\end{equation}
The Bochner formula for $\bar{\mathcal L}$ \citep[Equation~(1.16.4)]{bakry2014analysis} reads, for a function $g$ of the position,
\begin{equation}
    \frac12\, \bar{\mathcal L} |\nabla g|^2 = |\nabla^2 g|^2 + \nabla g \cdot \nabla \bar{\mathcal L} g + \nabla g \cdot \bigl( \Ric + \nabla^2 U \bigr) \nabla g .
    \label{eq: bochner formula}
\end{equation}

A function $\psi$ of $(x,v) \in T\mathcal X$ has a gradient $\nabla_x\psi$ in $x$ and a gradient $\nabla_v\psi$ in $v$, both with values in $T_x\mathcal X$. The velocity gradient $\nabla_v\psi(x,v)$ is the gradient of the function $\psi(x,\cdot)$ on the Euclidean space $T_x\mathcal X$, and $\Delta_v\psi(x,v)$ is the Laplacian of $\psi(x,\cdot)$. The position gradient $\nabla_x\psi(x,v)$ is defined by
\begin{equation}
    w \cdot \nabla_x\psi(x,v) = \frac{\D}{\D s}\, \psi\bigl( \mathfrak c(s), \mathfrak q(s) \bigr) \Big|_{s=0} , \qquad \forall w \in T_x\mathcal X ,
    \label{eq: position gradient}
\end{equation}
where $\mathfrak c$ is a curve in $\mathcal X$ with $\mathfrak c(0) = x$ and the derivative $w$ at $s = 0$, and $\mathfrak q(s) \in T_{\mathfrak c(s)}\mathcal X$ is the parallel transport of $v$ along $\mathfrak c$. The right side of \eqref{eq: position gradient} is the derivative of $\psi$ along the horizontal lift of $w$ at $(x,v)$, and $\nabla_v\psi$ is the gradient of $\psi$ along the vertical space \citep[Section~II.4]{sakai1996riemannian}. A function $g$ of the position has $\nabla_x g = \nabla g$ and $\nabla_v g = 0$. In $\mathbb R^d$, $\nabla_x$ and $\nabla_v$ are the partial gradients in the position and in the velocity.

Under $\mu$, given the position $x$, the velocity has the standard Gaussian distribution of $T_x\mathcal X$. A field $W$ on $T\mathcal X$ assigns to each state $(x,v)$ a vector $W(x,v) \in T_x\mathcal X$. Let $\nabla_v \cdot W$ be the divergence of $W(x,\cdot)$ on $T_x\mathcal X$, and let $\nabla_v^* \cdot W := -\nabla_v \cdot W + v \cdot W$. The Gaussian integration by parts in the velocity gives, for every function $\psi$ and every field $W$ on $T\mathcal X$,
\begin{equation}
    \mu( \psi\, \nabla_v^* \cdot W ) = \mu( \nabla_v\psi \cdot W ) .
    \label{eq: adjoints}
\end{equation}
A function $\psi \in L^2(\mu)$ has a velocity gradient in $L^2(\mu)$ if there is a field $\nabla_v\psi \in L^2(\mu)$ such that \eqref{eq: adjoints} holds for every smooth field $W$ on $T\mathcal X$ with compact support.

\subsection{Underdamped Langevin dynamics}
\label{sec: dynamics}

The underdamped Langevin dynamics \eqref{eq: uld} is a process $Z_t := (X_t,V_t)$ on $T\mathcal X$ with the position $X_t \in \mathcal X$ and the velocity $V_t \in T_{X_t}\mathcal X$. The position moves with the velocity. The equation of the velocity needs a moving frame, since the tangent space $T_{X_t}\mathcal X$ changes with $t$. Let the frame $\mathfrak e_1(t), \dots, \mathfrak e_d(t)$ be an orthonormal basis of $T_{X_t}\mathcal X$ that is parallel along the path of the position, and let $(B_t^1, \dots, B_t^d)_{t \Ge 0}$ be a standard Brownian motion in $\mathbb R^d$. We define the Brownian motion $(B_t)_{t \Ge 0}$ in the tangent spaces by its components $(B_t^i)_{t \Ge 0}$ in the frame, that is, for a field $W$ on $T\mathcal X$, the stochastic integral against $(B_t)_{t \Ge 0}$ is
\begin{equation}
    \int_0^t W(Z_s) \cdot \D B_s := \sum_{i=1}^d \int_0^t \mathfrak e_i(s) \cdot W(Z_s) \,\D B_s^i .
    \label{eq: stochastic integral}
\end{equation}
We define the velocity through its components in the frame: the components $V_t^i := \mathfrak e_i(t) \cdot V_t$ solve the system
\begin{equation}
    \D V_t^i = -\bigl( \mathfrak e_i(t) \cdot \nabla U(X_t) + \gamma V_t^i \bigr) \D t + \sqrt{2\gamma}\,\D B_t^i , \qquad i = 1, \dots, d .
    \label{eq: frame}
\end{equation}
Since the frame is parallel, $\D V_t^i$ are the components of the covariant differential $\D V_t$, and $\D B_t^i$ are the components of $\D B_t$ by \eqref{eq: stochastic integral}. Hence the system \eqref{eq: frame} is the second equation of \eqref{eq: uld}. For $U = 0$, \eqref{eq: frame} is the equation of the Ornstein--Uhlenbeck process, which \citet{jorgensen1978construction} constructs in a Riemannian manifold.

For a function $\psi$ on $T\mathcal X$, It\^o's formula with \eqref{eq: frame} reads
\begin{equation}
    \psi(Z_t) = \psi(Z_0) + \int_0^t \mathcal L\psi(Z_s) \,\D s + \sqrt{2\gamma} \int_0^t \nabla_v\psi(Z_s) \cdot \D B_s ,
    \label{eq: ito}
\end{equation}
with the generator of the underdamped Langevin dynamics \eqref{eq: uld}
\begin{equation}
    \mathcal L = \mathcal A + \mathcal S , \qquad
    \mathcal A := v \cdot \nabla_x - \nabla U \cdot \nabla_v , \qquad
    \mathcal S := -\gamma v \cdot \nabla_v + \gamma \Delta_v .
    \label{eq: generator}
\end{equation}

The operator $\mathcal A$ is the generator of the Hamiltonian flow of $U(x) + \frac12 |v|^2$ on $T\mathcal X$, and $\mathcal S$ is the generator of the Ornstein--Uhlenbeck process of the velocity in $T_x\mathcal X$ at a fixed position $x$. The Hamiltonian flow preserves $U(x) + \frac12 |v|^2$ and the measure $\D z$, so $\mathcal A$ is antisymmetric in $L^2(\mu)$, that is, $\mu( \psi_1 \mathcal A\psi_2 ) + \mu( \psi_2 \mathcal A\psi_1 ) = 0$ for functions $\psi_1$ and $\psi_2$ on $T\mathcal X$. By \eqref{eq: adjoints}, $\mathcal S$ is symmetric in $L^2(\mu)$ with the Dirichlet form $\gamma\, \mu( \nabla_v\psi_1 \cdot \nabla_v\psi_2 )$. Hence the adjoint of $\mathcal L$ is $\mathcal L^* = -\mathcal A + \mathcal S$, and $\psi_1$ and $\psi_2$ satisfy
\begin{equation}
    \mu\bigl( (-\mathcal L\psi_1)\, \psi_2 \bigr) = \mu( \psi_1\, \mathcal A\psi_2 ) + \gamma\, \mu( \nabla_v\psi_1 \cdot \nabla_v\psi_2 ) .
    \label{eq: generator form}
\end{equation}

The semigroup of the underdamped Langevin dynamics \eqref{eq: uld} is $(e^{t\mathcal L}\psi)(z) = \mathbb E^z\bigl[\psi(Z_t) \bigr]$, where $\mathbb E^z$ is the expectation under the law of the dynamics started at $z$. When the lower bound \eqref{eq: curvature} and the Poincaré inequality \eqref{eq: poincare} hold, the semigroup decays in $L^2(\mu)$: there are a constant $C$ and a rate $\lambda > 0$ such that, for every $t \Ge 0$ and every $\psi \in L^2(\mu)$ with $\mu(\psi) = 0$,
\begin{equation}
    \mu\bigl( ( e^{t\mathcal L}\psi )^2 \bigr) \Le C e^{-2\lambda t}\, \mu(\psi^2) .
    \label{eq: decay}
\end{equation}
The rate $\lambda$ is explicitly determined in Theorem~\ref{thm: decay}.

\subsection{Poisson equation}
\label{sec: poisson}

As in \eqref{eq: poisson solution}, the Poisson solution for a test function $f$ can be formally written as the integral
\begin{equation*}
    \phi = \int_0^\infty \bigl( e^{t\mathcal L}f - \mu(f) \bigr) \D t ,
\end{equation*}
and the convergence of the integral in $L^2(\mu)$ is guaranteed by the decay \eqref{eq: decay}. The strong Poisson equation \eqref{eq: poisson equation} requires the solution $\phi$ to be twice differentiable in the velocity and differentiable in the position, and the decay in $L^2(\mu)$ cannot supply this regularity.\footnote{For $\mathcal X = \mathbb R^d$ and a large damping rate $\gamma$, \citet{ye2026semigroup} prove that $\phi$ is a strong solution of \eqref{eq: poisson equation} when $U$ has bounded second and third derivatives and is strongly convex outside a ball.} Hence we need to define the weak solution.

The identity \eqref{eq: generator form} with $\psi_1 = \phi$ gives the weak form of the Poisson equation \eqref{eq: poisson equation}, with test functions in the space $C^1(T\mathcal X)$ of the continuously differentiable functions on $T\mathcal X$.

\begin{definition}
\label{def: weak solution}
    A function $\phi \in L^2(\mu)$ is a weak solution to the Poisson equation \eqref{eq: poisson equation}, if $\nabla_v\phi \in L^2(\mu)$, and for every $\psi \in C^1(T\mathcal X)$ with $\psi, \mathcal A\psi, \nabla_v\psi \in L^2(\mu)$,
    \begin{equation}
        \mu( \phi\, \mathcal A\psi ) + \gamma\, \mu( \nabla_v\phi \cdot \nabla_v\psi ) = \mu\bigl( ( f - \mu(f) )\, \psi \bigr) .
        \label{eq: weak poisson}
    \end{equation}
    A function $\phi$ that is smooth on $T\mathcal X$ is a strong solution if \eqref{eq: poisson equation} holds on $T\mathcal X$.
\end{definition}

For $\psi = \phi$, the antisymmetry of $\mathcal A$ gives $\mu( \phi \mathcal A\phi ) = 0$, and \eqref{eq: weak poisson} is the identity
\begin{equation}
    \gamma\, \mu\bigl( |\nabla_v\phi|^2 \bigr) = \mu\bigl( \phi\, ( f - \mu(f) ) \bigr) .
    \label{eq: energy}
\end{equation}
The asymptotic variance of $f$ is
\begin{equation}
    \sigma_f^2 = 2\gamma\, \mu\bigl( |\nabla_v\phi|^2 \bigr) .
    \label{eq: asymptotic variance}
\end{equation}
The name \emph{asymptotic variance} comes from the variance of the limit distribution $\mathcal N(0, \sigma_f^2)$ in the central limit theorem \eqref{eq: clt convergence}. The rigorous statement on the Poisson solution is in Theorem~\ref{thm: poisson}.

\section{Notations and assumptions}
\label{sec: notation}

This section fixes the remaining notations and states the assumptions of the paper.

The spaces $C^\infty(\mathcal X)$ and $C^\infty(T\mathcal X)$ consist of the smooth functions on $\mathcal X$ and on $T\mathcal X$. The spaces $C_c^\infty(\mathcal X)$ and $C_c^\infty(T\mathcal X)$ consist of the smooth functions with compact support on $\mathcal X$ and on $T\mathcal X$. For $p \Ge 1$, $L^p(\mu)$ is the space of the functions and of the fields $\psi$ on $T\mathcal X$ with $\mu\bigl( |\psi|^p \bigr) < \infty$, with the norm
\begin{equation*}
    \norm{\psi}_{L^p(\mu)} := \mu\bigl( |\psi|^p \bigr)^{\frac1p} .
\end{equation*}
Likewise, $L^p(\bar\mu)$ is the space of the functions and of the fields $g$ on $\mathcal X$ with $\mu\bigl( |g|^p \bigr) < \infty$. Regarded as functions of the position on $T\mathcal X$, they are the functions and the fields of the position in $L^p(\mu)$, with $\norm{g}_{L^p(\bar\mu)} = \norm{g}_{L^p(\mu)}$.
The process $(Z_t)_{t \Ge 0}$ is the underdamped Langevin dynamics \eqref{eq: uld} on $T\mathcal X$ with the initial state $Z_0 \sim \mu$. 

For probability distributions $\nu$ and $\nu'$ on $\mathbb R$ with finite first moments, the Wasserstein-1 distance is
\begin{equation*}
	\mathcal W_1(\nu,\nu') := \inf \bigl\{ \mathbb E|\xi - \xi'| : \mathrm{Law}(\xi) = \nu ,\ \mathrm{Law}(\xi') = \nu' \bigr\} ,
\end{equation*}
where the infimum runs over the random variables $(\xi,\xi')$ whose marginals are $\nu$ and $\nu'$. The letter $C$ denotes a universal constant, and it changes from line to line. In particular, no constant $C$ depends on $d$, $m$, $K$, $\gamma$, the test function $f$ or the time $t$.
Table~\ref{tab: notation} lists the main notations of the paper.

\begin{table}[!htb]
\centering
\begin{tabular}{@{}lp{120mm}@{}}
\toprule
\textbf{symbol} & \textbf{definition} \\
\midrule
\multicolumn{2}{@{}l}{$\blacktriangleright$ \textit{manifold}} \\
$\mathcal X$ & Riemannian manifold of dimension $d$ \\
$T_x\mathcal X$, $T\mathcal X$ & tangent space at $x$ and tangent bundle of $\mathcal X$ \\
$z = (x,v)$ & state in $T\mathcal X$ with the position $x \in \mathcal X$ and the velocity $v \in T_x\mathcal X$ \\
$\D x$, $\D v$, $\D z$ & Riemannian volume of $\mathcal X$, Lebesgue measure of $T_x\mathcal X$ and $\D z = \D x \D v$ \\
$\Ric(\cdot)$ & Ricci curvature of $\mathcal X$ \\
\midrule
\multicolumn{2}{@{}l}{$\blacktriangleright$ \textit{potential and dynamics}} \\
$U$ & potential function on $\mathcal X$ \\
$\gamma$ & damping rate \\
$(B_t)_{t \Ge 0}$ & Brownian motion in the tangent spaces, defined in \eqref{eq: stochastic integral} \\
$Z_t = (X_t,V_t)$ & underdamped Langevin dynamics, defined in \eqref{eq: uld} \\
$\mu$, $\bar\mu$ & Boltzmann distribution \eqref{eq: boltzmann} and its marginal distribution of the position \\
$L^p(\mu)$, $L^p(\bar\mu)$ & spaces of the functions and fields on $T\mathcal X$ and on $\mathcal X$ with finite $p$-th moment under $\mu$ and under $\bar\mu$ \\
$\mathcal L$ & generator of the underdamped Langevin dynamics, defined in \eqref{eq: generator} \\
$(e^{t\mathcal L})_{t \Ge 0}$ & semigroup of the underdamped Langevin dynamics \\
$\lambda$ & rate of the decay of $(e^{t\mathcal L})_{t \Ge 0}$ in $L^2(\mu)$, defined in \eqref{eq: l2 decay} \\
$\bar{\mathcal L}$ & diffusion operator on $\mathcal X$, defined in \eqref{eq: overdamped} \\
\midrule
\multicolumn{2}{@{}l}{$\blacktriangleright$ \textit{parameters}} \\
$K$ & $-K$ is the lower bound of $\Ric + \nabla^2 U$, defined in \eqref{eq: curvature} \\
$m$ & $m^{-1}$ is the smallest constant of the Poincaré inequality, defined in \eqref{eq: poincare} \\
\midrule
\multicolumn{2}{@{}l}{$\blacktriangleright$ \textit{test function and Poisson solution}} \\
$f$ & test function in $L^2(\bar\mu)$ \\
$\phi$ & Poisson solution of $f$, defined in \eqref{eq: poisson solution} \\
$\sigma_f^2$ & asymptotic variance of $f$, defined in \eqref{eq: asymptotic variance} \\
$\mathbb V(\gamma)$ & normalized variance, defined in \eqref{eq: normalized variance} \\
\bottomrule
\end{tabular}
\caption{Main notations of the paper.}
\label{tab: notation}
\end{table}

The first assumption concerns the manifold.

\begin{assumption}[manifold]
\label{ass: manifold}
    The manifold $\mathcal X$ is a complete, connected and oriented Riemannian manifold of dimension $d \Ge 1$ without boundary.
\end{assumption}

Under Assumption~\ref{ass: manifold}, the underdamped Langevin dynamics \eqref{eq: uld} is well posed for every potential $U \in C^\infty(\mathcal X)$ that is bounded from below with $e^{-U}$ integrable and for every $\gamma > 0$. The solution exists for all times from every initial state by the Lyapunov criterion of nonexplosion \citep[Theorem~3.5]{khasminskii2012stochastic}, and it is unique in law and is a Markov process, since the coefficients of \eqref{eq: uld} are smooth. The Boltzmann distribution $\mu$ is invariant \citep[Section~4.1]{eberle2026space}. The generator of the semigroup in $L^2(\mu)$ satisfies the condition of \citep[Theorem~1.1]{hormander1967hypoelliptic}, and the smooth functions with compact support on $T\mathcal X$ form a core of this generator.

The second assumption concerns the potential.

\begin{assumption}[curvature]
\label{ass: curvature}
    Let $K \Ge 0$. The potential $U \in C^\infty(\mathcal X)$ is bounded from below and satisfies the lower bound \eqref{eq: curvature},
    \begin{equation*}
        w \cdot \bigl( \Ric + \nabla^2 U \bigr)(x)\, w \Ge -K |w|^2 , \qquad \forall x \in \mathcal X ,\quad w \in T_x\mathcal X .
    \end{equation*}
\end{assumption}

For $U = 0$, the lower bound \eqref{eq: curvature} holds with $K = 0$ on a flat torus and on a round sphere. In $\mathbb R^d$, \eqref{eq: curvature} holds with $K = 0$ for every convex potential $U$.

The third assumption concerns the distribution $\bar\mu$ of the position.

\begin{assumption}[Poincaré]
\label{ass: poincare}
    Let $m > 0$. The function $e^{-U}$ is integrable with respect to the Riemannian volume of $\mathcal X$. The Poincaré inequality \eqref{eq: poincare} holds, namely, for every $g \in C^1(\mathcal X) \cap L^2(\bar\mu)$ with $\nabla g \in L^2(\bar\mu)$,
    \begin{equation*}
        \mu\bigl( (g - \mu(g))^2 \bigr) \Le \frac1m\, \mu\bigl( |\nabla g|^2 \bigr) .
    \end{equation*}
    The constant $m^{-1}$ is the smallest constant of \eqref{eq: poincare}, and $\bar{\mathcal L}$ in \eqref{eq: overdamped} has a discrete spectrum in $L^2(\bar\mu)$.
\end{assumption}

The constant $m$ is the spectral gap of the diffusion operator $\bar{\mathcal L}$. The spectrum of $\bar{\mathcal L}$ is discrete when $\mathcal X$ is compact, and in $\mathbb R^d$ when $\frac14 |\nabla U(x)|^2 - \frac12 \Delta U(x) \to \infty$ as $|x| \to \infty$ \citep[Corollary~4.10.9]{bakry2014analysis}. The discrete spectrum and the orientation of $\mathcal X$ enter the results only through the decay of the semigroup in Theorem~\ref{thm: decay}.

Throughout the paper, Assumptions~\ref{ass: manifold}, \ref{ass: curvature} and~\ref{ass: poincare} hold, and the statements do not repeat them.

\section{Main results}
\label{sec: results}

This section states the results of the paper. Section~\ref{sec: decay} states the decay of the semigroup $(e^{t\mathcal L})_{t \Ge 0}$ in $L^2(\mu)$. Section~\ref{sec: regularity} states the regularity of the Poisson solution, bounds the asymptotic variance and bounds the mean-square error of the time average. Section~\ref{sec: variance} defines the normalized variance and states its properties. Section~\ref{sec: mixing} bounds the mixing time from below and from above by the normalized variance. Section~\ref{sec: clt} states the central limit theorem in the Wasserstein-1 distance. The proofs are in Appendices~\ref{app: decay}--\ref{app: clt proof}, and Table~\ref{tab: proofs} lists the appendix of each result.

\begin{table}[!htb]
\centering
\begin{tabular}{@{}lll@{}}
\toprule
\textbf{section} & \textbf{appendix} & \textbf{results} \\
\midrule
Section~\ref{sec: decay} & Appendix~\ref{app: decay} & Theorem~\ref{thm: decay} \\
Section~\ref{sec: regularity} & Appendix~\ref{app: poisson} & Theorem~\ref{thm: poisson}, Corollaries~\ref{cor: variance} and~\ref{cor: second moment} \\
Section~\ref{sec: variance} & Appendix~\ref{app: variance} & Theorem~\ref{thm: properties}, Lemmas~\ref{lem: quadratic} and~\ref{lem: torus} \\
Section~\ref{sec: mixing} & Appendix~\ref{app: mixing} & Theorem~\ref{thm: mixing} \\
Section~\ref{sec: clt} & Appendix~\ref{app: clt proof} & Theorem~\ref{thm: clt} \\
\bottomrule
\end{tabular}
\caption{Results of Section~\ref{sec: results} and the appendices of their proofs.}
\label{tab: proofs}
\end{table}

\subsection{Decay of the semigroup in $L^2(\mu)$}
\label{sec: decay}

We prove the decay of the semigroup $(e^{t\mathcal L})_{t \Ge 0}$ in $L^2(\mu)$ with a new rate.

\begin{figure}[!htb]
\centering
\includegraphics[width=\textwidth]{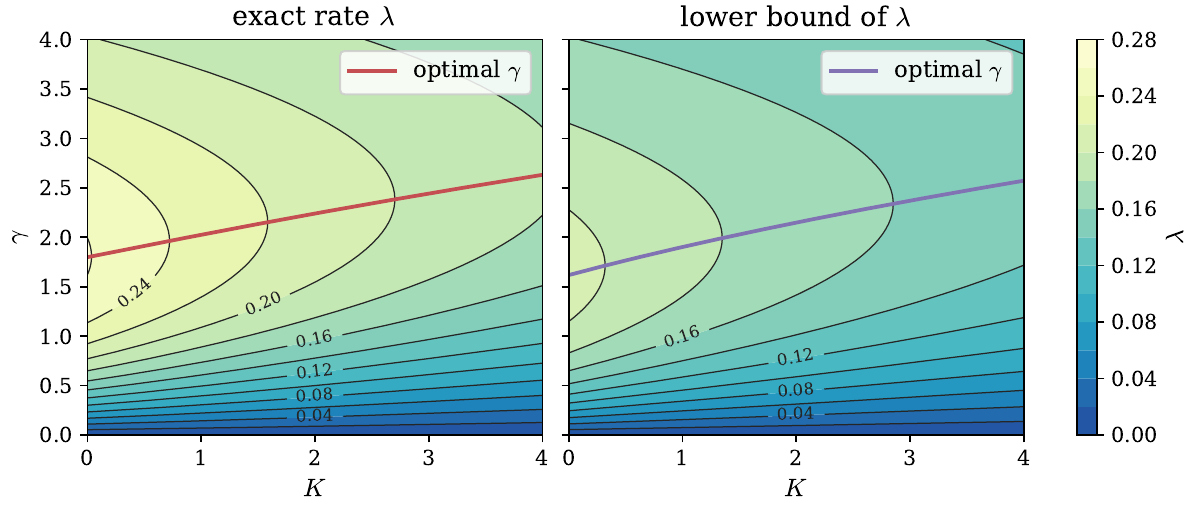}
\caption{Contours of the rate $\lambda$ of \eqref{eq: l2 decay} (left) and of the lower bound \eqref{eq: rate lower} (right) for $m = 1$, as functions of $K$ and $\gamma$, with the same levels and colors. The red curve is the damping rate $\gamma$ that maximizes $\lambda$ for each $K$, and the purple curve is the damping rate $\gamma$ that maximizes the lower bound \eqref{eq: rate lower} for each $K$.}
\label{fig: rate}
\end{figure}

\begin{theorem}
\label{thm: decay}
    Let $\gamma > 0$, and let $\beta = \beta(\gamma,m,K)$ be the unique positive root of
    \begin{equation}
        \beta^2 \Bigl( 1 + \frac Km + \Bigl( \frac{\gamma}{\sqrt m} + \beta \Bigr)^2 \Bigr) = 1 .
        \label{eq: root}
    \end{equation}
    Then every $\psi \in L^2(\mu)$ with $\mu(\psi) = 0$ and every $t \Ge 0$ satisfy
    \begin{equation}
        \mu\bigl( ( e^{t\mathcal L}\psi )^2 \bigr) \Le e^{1 - 2\lambda t}\, \mu(\psi^2) , \qquad
        \lambda := \frac{\gamma \beta^2}{1 + \beta^2} .
        \label{eq: l2 decay}
    \end{equation}
    Consequently, the rate $\lambda$ satisfies
    \begin{equation}
        \lambda \Ge \frac{m\gamma}{\gamma^2 + 2 \sqrt{\frac{\sqrt5 - 1}{2}} \gamma\sqrt m + \frac{3 + \sqrt5}{2} m + K} .
        \label{eq: rate lower}
    \end{equation}
\end{theorem}
The decay \eqref{eq: l2 decay} is equivalently written as, for every $\psi \in L^2(\mu)$ with $\mu(\psi) = 0$,
\begin{equation*}
    \norm{e^{t\mathcal L}\psi}_{L^2(\mu)} \Le e^{\frac12 - \lambda t} \norm{\psi}_{L^2(\mu)} , \qquad \forall t\Ge0.
\end{equation*}
Figure~\ref{fig: rate} shows the rate $\lambda$ of \eqref{eq: l2 decay} and the lower bound \eqref{eq: rate lower} for $m = 1$ as functions of $K$ and $\gamma$, with the damping rate $\gamma$ that maximizes each of them for every $K$.

The lower bound \eqref{eq: rate lower} has the form of the rates known in $\mathbb R^d$, where $\Ric = 0$. \citet{cao2023explicit} prove the decay at the rate \eqref{eq: euclidean rate} with a universal constant $C$. \citet{fan2026sharp} prove the decay \eqref{eq: decay} with the rate
\begin{equation*}
    \lambda = \frac{m\gamma}{\gamma^2 + 2 \gamma \bigl( 2\sqrt m + \sqrt{2m + 2K} \bigr)  + 6m + 2K} \Ge \frac{m\gamma}{6 ( \gamma + \sqrt m + \sqrt K )^2} ,
\end{equation*}
by the modified $L^2$ method of \citet{dolbeault2015hypocoercivity}, and \citet{du2026sharp} extend this method to bounded convex domains of $\mathbb R^d$ with specular reflection of the velocity. On Riemannian manifolds with a locally convex boundary, \citet[Theorem~16]{eberle2026space} prove the decay of the averages of $\norm{e^{s\mathcal L}\psi}_{L^2(\mu)}$ over time windows, under a discrete spectrum of $\bar{\mathcal L}$. For every $\gamma > 0$, the lower bound \eqref{eq: rate lower} is larger than the rates of \citet{fan2026sharp} and of \citet{eberle2026space}.

The proof of Theorem~\ref{thm: decay} in Appendix~\ref{app: decay} follows \citep[Theorem~16]{eberle2026space}: a divergence equation is solved in the eigenfunctions of $\bar{\mathcal L}$, and its solution gives a test function for the equation $\partial_s ( e^{s\mathcal L}\psi ) = \mathcal L ( e^{s\mathcal L}\psi )$.

\begin{definition}
\label{def: mixing time}
    The mixing time is
    \begin{equation}
        t_{\mathrm{mix}}(\gamma) := \inf \Bigl\{ t \Ge 0 : \norm{e^{t\mathcal L}\psi}_{L^2(\mu)} \Le e^{-1} \norm{\psi}_{L^2(\mu)} \text{ for every } \psi \in L^2(\mu) \text{ with } \mu(\psi) = 0 \Bigr\} .
        \label{eq: mixing time}
    \end{equation}
\end{definition}
The decay \eqref{eq: l2 decay} and the lower bound \eqref{eq: rate lower} of $\lambda$ give an upper bound of the mixing time: for every $\gamma > 0$,
\begin{equation}
    t_{\mathrm{mix}}(\gamma) \Le \frac{3}{2\lambda} \Le \frac{3 \bigl( \gamma^2 + 2 \sqrt{\frac{\sqrt5 - 1}{2}} \gamma\sqrt m + \frac{3 + \sqrt5}{2} m + K \bigr)}{2 m\gamma} .
    \label{eq: mixing decay}
\end{equation}
In $\mathbb R^d$, when the Hessian of $U$ is bounded from below by $-K$, \citet[Theorem~27]{eberle2026relaxation} prove $t_{\mathrm{mix}}(\gamma) \Le C \max\bigl( \frac\gamma m , \frac{m + K}{m\gamma} \bigr)$, which has the form of the upper bound \eqref{eq: mixing decay}.

\subsection{Regularity of the Poisson solution}
\label{sec: regularity}

With a given test function $f$, we show the regularity of the Poisson solution \eqref{eq: poisson solution} and bound the asymptotic variance \eqref{eq: asymptotic variance} from above in Theorem~\ref{thm: poisson}. Corollary~\ref{cor: variance} complements these bounds with a bound of the Poisson solution in $L^4(\mu)$ and a lower bound of the asymptotic variance.

\begin{theorem}
\label{thm: poisson}
    Let $\gamma > 0$ and $f \in L^2(\bar\mu)$. Then the integral \eqref{eq: poisson solution} converges in $L^2(\mu)$, and $\phi$ is a weak solution of \eqref{eq: poisson equation}. If $f \in C^\infty(\mathcal X)$, then $\phi$ is a strong solution. For the asymptotic variance $\sigma_f^2 = 2\gamma\, \mu\bigl( |\nabla_v\phi|^2 \bigr)$,
    \begin{equation}
        \sigma_f^2 = 2\gamma \norm{\nabla_v \phi}_{L^2(\mu)}^2 \Le \frac{2 ( \gamma^2 + m + K )}{m\gamma}\, \mu\bigl( ( f - \mu(f) )^2 \bigr) ,
        \label{eq: variance upper}
    \end{equation}
    and the $L^2(\mu)$ norm of $\phi$ is bounded by
    \begin{equation}
        \norm{\phi}_{L^2(\mu)}^2 \Le \frac{\gamma^2 + 2m + K}{2m\gamma}\, \sigma_f^2 .
        \label{eq: poisson moments}
    \end{equation}
    For every $t \Ge 0$, It\^o's formula \eqref{eq: poisson decomposition} holds almost surely along the dynamics $(Z_t)_{t \Ge 0}$ from $\mu$-almost every initial state, and simultaneously for all $t \Ge 0$ after passing to a continuous modification of $t \mapsto \phi(Z_t)$,
    \begin{equation}
        \int_0^t \bigl( f(X_s) - \mu(f) \bigr) \D s = \phi(Z_0) - \phi(Z_t) + \sqrt{2\gamma} \int_0^t \nabla_v\phi(Z_s) \cdot \D B_s .
        \label{eq: martingale decomposition}
    \end{equation}
\end{theorem}
The bounds \eqref{eq: variance upper} and \eqref{eq: poisson moments} bound the $L^2(\mu)$ norms of $\nabla_v\phi$ and $\phi$ by the variance of $f$ under $\mu$. Along the underdamped Langevin dynamics \eqref{eq: uld} with a stationary start, define the martingale
\begin{equation}
    N_t := \sqrt{2\gamma} \int_0^t \nabla_v\phi(Z_s) \cdot \D B_s ,
    \label{eq: martingale}
\end{equation}
where the stochastic integral against $(B_t)_{t \Ge 0}$ is \eqref{eq: stochastic integral} for the field $\nabla_v\phi$, and write It\^o's formula \eqref{eq: martingale decomposition} as
\begin{equation}
    \frac1t \int_0^t f(X_s) \,\D s - \mu(f) = \frac{\phi(Z_0) - \phi(Z_t)}{t} + \frac{N_t}{t} .
    \label{eq: time average}
\end{equation}
By \eqref{eq: asymptotic variance}, $\mathbb E[N_t^2] = \sigma_f^2\, t$, so as $t \to \infty$, formally, the right side of \eqref{eq: time average} is dominated by $N_t/t$, whose second moment is $\sigma_f^2/t$. Hence the asymptotic variance $\sigma_f^2$ measures the error of the time average.

\begin{corollary}
\label{cor: variance}
    Let $\gamma > 0$ and $f \in L^2(\bar\mu)$, and let $\phi$ be the Poisson solution \eqref{eq: poisson solution} of $f$.
    \begin{enumerate}
        \renewcommand{\labelenumi}{(\roman{enumi})}
        \item If $f \in L^4(\bar\mu)$, then
        \begin{equation}
            \norm{\phi}_{L^4(\mu)} \Le \frac{12 ( \gamma + \sqrt m + \sqrt K )^2}{m\gamma}\, \norm{f - \mu(f)}_{L^4(\mu)} .
            \label{eq: poisson fourth}
        \end{equation}
        \item If $\nabla f \in L^2(\bar\mu)$ and $\mu\bigl( |\nabla f|^2 \bigr) > 0$, then
        \begin{equation}
            \sigma_f^2 \Ge 2\gamma\, \frac{\mu\bigl( ( f - \mu(f) )^2 \bigr)^2}{\mu\bigl( |\nabla f|^2 \bigr)} .
            \label{eq: variance below}
        \end{equation}
    \end{enumerate}
\end{corollary}
Since $m^{-1}$ is the smallest constant of \eqref{eq: poincare}, the bounds \eqref{eq: variance upper} and \eqref{eq: variance below} give
\begin{equation}
    \frac{2\gamma}{m^2} \Le \sup_{f \in C^1(\mathcal X)} \frac{\sigma_f^2}{\mu\bigl( |\nabla f|^2 \bigr)} \Le \frac{2 ( \gamma^2 + m + K )}{m^2\gamma} .
    \label{eq: variance}
\end{equation}

\begin{remark}
\label{rem: overdamped}
    The bounds \eqref{eq: variance upper} and \eqref{eq: variance below} have an alternate form in terms of the diffusion operator $\bar{\mathcal L}$ in \eqref{eq: overdamped}. For $f \in L^2(\bar\mu)$, let $\bar\phi$ be the solution of $-\bar{\mathcal L}\bar\phi = f - \mu(f)$ on $\mathcal X$ with $\mu(\bar\phi) = 0$. Then
    \begin{equation}
        2\gamma\, \mu\bigl( |\nabla\bar\phi|^2 \bigr) \Le \sigma_f^2 \Le 2\gamma\, \mu\bigl( |\nabla\bar\phi|^2 \bigr) + \frac2\gamma \mu\bigl( |\nabla^2\bar\phi|^2 \bigr).
        \label{eq: variance overdamped}
    \end{equation}
    In $\mathbb R^d$, \citet{dechant2023thermodynamic} derive the first inequality of \eqref{eq: variance overdamped} and the bound \eqref{eq: variance below} formally, as lower bounds of the correlation time of $f$. For nonreversible perturbations of the overdamped Langevin dynamics, the asymptotic variance is at most the asymptotic variance of the reversible dynamics \citep{duncan2016variance}.
\end{remark}

We next bound the mean-square error of the time average. By \eqref{eq: poisson moments}, $\phi(Z_0) - \phi(Z_t)$ of \eqref{eq: time average} satisfies
\begin{equation*}
    \mathbb E\bigl[ ( \phi(Z_0) - \phi(Z_t) )^2 \bigr] \Le 4 \norm{\phi}_{L^2(\mu)}^2 \Le \frac{2 ( \gamma^2 + 2m + K )}{m\gamma}\, \sigma_f^2 .
\end{equation*}
Consequently,
\begin{equation*}
    t\, \mathbb E\biggl[ \Bigl( \frac1t \int_0^t f(X_s) \,\D s - \mu(f) \Bigr)^2 \biggr] \to \sigma_f^2 \qquad \text{as } t \to \infty .
\end{equation*}
Corollary~\ref{cor: second moment} gives the explicit rate of the convergence.

\begin{corollary}
\label{cor: second moment}
    Let $\gamma > 0$ and $f \in L^2(\bar\mu)$. Then
    \begin{equation}
        \Bigl( \frac1t - \frac{2 ( \gamma^2 + 2m + K )}{m\gamma t^2} \Bigr) \sigma_f^2 \Le \mathbb E\biggl[ \Bigl( \frac1t \int_0^t f(X_s) \,\D s - \mu(f) \Bigr)^2 \biggr] \Le \Bigl( \frac1t + \frac{2}{\gamma t^2} \Bigr) \sigma_f^2 , \qquad \forall t > 0 .
        \label{eq: second moment}
    \end{equation}
\end{corollary}
The asymptotic variance $\sigma_f^2$ governs the convergence of the time averages in mean square, and the upper bound of \eqref{eq: second moment} has no explicit dependence on $K$.

For general stationary Markov processes, \citet{xu2026singular} bounds the difference between the variance of the time average of a function $\psi$ and its leading term by a multiple of $\mu(\psi^2) t^{-2}$ divided by the square of a gap of the generator. \citet{paulin2015concentration} proves a bound of the same form in discrete time. In \eqref{eq: second moment}, the difference between the mean-square error and $\sigma_f^2 / t$ is bounded by a multiple of $\sigma_f^2 t^{-2}$. For bounded functions, \citet{birrell2019concentration} prove concentration inequalities for the time averages under a hypocoercive estimate, with the variance of the function divided by the hypocoercive rate in place of $\sigma_f^2$.

\subsection{Properties of the normalized variance}
\label{sec: variance}

The bound \eqref{eq: variance upper} controls the ratio of the asymptotic variance $\sigma_f^2$ to the variance $\mu\bigl( ( f - \mu(f) )^2 \bigr)$ of the test function. We define the supremum of this ratio as the normalized variance $\mathbb V(\gamma)$, state its properties as a function of the damping rate $\gamma$, and give its exact value for the quadratic potential and for the flat torus.

\begin{definition}
\label{def: normalized variance}
    The normalized variance is
    \begin{equation}
        \mathbb V(\gamma) := \sup_{f \in L^2(\bar\mu)} \frac{\sigma_f^2}{\mu\bigl( ( f - \mu(f) )^2 \bigr)} = \sup_{f \in L^2(\bar\mu)} \frac{2\gamma\, \mu\bigl( |\nabla_v\phi|^2 \bigr)}{\mu\bigl( ( f - \mu(f) )^2 \bigr)} ,
        \label{eq: normalized variance}
    \end{equation}
    where the supremum is over the nonconstant functions $f \in L^2(\bar\mu)$, and $\phi$ is the Poisson solution \eqref{eq: poisson solution} of $f$.
\end{definition}

In terms of $\mathbb V(\gamma)$, \eqref{eq: second moment} implies
\begin{equation}
	\sup_{\norm{f - \mu(f)}_{L^2(\bar\mu)} \Le 1}
    \mathbb E\biggl[ \Bigl( \frac1t \int_0^t f(X_s) \,\D s - \mu(f) \Bigr)^2 \biggr] \Le \Bigl( \frac1t + \frac{2}{\gamma t^2} \Bigr) \mathbb V(\gamma) , \qquad \forall t > 0 .
    \label{eq: mean square}
\end{equation}
Since $m^{-1}$ is the smallest constant of \eqref{eq: poincare}, the bounds \eqref{eq: variance upper} and \eqref{eq: variance below} give
\begin{equation}
    \frac{2\gamma}{m} \Le \mathbb V(\gamma) \Le \frac{2 ( \gamma^2 + m + K )}{m\gamma} .
    \label{eq: normalized bounds}
\end{equation}
By \eqref{eq: variance below}, $\sigma_f^2 > 0$ for every nonconstant $f$ with $\nabla f \in L^2(\bar\mu)$. The upper bound of \eqref{eq: normalized bounds} is smallest at $\gamma = \sqrt{m + K}$. For $K = 0$, \eqref{eq: normalized bounds} reads
\begin{equation}
    \frac{2\gamma}{m} \Le \mathbb V(\gamma) \Le \frac{2 ( \gamma^2 + m )}{m\gamma} ,
    \label{eq: variance convex}
\end{equation}
and by choosing $\gamma = \sqrt m$,
\begin{equation*}
    \frac{2}{\sqrt m} \Le \mathbb V(\sqrt m) \Le \frac{4}{\sqrt m} .
\end{equation*}

Theorem~\ref{thm: properties} states two basic properties of the normalized variance as a function of $\gamma$.

\begin{theorem}
\label{thm: properties}
    The normalized variance $\mathbb V(\gamma)$ has the following properties.
    \begin{enumerate}
        \renewcommand{\labelenumi}{(\roman{enumi})}
        \item The function $\gamma \mapsto \gamma \mathbb V(\gamma)$ is nondecreasing and the function $\gamma \mapsto \gamma^{-1} \mathbb V(\gamma)$ is nonincreasing on $(0,\infty)$.
        \item If $U$ is not constant, there is a constant $c > 0$ such that $\mathbb V(\gamma) \Ge \frac c\gamma$ for every $\gamma > 0$.
    \end{enumerate}
\end{theorem}

Lemmas~\ref{lem: quadratic} and~\ref{lem: torus} give the exact value of $\mathbb V(\gamma)$ for the quadratic potential and for the flat torus.

\begin{lemma}
\label{lem: quadratic}
    Let $\mathcal X = \mathbb R^d$ and $U(x) = \frac m2 |x|^2$ with $m > 0$, so that $m^{-1}$ is the smallest constant of \eqref{eq: poincare}. Then
    \begin{equation}
        \mathbb V(\gamma) = \max\Bigl( \frac{2\gamma}{m} , \frac\gamma m + \frac1\gamma \Bigr) .
        \label{eq: variance quadratic}
    \end{equation}
    In particular, $\mathbb V(\gamma)$ is smallest at the damping rate $\gamma^*$ with
    \begin{equation}
        \gamma^* = \sqrt m , \qquad
        \mathbb V(\gamma^*) = \frac{2}{\sqrt m} .
        \label{eq: optimal damping}
    \end{equation}
\end{lemma}

\begin{lemma}
\label{lem: torus}
    Let $\mathcal X = \bigl( \mathbb R / \frac{2\pi}{\sqrt m} \mathbb Z \bigr)^d$ and $U = 0$ with $m > 0$, so that $m^{-1}$ is the smallest constant of \eqref{eq: poincare}. Then
    \begin{equation}
        \mathbb V(\gamma) = \frac2\gamma \int_0^\infty \exp\Bigl( -\frac{m}{\gamma^2} \bigl( s - 1 + e^{-s} \bigr) \Bigr) \D s ,
        \label{eq: torus formula}
    \end{equation}
    and thus
    \begin{equation}
        \max\Bigl( \frac{2\gamma}{m} , \sqrt{\frac{2\pi}{m}} \Bigr) \Le \mathbb V(\gamma) \Le \frac{2\gamma}{m} + \sqrt{\frac{2\pi}{m}} .
        \label{eq: variance torus}
    \end{equation}
\end{lemma}

By \eqref{eq: variance quadratic}, the quadratic potential has $\mathbb V(\gamma) \Ge \frac{\gamma^2 + m}{m\gamma}$, which is half of the upper bound of \eqref{eq: variance convex}, so that the upper bound of \eqref{eq: variance convex} is sharp in $m$ and $\gamma$.

For the quadratic potential of Lemma~\ref{lem: quadratic}, \citet[Section~5]{skeel2021choice} compute the autocorrelation time of the Hermite polynomials of a coordinate of the position, which is the asymptotic variance $\sigma_f^2$ divided by $2 \mu\bigl( ( f - \mu(f) )^2 \bigr)$. This time is $\frac\gamma m$ for the polynomial of degree $1$ and $\frac{\gamma}{2m} + \frac{1}{2\gamma}$ for the polynomial of degree $2$. Twice these times are the two expressions of \eqref{eq: variance quadratic}.

\subsection{Relation to the mixing time}
\label{sec: mixing}

The decay \eqref{eq: l2 decay} of the semigroup provides the upper bound \eqref{eq: mixing decay} of the mixing time \eqref{eq: mixing time} in terms of $m$ and $K$. Theorem~\ref{thm: mixing} bounds the mixing time from below and from above in terms of the normalized variance $\mathbb V(\gamma)$ and the damping rate $\gamma$ alone, with no explicit dependence on $m$ and $K$.

\begin{theorem}
\label{thm: mixing}
    Let $\gamma > 0$. Then the mixing time $t_{\mathrm{mix}}(\gamma)$ of \eqref{eq: mixing time} satisfies
    \begin{equation}
        \frac{1 - e^{-1}}{2} \max\Bigl( \mathbb V(\gamma) , \frac1\gamma \Bigr) \Le t_{\mathrm{mix}}(\gamma) \Le \Bigl( 1 + \frac\pi2 \Bigr) \Bigl( \mathbb V(\gamma) + \frac3\gamma \Bigr) .
        \label{eq: mixing bounds}
    \end{equation}
\end{theorem}

The lower bound of \eqref{eq: mixing bounds} is related to the singular value gap of a Markov semigroup with the generator $\mathcal L$ and the invariant distribution $\mu$, which is the infimum of $\norm{\mathcal L\varrho}_{L^2(\mu)} / \norm{\varrho}_{L^2(\mu)}$ over the functions $\varrho \neq 0$ in the domain of $\mathcal L$ with $\mu(\varrho) = 0$. \citet{chatterjee2025spectral} introduced this gap for finite Markov chains. Up to universal constants, the inverse of the singular value gap bounds the mixing time from below \citep[Theorem~7]{eberle2026relaxation} and the asymptotic variance of a function of unit variance from above \citep{xu2026singular}. Together, these two bounds give a lower bound of the mixing time by $\mathbb V(\gamma)$, which the proof of Theorem~\ref{thm: mixing} obtains directly from the definition \eqref{eq: mixing time}.

With the first inequality of \eqref{eq: normalized bounds}, the lower bound of \eqref{eq: mixing bounds} gives
\begin{equation}
    t_{\mathrm{mix}}(\gamma) \Ge ( 1 - e^{-1} ) \max\Bigl( \frac\gamma m , \frac{1}{2\gamma} \Bigr) .
    \label{eq: mixing lower}
\end{equation}
We compare \eqref{eq: mixing lower} with the bounds of the mixing time in the literature. For the underdamped Langevin dynamics in $\mathbb R^d$, \citet[Theorem~26]{eberle2026relaxation} prove lower bounds of the mixing time of the orders $\gamma^{-1}$ and $\gamma$. The bound \eqref{eq: mixing lower} gives both orders on Riemannian manifolds, with the Poincaré constant $m$.

For $K = 0$, \eqref{eq: mixing decay} and \eqref{eq: mixing lower} give
\begin{equation*}
    ( 1 - e^{-1} ) \max\Bigl( \frac\gamma m , \frac{1}{2\gamma} \Bigr) \Le t_{\mathrm{mix}}(\gamma) \Le 5 \Bigl( \frac\gamma m + \frac1\gamma \Bigr) ,
\end{equation*}
so that the mixing time is determined up to a universal constant for every damping rate $\gamma$. For an interval with $U = 0$ and specular reflection of the velocity, \citet[Proposition~13]{du2026sharp} bound the exponential decay rate of the semigroup from above by $\min\bigl( \frac m\gamma , 2\gamma \bigr)$.

The upper bound of \eqref{eq: mixing bounds} is new. The normalized variance $\mathbb V(\gamma)$ and the mixing time $t_{\mathrm{mix}}(\gamma)$ control rates of two different kinds: $\mathbb V(\gamma)$ controls the convergence of the time averages of the test functions in $L^2(\bar\mu)$, as in \eqref{eq: mean square}, and $t_{\mathrm{mix}}(\gamma)$ controls the contraction rate of the semigroup $(e^{t\mathcal L})_{t \Ge 0}$ in $L^2(\mu)$.
By \eqref{eq: mixing bounds}, the mixing time is determined by $\mathbb V(\gamma) + \gamma^{-1}$ up to a universal constant. The measurement in the time averages is looser than the mixing time, since it involves only the functions of the position, and $\gamma$ controls the damping in the velocity only. This also explains why the mixing time is determined by $\mathbb V(\gamma) + \gamma^{-1}$ and not by $\mathbb V(\gamma)$ alone. The two measurements are nevertheless equivalent in the sense of \eqref{eq: mixing bounds}.

\subsection{Central limit theorem}
\label{sec: clt}

We state the central limit theorem in the Wasserstein-1 distance. 
For notational convenience, define
\begin{equation}
	\alpha := \frac{\gamma + \sqrt m + \sqrt K}{\gamma} .
	\label{eq: ratio}
\end{equation}
For $f \in C^1(\mathcal X)$ with $\norm{\nabla f}_\infty < \infty$ and $t > 0$, let 
\begin{equation}
    \mathbb W(f,t) := \mathcal W_1\biggl( \mathrm{Law}\Bigl( \frac{1}{\sqrt t} \int_0^t \bigl( f(X_s) - \mu(f) \bigr) \D s \Bigr) , \mathcal N(0, \sigma_f^2) \biggr) .
    \label{eq: wasserstein}
\end{equation}
Theorem~\ref{thm: clt} bounds $\mathbb W(f,t)$ at the rate $t^{-\frac12}$ along the stationary dynamics, explicitly in $\gamma$, $m$ and $K$.

\begin{theorem}
\label{thm: clt}
    Let $\gamma > 0$ and $f \in C^1(\mathcal X)$ with $\norm{\nabla f}_\infty < \infty$. Then with $\alpha = \frac{\gamma + \sqrt m + \sqrt K}{\gamma}$,
    \begin{equation}
        \mathbb W(f,t) \Le \frac{C \gamma \alpha^2}{m^{\frac32}}\, ( 1 + \log\alpha )\, \frac{\norm{\nabla f}_\infty}{\sqrt t} , \qquad \forall t>0 .
        \label{eq: clt}
    \end{equation}
\end{theorem}
For $K = 0$, $\alpha = \frac{\gamma + \sqrt m}{\gamma}$ and \eqref{eq: clt} reads
\begin{equation}
    \mathbb W(f,t) \Le \frac{C ( \gamma + \sqrt m )^2}{m^{\frac32}\gamma} \Bigl( 1 + \log\frac{\gamma + \sqrt m}{\gamma} \Bigr) \frac{\norm{\nabla f}_\infty}{\sqrt t} ,
    \label{eq: clt convex}
\end{equation}
and by choosing $\gamma = \sqrt m$,
\begin{equation*}
    \mathbb W(f,t) \Le \frac{C}{m}\, \frac{\norm{\nabla f}_\infty}{\sqrt t} .
\end{equation*}

We compare Theorem~\ref{thm: clt} with related rates of convergence. For a continuous martingale, the Wasserstein-1 distance to the Gaussian distribution is at most twice the square root of the $L^1$ distance between the normalized quadratic variation and $1$ \citep[Theorem~3.2]{mourrat2014kantorovich}. For $(N_t)_{t \Ge 0}$ of \eqref{eq: martingale}, the normalized quadratic variation is the time average of $2\gamma |\nabla_v\phi|^2 / \sigma_f^2$. When this average has positive asymptotic variance, its $L^1$ distance to $1$ is of order $t^{-\frac12}$, giving the rate $t^{-\frac14}$. For a kinetic diffusion with a periodic force, \citet{hairer2004periodic} prove that the position at the time scale $\varepsilon^{-2}$, multiplied by $\varepsilon$, converges to a Brownian motion at a rate $\varepsilon^r$ with $r < \frac12$, in a weak Wasserstein distance on the space of paths. \citet{kusuoka2000malliavin} obtain asymptotic expansions of the distribution of additive functionals of mixing diffusions by the Malliavin calculus, and \citet{fukasawa2008edgeworth} proves an Edgeworth expansion for ergodic diffusions. In discrete time, \citet{shen2025high} prove a rate in the Wasserstein-1 distance for the sums of the Euler--Maruyama chain of the overdamped Langevin dynamics with a strongly convex potential, and \citet{fang2026martingale} prove rates close to $n^{-\frac12}$ for martingales with $n$ steps under a lower bound of the conditional covariances. Theorem~\ref{thm: clt} proves the rate $t^{-\frac12}$ in continuous time for the underdamped Langevin dynamics, and assumes no lower bound of $\sigma_f^2$.

\section{Conclusion}
\label{sec: conclusion}

This paper proves quantitative forms of the central limit theorem for the time averages of the underdamped Langevin dynamics \eqref{eq: uld} on the tangent bundle of a Riemannian manifold, for test functions of the position with bounded gradient and every damping rate $\gamma > 0$. Theorem~\ref{thm: decay} proves the decay of the semigroup in $L^2(\mu)$ with a new rate. Theorem~\ref{thm: poisson} proves that the Poisson solution is a weak solution of the Poisson equation, with an explicit bound in $L^2(\mu)$, and bounds the asymptotic variance from above. Corollary~\ref{cor: second moment} bounds the mean-square error of the time average by the asymptotic variance with no explicit dependence on $K$. Theorem~\ref{thm: clt} proves the central limit theorem at the rate $t^{-\frac12}$ in the Wasserstein-1 distance, with a coefficient that depends neither on the dimension nor on the asymptotic variance. For the quadratic potential, Lemma~\ref{lem: quadratic} shows that the upper bound of \eqref{eq: variance convex} is sharp in $m$ and $\gamma$. By Theorem~\ref{thm: mixing}, the mixing time is determined by $\mathbb V(\gamma) + \frac1\gamma$ up to a universal constant.

Several directions are left for future work. The results remain to be extended to Riemannian manifolds with boundary, where the velocity is reflected at the boundary. Another direction is to select the optimal damping rate $\gamma$ by minimizing the asymptotic variance $\sigma_f^2$ or the normalized variance $\mathbb V(\gamma)$ in \eqref{eq: normalized variance}.

\subsection*{Declaration of AI use}
The authors have used Claude Fable 5.2 and GPT-6.2 Astra to check the correctness of the proofs.

\bibliographystyle{plainnat}
\providecommand{\natexlab}[1]{}
\bibliography{references}

\appendix

\section{Proofs for the decay of the semigroup}
\label{app: decay}

This appendix proves Theorem~\ref{thm: decay}. Throughout the appendices, $\gamma > 0$ is fixed, $\mathrm{Dom}(\mathcal L)$ is the domain of the generator of the semigroup $(e^{t\mathcal L})_{t \Ge 0}$ in $L^2(\mu)$, and $\mathcal L\varrho$ denotes the generator applied to $\varrho \in \mathrm{Dom}(\mathcal L)$.

\subsection{Well-posedness and generator of the dynamics}
\label{app: generator}

First, we prove that the underdamped Langevin dynamics \eqref{eq: uld} is well posed on the tangent bundle $T\mathcal X$, and we identify the generator of its semigroup on $C_c^\infty(T\mathcal X)$ as \eqref{eq: generator},
\begin{equation*}
    \mathcal L = \mathcal A + \mathcal S = v \cdot \nabla_x - \nabla U \cdot \nabla_v - \gamma v \cdot \nabla_v + \gamma \Delta_v .
\end{equation*}
For $\psi_1 \in C_c^\infty(T\mathcal X)$, $\psi_2 \in C^1(T\mathcal X)$ and a continuous field $W$ on $T\mathcal X$ that is continuously differentiable in the velocity, \eqref{eq: adjoints} holds for $\psi_1$ and $W$, and \eqref{eq: generator form} holds for $\psi_1$ and $\psi_2$, since $\psi_1$ has compact support.

\begin{lemma}
\label{lem: dynamics}
    The solution of \eqref{eq: uld} exists for all times from every initial state, is unique in law and is a Markov process with the invariant distribution $\mu$. Its semigroup $(e^{t\mathcal L})_{t \Ge 0}$ is strongly continuous in $L^p(\mu)$ for every $p \in [1,\infty)$ and contracts $L^p(\mu)$ for every $p \Ge 1$, that is, $\norm{e^{t\mathcal L}\psi}_{L^p(\mu)} \Le \norm{\psi}_{L^p(\mu)}$ for every $\psi \in L^p(\mu)$ and $t \Ge 0$. Every $\varrho \in C_c^\infty(T\mathcal X)$ belongs to $\mathrm{Dom}(\mathcal L)$, and the generator maps $\varrho$ to the function $\mathcal L\varrho$ of \eqref{eq: generator}. Moreover, $C_c^\infty(T\mathcal X)$ is a core of the generator, that is, for every $\varrho \in \mathrm{Dom}(\mathcal L)$ there are $\varrho_n \in C_c^\infty(T\mathcal X)$ such that, as $n \to \infty$, $\varrho_n \to \varrho$ and $\mathcal L\varrho_n \to \mathcal L\varrho$ in $L^2(\mu)$.
\end{lemma}

\begin{proof}[Proof of Lemma~\ref{lem: dynamics}]
    The coefficients of \eqref{eq: uld} are smooth, so the solution is unique in law and is a Markov process until its explosion. Since $\mathcal X$ is complete, there is a smooth function $\tilde\rho \Ge 0$ on $\mathcal X$ with $\norm{\nabla\tilde\rho}_\infty \Le 1$ and $\tilde\rho(x) \to \infty$ as $x \to \infty$. Since $U$ is bounded from below, the function $V(x,v) := U(x) - \inf U + \frac12 |v|^2 + \tilde\rho(x)$ is nonnegative and proper on $T\mathcal X$ and satisfies
    \begin{equation*}
        \mathcal L V = \gamma ( d - |v|^2 ) + v \cdot \nabla\tilde\rho \Le \gamma d + \frac12 + V .
    \end{equation*}
    Hence the solution exists for all times by the Lyapunov criterion of nonexplosion \citep[Theorem~3.5]{khasminskii2012stochastic}. The Boltzmann distribution $\mu$ is invariant \citep[Section~4.1]{eberle2026space}. Hence the semigroup contracts $L^p(\mu)$ for every $p \Ge 1$. It is strongly continuous in $L^p(\mu)$ for every $p \in [1,\infty)$, since the paths of the dynamics are continuous and $C_c^\infty(T\mathcal X)$ is dense in $L^p(\mu)$. For $\varrho \in C_c^\infty(T\mathcal X)$ and $z \in T\mathcal X$, taking $\mathbb E^z$ in It\^o's formula \eqref{eq: ito} gives
    \begin{equation*}
        (e^{t\mathcal L}\varrho)(z) - \varrho(z) = \int_0^t (e^{s\mathcal L}\mathcal L\varrho)(z) \,\D s .
    \end{equation*}
    By the strong continuity of the semigroup, $\varrho \in \mathrm{Dom}(\mathcal L)$, and the generator maps $\varrho$ to $\mathcal L\varrho$ of \eqref{eq: generator}.

    Finally, to prove that $C_c^\infty(T\mathcal X)$ is a core, we first prove that the functions $\varrho - \mathcal L\varrho$, $\varrho \in C_c^\infty(T\mathcal X)$, are dense in $L^2(\mu)$. For this density, we only need to prove that $\psi = 0$ for every $\psi \in L^2(\mu)$ with
    \begin{equation}
        \mu\bigl( \psi\, ( \varrho - \mathcal L\varrho ) \bigr) = 0 , \qquad \forall \varrho \in C_c^\infty(T\mathcal X) .
        \label{eq: orthogonal}
    \end{equation}
    Since the adjoint of $\mathcal L$ in $L^2(\mu)$ is $\mathcal L^* = -\mathcal A + \mathcal S$, such a function $\psi$ solves the equation
    \begin{equation*}
        \psi + \mathcal A\psi - \mathcal S\psi = 0
    \end{equation*}
    in the sense of distributions. The operator $1 + \mathcal A - \mathcal S$ satisfies the condition of \citep[Theorem~1.1]{hormander1967hypoelliptic}, so $\psi \in C^\infty(T\mathcal X)$. Since $\mathcal X$ is complete, there are $\eta_n \in C_c^\infty(\mathcal X)$, $n \Ge 1$, with $0 \Le \eta_n \Le 1$, $\norm{\nabla\eta_n}_\infty \Le n^{-2}$ and $\eta_n \to 1$ pointwise as $n \to \infty$. Let $\vartheta \colon [0,\infty) \to [0,1]$ be smooth with $\vartheta = 1$ on $[0,1]$ and $\vartheta = 0$ on $[2,\infty)$, and let $\zeta_n$ be the functions
    \begin{equation*}
        \zeta_n(x,v) :=\eta_n(x)\, \vartheta\Bigl( \frac{U(x) + \frac12 |v|^2 - \inf U}{n^2} \Bigr) .
    \end{equation*}
    Then $\zeta_n \in C_c^\infty(T\mathcal X)$, $0 \Le \zeta_n \Le 1$, $\zeta_n \to 1$ pointwise as $n \to \infty$, and $|v| \Le 2n$ on the support of $\zeta_n$. Since $\mathcal A \bigl( U + \frac12 |v|^2 \bigr) = 0$ and $\norm{\nabla\eta_n}_\infty \Le n^{-2}$,
    \begin{equation}
        |\mathcal A\zeta_n| \Le \frac2n , \qquad
        |\nabla_v\zeta_n| \Le \frac{C}{n} .
        \label{eq: cutoff bounds}
    \end{equation}
    By \eqref{eq: orthogonal} for $\varrho = \zeta_n^2 \psi$ and \eqref{eq: generator form} for $\psi_1 = \zeta_n^2 \psi$ and $\psi_2 = \psi$,
    \begin{equation*}
        0 = \mu\bigl( \psi\, ( \zeta_n^2 \psi - \mathcal L ( \zeta_n^2 \psi ) ) \bigr) = \mu( \zeta_n^2 \psi^2 ) + \mu( \zeta_n^2 \psi\, \mathcal A\psi ) + \gamma\, \mu\bigl( \nabla_v ( \zeta_n^2 \psi ) \cdot \nabla_v\psi \bigr) .
    \end{equation*}
    By the antisymmetry of $\mathcal A$ and the product rule,
    \begin{equation*}
        \mu( \zeta_n^2 \psi\, \mathcal A\psi ) = -\mu( \psi^2\, \zeta_n\, \mathcal A\zeta_n ) , \qquad
        \nabla_v ( \zeta_n^2 \psi ) \cdot \nabla_v\psi \Ge -\psi^2 |\nabla_v\zeta_n|^2 .
    \end{equation*}
    Hence \eqref{eq: cutoff bounds} gives
    \begin{equation*}
        \mu( \zeta_n^2 \psi^2 ) \Le \frac{C ( 1 + \gamma )}{n}\, \mu(\psi^2) ,
    \end{equation*}
    and the limit $n \to \infty$ gives $\psi = 0$, which proves the density.

    To construct $\varrho_n$, note that the semigroup contracts $L^2(\mu)$, so every $\varrho \in \mathrm{Dom}(\mathcal L)$ satisfies $\mu( \varrho\, \mathcal L\varrho ) \Le 0$, and thus
    \begin{equation}
        \norm{\varrho}_{L^2(\mu)} \Le \norm{\varrho - \mathcal L\varrho}_{L^2(\mu)} .
        \label{eq: dissipative}
    \end{equation}
    For $\varrho \in \mathrm{Dom}(\mathcal L)$, the density gives $\varrho_n \in C_c^\infty(T\mathcal X)$ with $\varrho_n - \mathcal L\varrho_n \to \varrho - \mathcal L\varrho$ in $L^2(\mu)$. Hence \eqref{eq: dissipative} for $\varrho_n - \varrho$ in place of $\varrho$ gives $\varrho_n \to \varrho$, and thus $\mathcal L\varrho_n \to \mathcal L\varrho$ in $L^2(\mu)$.
\end{proof}

\subsection{Velocity average and the domain of the generator}
\label{app: domain}

Next, we define the velocity average of a function of $L^2(\mu)$, and in Lemma~\ref{lem: domain} we extend \eqref{eq: adjoints}, \eqref{eq: generator form} and It\^o's formula \eqref{eq: ito} from $C_c^\infty(T\mathcal X)$ to $\mathrm{Dom}(\mathcal L)$. For $\psi \in L^2(\mu)$, the velocity average of $\psi$ is the function $\langle\psi\rangle$ of the position with
\begin{equation}
    \langle\psi\rangle(x) := (2\pi)^{-\frac d2} \int_{T_x\mathcal X} \psi(x,v)\, e^{-\frac12 |v|^2} \D v .
    \label{eq: velocity average}
\end{equation}
Every $\psi \in L^2(\mu)$ and every $g \in L^2(\bar\mu)$ satisfy
\begin{equation}
    \mu( \psi g ) = \mu( \langle\psi\rangle g ) ,
    \label{eq: average identity}
\end{equation}
and thus
\begin{equation}
    \mu( \psi^2 ) = \mu\bigl( \langle\psi\rangle^2 \bigr) + \mu\bigl( ( \psi - \langle\psi\rangle )^2 \bigr) .
    \label{eq: average split}
\end{equation}
For a vector $a \in T_x\mathcal X$ and a symmetric linear map $A$ of $T_x\mathcal X$, the Gaussian distribution of the velocity gives
\begin{equation}
    (2\pi)^{-\frac d2} \int_{T_x\mathcal X} | a - A v |^2 e^{-\frac12 |v|^2} \D v = |a|^2 + |A|^2 .
    \label{eq: gaussian moments}
\end{equation}

\begin{lemma}
\label{lem: domain}
    Let $\varrho \in \mathrm{Dom}(\mathcal L)$. Then $\varrho$ has a velocity gradient $\nabla_v\varrho \in L^2(\mu)$, and every continuous field $W$ on $T\mathcal X$ that is continuously differentiable in the velocity, with $W, \nabla_v^* \cdot W \in L^2(\mu)$, satisfies
    \begin{equation}
        \mu( \varrho\, \nabla_v^* \cdot W ) = \mu( \nabla_v\varrho \cdot W ) .
        \label{eq: domain adjoints}
    \end{equation}
    The function $\varrho$ and its velocity gradient satisfy
    \begin{equation}
        \gamma\, \mu\bigl( |\nabla_v\varrho|^2 \bigr) = -\mu( \varrho\, \mathcal L\varrho ) ,
        \label{eq: dissipation}
    \end{equation}
    and the velocity average $\langle\varrho\rangle$ of \eqref{eq: velocity average} satisfies
    \begin{equation}
        \mu\bigl( ( \varrho - \langle\varrho\rangle )^2 \bigr) \Le \mu\bigl( |\nabla_v\varrho|^2 \bigr) .
        \label{eq: gaussian poincare}
    \end{equation}
    Every $\psi \in C^1(T\mathcal X)$ with $\psi, \mathcal A\psi, \nabla_v\psi \in L^2(\mu)$ satisfies
    \begin{equation}
        \mu\bigl( (-\mathcal L\varrho)\, \psi \bigr) = \mu( \varrho\, \mathcal A\psi ) + \gamma\, \mu( \nabla_v\varrho \cdot \nabla_v\psi ) .
        \label{eq: domain form}
    \end{equation}
    Moreover, the dynamics $(Z_t)_{t \Ge 0}$ with $Z_0 \sim \mu$ satisfies, for every $t \Ge 0$ almost surely,
    \begin{equation}
        \varrho(Z_t) = \varrho(Z_0) + \int_0^t \mathcal L\varrho(Z_s) \,\D s + \sqrt{2\gamma} \int_0^t \nabla_v\varrho(Z_s) \cdot \D B_s ,
        \label{eq: domain ito}
    \end{equation}
    and the right side of \eqref{eq: domain ito} defines a continuous semimartingale modification of $t \mapsto \varrho(Z_t)$.
    Finally, if $\mathcal L\varrho \in C^\infty(T\mathcal X)$, then $\varrho$ is smooth, and $\mathcal L\varrho$ is given by \eqref{eq: generator}.
\end{lemma}

\begin{proof}[Proof of Lemma~\ref{lem: domain}]
    Let $\varrho_n \in C_c^\infty(T\mathcal X)$ be the functions of Lemma~\ref{lem: dynamics} with $\varrho_n \to \varrho$ and $\mathcal L\varrho_n \to \mathcal L\varrho$ in $L^2(\mu)$. By the antisymmetry of $\mathcal A$, \eqref{eq: generator form} for $\psi_1 = \psi_2 = \varrho_n - \varrho_j$ gives
    \begin{equation*}
        \gamma\, \mu\bigl( |\nabla_v ( \varrho_n - \varrho_j )|^2 \bigr) = -\mu\bigl( ( \varrho_n - \varrho_j )\, \mathcal L ( \varrho_n - \varrho_j ) \bigr) \to 0
    \end{equation*}
    as $n, j \to \infty$, so $\nabla_v\varrho_n$ converges in $L^2(\mu)$. By \eqref{eq: adjoints} for $\psi = \varrho_n$, the limit of $\nabla_v\varrho_n$ is a velocity gradient $\nabla_v\varrho$ of $\varrho$ and satisfies \eqref{eq: domain adjoints}. The limit of \eqref{eq: generator form} for $\psi_1 = \psi_2 = \varrho_n$ is \eqref{eq: dissipation}. With $\varrho_n$ in place of $\varrho$, \eqref{eq: gaussian poincare} is the Poincaré inequality of the Gaussian distribution in the velocity \citep[Proposition~4.1.1]{bakry2014analysis}, and \eqref{eq: gaussian poincare} for $\varrho$ follows in the limit $n \to \infty$, since $\langle\varrho_n\rangle \to \langle\varrho\rangle$ in $L^2(\mu)$ by \eqref{eq: average split}. The limit of \eqref{eq: generator form} for $\psi_1 = \varrho_n$ and $\psi_2 = \psi$ is \eqref{eq: domain form}. Since $Z_s \sim \mu$ for every $s$, each term of \eqref{eq: ito} for $\psi = \varrho_n$ converges in mean square, and the limit is \eqref{eq: domain ito}. Since the time integral and the stochastic integral on the right side of \eqref{eq: domain ito} have continuous paths, the right side defines a continuous semimartingale modification of $t \mapsto \varrho(Z_t)$.

    Finally, suppose that $\mathcal L\varrho \in C^\infty(T\mathcal X)$. For $\psi \in C_c^\infty(T\mathcal X)$, \eqref{eq: domain form} and \eqref{eq: domain adjoints} for $W = \nabla_v\psi$ give
    \begin{equation*}
        \mu( \psi\, \mathcal L\varrho ) = \mu\bigl( \varrho\, ( -\mathcal A\psi + \mathcal S\psi ) \bigr) ,
    \end{equation*}
    so the operator \eqref{eq: generator} applied to $\varrho$ in the sense of distributions is the smooth function $\mathcal L\varrho$. This operator satisfies the condition of \citep[Theorem~1.1]{hormander1967hypoelliptic}, so $\varrho$ is smooth, and $\mathcal L\varrho$ is then given by \eqref{eq: generator}.
\end{proof}

\subsection{Eigenfunctions of the diffusion operator}
\label{app: eigenfunctions}

Next, in Lemma~\ref{lem: eigenfunctions}, we give eigenfunctions of $\bar{\mathcal L}$ in \eqref{eq: overdamped} that form an orthonormal basis of $L^2(\bar\mu)$, and we bound the Hessian of the functions $g$ on $\mathcal X$ with $g, \bar{\mathcal L} g \in L^2(\bar\mu)$. Formally, the spectral gap of $\bar{\mathcal L}$ is $m$ by the Poincaré inequality \eqref{eq: poincare}, so the nontrivial eigenvalues of $-\bar{\mathcal L}$ are at least $m$.

\begin{lemma}
\label{lem: eigenfunctions}
    There are functions $\mathfrak b_k \in C^\infty(\mathcal X)$ and numbers $\omega_k \Ge \sqrt m$, $k \Ge 1$, such that the constant function $1$ and the functions $\mathfrak b_k$ form an orthonormal basis of $L^2(\bar\mu)$, and for $j, k \Ge 1$,
    \begin{equation}
        -\bar{\mathcal L}\mathfrak b_k = \omega_k^2 \mathfrak b_k , \qquad
        \mu( \nabla\mathfrak b_j \cdot \nabla\mathfrak b_k ) = \omega_k^2 \mathbf 1_{j = k} .
        \label{eq: eigenfunctions}
    \end{equation}
    Moreover, every $g \in C^\infty(\mathcal X)$ with $g, \bar{\mathcal L} g \in L^2(\bar\mu)$ has a Hessian $\nabla^2 g \in L^2(\bar\mu)$ with
    \begin{equation}
        \mu\bigl( |\nabla^2 g|^2 \bigr) \Le \mu\bigl( ( \bar{\mathcal L} g )^2 \bigr) + K \mu\bigl( |\nabla g|^2 \bigr) .
        \label{eq: hessian bound}
    \end{equation}
\end{lemma}

\begin{proof}[Proof of Lemma~\ref{lem: eigenfunctions}]
    Since $\mathcal X$ is complete, and the lower bound \eqref{eq: curvature} and the Poincaré inequality \eqref{eq: poincare} hold, the diffusion operator $\bar{\mathcal L}$ of \eqref{eq: overdamped} on $C_c^\infty(\mathcal X)$ is essentially self-adjoint in $L^2(\bar\mu)$ \citep[Corollary~4]{eberle2026space}. By the same corollary, every function $g$ in the domain of its self-adjoint extension has $\nabla g, \nabla^2 g \in L^2(\bar\mu)$ and satisfies \eqref{eq: hessian bound}. The bound \eqref{eq: hessian bound} is the integral of the Bochner formula \eqref{eq: bochner formula} against $\bar\mu$ with the lower bound \eqref{eq: curvature}. The identity \eqref{eq: position form} extends from $C_c^\infty(\mathcal X)$ to the functions $g_1$ and $g_2$ of this domain by density.

    By Assumption~\ref{ass: poincare}, the spectrum of the extension is discrete, so $L^2(\bar\mu)$ has an orthonormal basis that consists of the constant function $1$ and of eigenfunctions $\mathfrak b_k$, $k \Ge 1$, of the extension \citep[Section~2]{eberle2026space}. The functions $\mathfrak b_k$ satisfy $-\bar{\mathcal L}\mathfrak b_k = \omega_k^2 \mathfrak b_k$ with numbers $\omega_k > 0$, and they are smooth, since $\bar{\mathcal L}$ is an elliptic operator with smooth coefficients \citep[Corollary~8.11]{gilbarg2001elliptic}. Then \eqref{eq: position form} for $g_1 = \mathfrak b_j$ and $g_2 = \mathfrak b_k$ gives the identity $\mu( \nabla\mathfrak b_j \cdot \nabla\mathfrak b_k ) = \omega_k^2 \mathbf 1_{j = k}$ of \eqref{eq: eigenfunctions}, and the Poincaré inequality \eqref{eq: poincare} for $g = \mathfrak b_k$ gives $\omega_k^2 \Ge m$. Finally, every $g \in C^\infty(\mathcal X)$ with $g, \bar{\mathcal L} g \in L^2(\bar\mu)$ satisfies $\mu( g\, \bar{\mathcal L} h ) = \mu( h\, \bar{\mathcal L} g )$ for every $h \in C_c^\infty(\mathcal X)$, so $g$ belongs to the domain of the extension and satisfies \eqref{eq: hessian bound}.
\end{proof}

\subsection{Divergence equation for one eigenfunction}
\label{app: mode}

We solve the divergence equation of \citep[Theorem~5]{eberle2026space} with the right side $\mathfrak a(s) \mathfrak b_k$, where $\mathfrak a$ is a function of the time and $\mathfrak b_k$ is an eigenfunction of Lemma~\ref{lem: eigenfunctions}. For $T > 0$ and the time $s \in (0,T)$, the divergence equation for $\mathfrak h(s)$ and $\mathfrak g(s)$ with $\mathfrak h(0) = \mathfrak h(T) = 0$ and $\mathfrak g(0) = \mathfrak g(T) = 0$ reads
\begin{equation}
    \partial_s \bigl( \mathfrak h(s) \mathfrak b_k \bigr) - \bar{\mathcal L} \bigl( \mathfrak g(s) \mathfrak b_k \bigr) = \mathfrak a(s) \mathfrak b_k ,
    \label{eq: mode divergence}
\end{equation}
with the derivative $\partial_s$ in the time. Since $-\bar{\mathcal L}\mathfrak b_k = \omega_k^2 \mathfrak b_k$ by \eqref{eq: eigenfunctions}, the left side of \eqref{eq: mode divergence} is $\bigl( \mathfrak h'(s) + \omega_k^2 \mathfrak g(s) \bigr) \mathfrak b_k$. Hence, to solve \eqref{eq: mode divergence}, we only need to solve the equation in the time
\begin{equation*}
    \mathfrak h'(s) + \omega_k^2 \mathfrak g(s) = \mathfrak a(s) , \quad s \in (0,T) , \qquad
    \mathfrak h(0) = \mathfrak h(T) = 0 , \quad \mathfrak g(0) = \mathfrak g(T) = 0 .
\end{equation*}
In Lemma~\ref{lem: mode}, for every $\omega > 0$ in place of $\omega_k$ and every $T \Ge T_0$ with some $T_0 > 0$, we solve the equation in the time, and we bound $\int_0^T \mathfrak g(s)^2 \D s$ and $\int_0^T \bigl( \mathfrak g'(s) - \mathfrak h(s) - \gamma \mathfrak g(s) \bigr)^2 \D s$ by $\int_0^T \mathfrak a(s)^2 \D s$ in \eqref{eq: mode bounds}.

For $T > 0$, $L^2(0,T)$ is the space of the square integrable functions on $(0,T)$ with the norm
\begin{equation*}
    \norm{\mathfrak a}_{L^2(0,T)} := \Bigl( \int_0^T \mathfrak a(s)^2 \,\D s \Bigr)^{\frac12} .
\end{equation*}
The space $H^1(0,T)$ consists of the functions $\mathfrak y$ on $[0,T]$ with $\mathfrak y' \in L^2(0,T)$, $H^2(0,T)$ of the functions $\mathfrak y$ with $\mathfrak y', \mathfrak y'' \in L^2(0,T)$, and $H_0^1(0,T)$ of the functions $\mathfrak y \in H^1(0,T)$ with $\mathfrak y(0) = \mathfrak y(T) = 0$.

For $\omega > 0$, let $\beta_\omega$ be the unique positive root of
\begin{equation}
    \beta_\omega^2 \Bigl( 1 + \frac{K}{\omega^2} + \Bigl( \frac\gamma\omega + \beta_\omega \Bigr)^2 \Bigr) = 1 ,
    \label{eq: mode root}
\end{equation}
which is \eqref{eq: root} with $\omega^2$ in place of $m$, so that $\beta_\omega = \beta(\gamma,\omega^2,K)$ and $\beta_{\sqrt m} = \beta(\gamma,m,K)$. As a function of $\beta_\omega$ and $\omega$, the left side of \eqref{eq: mode root} is increasing in $\beta_\omega$ and nonincreasing in $\omega$, so
\begin{equation}
    \beta_\omega \Ge \beta(\gamma,m,K) , \qquad \forall \omega \Ge \sqrt m .
    \label{eq: root monotone}
\end{equation}

\begin{lemma}
\label{lem: mode}
    Let $\omega > 0$, let $\beta_\omega$ be the root of \eqref{eq: mode root} and let $\varepsilon \in ( 0, \beta_\omega^2 )$. Then there is $T_0 > 0$ such that, for every $T \Ge T_0$ and every $\mathfrak a \in L^2(0,T)$, there are $\mathfrak h, \mathfrak g \in H_0^1(0,T)$ such that
    \begin{equation}
        \mathfrak h'(s) + \omega^2 \mathfrak g(s) = \mathfrak a(s)
        \label{eq: mode equation}
    \end{equation}
    for almost every $s \in (0,T)$, and
    \begin{equation}
        \omega^2 \norm{\mathfrak g' - \mathfrak h - \gamma \mathfrak g}_{L^2(0,T)}^2 + ( \omega^4 + K \omega^2 ) \norm{\mathfrak g}_{L^2(0,T)}^2 \Le \frac{1}{\beta_\omega^2 - \varepsilon} \norm{\mathfrak a}_{L^2(0,T)}^2 .
        \label{eq: mode bounds}
    \end{equation}
\end{lemma}

\begin{proof}[Proof of Lemma~\ref{lem: mode}]
    We first find $T_0 > 0$ such that, for every $T \Ge T_0$, every $\mathfrak y \in H^2(0,T)$ satisfies
    \begin{equation}
        \norm{\mathfrak y'' + \gamma \mathfrak y' + \omega^2 \mathfrak y}_{L^2(0,T)}^2 + ( \omega^2 + K ) \norm{\mathfrak y'}_{L^2(0,T)}^2 \Ge ( \beta_\omega^2 - \varepsilon ) \omega^2 ( \omega^2 + K ) \norm{\mathfrak y}_{L^2(0,T)}^2 .
        \label{eq: mode coercivity}
    \end{equation}
    For $T > 0$, let $\mathfrak w$ be the linear function
    \begin{equation}
        \mathfrak w(s) := \frac{\omega ( 1 - \beta_\omega^2 )}{\beta_\omega} \Bigl( 1 - \frac{2s}{T} \Bigr) , \qquad s \in [0,T] .
        \label{eq: mode weight}
    \end{equation}
    Integration by parts and \eqref{eq: mode root} give, for every $\mathfrak y \in H^2(0,T)$,
    \begin{equation}
    \begin{aligned}
        & \norm{\mathfrak y'' + \gamma \mathfrak y' + \omega^2 \mathfrak y}_{L^2(0,T)}^2 + ( \omega^2 + K ) \norm{\mathfrak y'}_{L^2(0,T)}^2 - \bigl\lVert \mathfrak y'' + \mathfrak w \mathfrak y' + \omega \beta_\omega ( \gamma + \omega \beta_\omega ) \mathfrak y \bigr\rVert_{L^2(0,T)}^2 \\
        & \qquad = \int_0^T \bigl( \mathfrak w(0)^2 - \mathfrak w^2 + \mathfrak w' \bigr) ( \mathfrak y' )^2 \,\D s + \int_0^T \bigl( \beta_\omega^2 \omega^2 ( \omega^2 + K ) + \omega \beta_\omega ( \gamma + \omega \beta_\omega ) \mathfrak w' \bigr) \mathfrak y^2 \,\D s \\
        & \qquad \quad + \Bigl[ 2 \bigl( \omega^2 - \omega \beta_\omega ( \gamma + \omega \beta_\omega ) \bigr) \mathfrak y \mathfrak y' + ( \gamma - \mathfrak w ) ( \mathfrak y' )^2 + \bigl( \gamma \omega^2 - \omega \beta_\omega ( \gamma + \omega \beta_\omega ) \mathfrak w \bigr) \mathfrak y^2 \Bigr]_{s=0}^{s=T} .
    \end{aligned}
        \label{eq: mode identity}
    \end{equation}
    Since $\beta_\omega ( \gamma + \omega \beta_\omega ) < \omega$ by \eqref{eq: mode root}, \eqref{eq: mode weight} gives $\mathfrak w(T) = -\mathfrak w(0)$ and
    \begin{equation*}
        \mathfrak w(0) - \gamma = \frac{\omega - \beta_\omega ( \gamma + \omega \beta_\omega )}{\beta_\omega} > 0 , \qquad
        \omega \beta_\omega ( \gamma + \omega \beta_\omega ) \mathfrak w(0) - \gamma \omega^2 = \omega^2 \beta_\omega^2 \bigl( \mathfrak w(0) - \gamma \bigr) .
    \end{equation*}
    Hence $\omega^2 - \omega \beta_\omega ( \gamma + \omega \beta_\omega ) = \omega \beta_\omega \bigl( \mathfrak w(0) - \gamma \bigr)$, and the boundary term of \eqref{eq: mode identity} is
    \begin{equation*}
        \bigl( \mathfrak w(0) - \gamma \bigr) \bigl( \mathfrak y'(T) + \omega \beta_\omega \mathfrak y(T) \bigr)^2 + 2\gamma \bigl( \mathfrak y'(T)^2 + \omega^2 \mathfrak y(T)^2 \bigr) + \bigl( \mathfrak w(0) - \gamma \bigr) \bigl( \mathfrak y'(0) - \omega \beta_\omega \mathfrak y(0) \bigr)^2 \Ge 0 .
    \end{equation*}
    Since $\mathfrak w(s)^2 \Le \mathfrak w(0)^2$ and $\mathfrak w'(s) = -\frac{2 \mathfrak w(0)}{T}$ for $s \in [0,T]$, \eqref{eq: mode identity} gives
    \begin{equation*}
    \begin{aligned}
        & \Bigl( 1 + \frac{2 \mathfrak w(0)}{( \omega^2 + K ) T} \Bigr) \Bigl( \norm{\mathfrak y'' + \gamma \mathfrak y' + \omega^2 \mathfrak y}_{L^2(0,T)}^2 + ( \omega^2 + K ) \norm{\mathfrak y'}_{L^2(0,T)}^2 \Bigr) \\
        & \qquad \Ge \Bigl( \beta_\omega^2 \omega^2 ( \omega^2 + K ) - \frac{2 \omega \beta_\omega ( \gamma + \omega \beta_\omega ) \mathfrak w(0)}{T} \Bigr) \norm{\mathfrak y}_{L^2(0,T)}^2 .
    \end{aligned}
    \end{equation*}
    Since $\mathfrak w(0)$ does not depend on $T$, there is $T_0 > 0$ such that \eqref{eq: mode coercivity} holds for every $T \Ge T_0$.

    To construct $\mathfrak h$ and $\mathfrak g$, let $T \Ge T_0$ and $\mathfrak a \in L^2(0,T)$. By \eqref{eq: mode coercivity}, the functional
    \begin{equation*}
        \mathfrak y \mapsto \frac12 \norm{\mathfrak y'' + \gamma \mathfrak y' + \omega^2 \mathfrak y}_{L^2(0,T)}^2 + \frac{\omega^2 + K}{2} \norm{\mathfrak y'}_{L^2(0,T)}^2 - ( \omega^2 + K ) \int_0^T \mathfrak a \mathfrak y \,\D s
    \end{equation*}
    has a minimizer $\mathfrak u$ in $H^2(0,T)$, which satisfies, for every $\mathfrak y \in H^2(0,T)$,
    \begin{equation}
        \int_0^T \bigl( ( \mathfrak u'' + \gamma \mathfrak u' + \omega^2 \mathfrak u ) ( \mathfrak y'' + \gamma \mathfrak y' + \omega^2 \mathfrak y ) + ( \omega^2 + K ) \mathfrak u' \mathfrak y' \bigr) \D s = ( \omega^2 + K ) \int_0^T \mathfrak a \mathfrak y \,\D s .
        \label{eq: mode minimizer}
    \end{equation}
    Let $\mathfrak g$ and $\mathfrak h$ be the functions
    \begin{equation*}
        \mathfrak g(s) :=\frac{\mathfrak u''(s) + \gamma \mathfrak u'(s) + \omega^2 \mathfrak u(s)}{\omega^2 + K} , \qquad
        \mathfrak h(s) := \int_0^s \bigl( \mathfrak a(r) - \omega^2 \mathfrak g(r) \bigr) \D r .
    \end{equation*}
    Then $\mathfrak h$ belongs to $H^1(0,T)$ with $\mathfrak h(0) = 0$ and satisfies \eqref{eq: mode equation}, so an integration by parts turns \eqref{eq: mode minimizer} into
    \begin{equation}
        \int_0^T \bigl( \mathfrak g \mathfrak y'' + ( \gamma \mathfrak g + \mathfrak h + \mathfrak u' ) \mathfrak y' \bigr) \D s = \mathfrak h(T) \mathfrak y(T) .
        \label{eq: mode weak}
    \end{equation}
    For $\mathfrak y = 1$, \eqref{eq: mode weak} gives $\mathfrak h(T) = 0$. Since $\mathfrak y'$ ranges over $H^1(0,T)$, \eqref{eq: mode weak} then gives $\mathfrak g \in H_0^1(0,T)$ with
    \begin{equation*}
        \mathfrak g'(s) - \mathfrak h(s) - \gamma \mathfrak g(s) = \mathfrak u'(s) .
    \end{equation*}
    By the Cauchy--Schwarz inequality, \eqref{eq: mode minimizer} and \eqref{eq: mode coercivity} for $\mathfrak y = \mathfrak u$ give
    \begin{equation*}
    \begin{aligned}
        & \omega^2 \norm{\mathfrak g' - \mathfrak h - \gamma \mathfrak g}_{L^2(0,T)}^2 + ( \omega^4 + K \omega^2 ) \norm{\mathfrak g}_{L^2(0,T)}^2 = \omega^2 \int_0^T \mathfrak a \mathfrak u \,\D s \\
        & \qquad \Le \frac{\norm{\mathfrak a}_{L^2(0,T)}}{( \beta_\omega^2 - \varepsilon )^{\frac12}} \Bigl( \omega^2 \norm{\mathfrak g' - \mathfrak h - \gamma \mathfrak g}_{L^2(0,T)}^2 + ( \omega^4 + K \omega^2 ) \norm{\mathfrak g}_{L^2(0,T)}^2 \Bigr)^{\frac12} ,
    \end{aligned}
    \end{equation*}
    and thus \eqref{eq: mode bounds}.
\end{proof}

\subsection{Time integral of the squared norm of the semigroup}
\label{app: window}

In Lemma~\ref{lem: window}, we bound $\int_0^\infty \mu\bigl( ( e^{s\mathcal L}\psi )^2 \bigr) \D s$ for $\psi \in \mathrm{Dom}(\mathcal L)$ with $\mu(\psi) = 0$. For $s \Ge 0$, $e^{s\mathcal L}\psi \in \mathrm{Dom}(\mathcal L)$, $\mu( e^{s\mathcal L}\psi ) = 0$ and
\begin{equation}
    \partial_s ( e^{s\mathcal L}\psi ) = \mathcal L ( e^{s\mathcal L}\psi )
    \label{eq: semigroup equation}
\end{equation}
in $L^2(\mu)$, so \eqref{eq: dissipation} gives
\begin{equation}
    \frac{\D}{\D s}\, \mu\bigl( ( e^{s\mathcal L}\psi )^2 \bigr) = 2 \mu\bigl( e^{s\mathcal L}\psi\, \mathcal L ( e^{s\mathcal L}\psi ) \bigr) = -2\gamma\, \mu\bigl( |\nabla_v ( e^{s\mathcal L}\psi )|^2 \bigr) .
    \label{eq: energy derivative}
\end{equation}
The proof of Lemma~\ref{lem: window} tests \eqref{eq: semigroup equation} against the test function $\mathfrak h + v \cdot \nabla\mathfrak g$ of \citep[Lemma~14]{eberle2026space}, where $\mathfrak h$ and $\mathfrak g$ are sums of solutions of the divergence equation \eqref{eq: mode divergence}.

\begin{lemma}
\label{lem: window}
    Every $\psi \in \mathrm{Dom}(\mathcal L)$ with $\mu(\psi) = 0$ satisfies, with the rate $\lambda$ of \eqref{eq: l2 decay},
    \begin{equation}
        \int_0^\infty \mu\bigl( ( e^{s\mathcal L}\psi )^2 \bigr) \D s \Le \frac{\mu(\psi^2)}{2\lambda} .
        \label{eq: window}
    \end{equation}
\end{lemma}

\begin{proof}[Proof of Lemma~\ref{lem: window}]
    To prove \eqref{eq: window}, let $\mathfrak b_k$ and $\omega_k$, $k \Ge 1$, be as in Lemma~\ref{lem: eigenfunctions}, and let $\mathfrak a_k$ be the coefficients
    \begin{equation}
        \mathfrak a_k(s) :=\mu( \mathfrak b_k\, e^{s\mathcal L}\psi ) , \qquad k \Ge 1 , \quad s \Ge 0 .
        \label{eq: coefficients}
    \end{equation}
    By \eqref{eq: average identity}, $\mu\bigl( \langle e^{s\mathcal L}\psi\rangle \bigr) = 0$ and $\mu\bigl( \mathfrak b_k \langle e^{s\mathcal L}\psi\rangle \bigr) = \mathfrak a_k(s)$, so the basis of Lemma~\ref{lem: eigenfunctions} gives
    \begin{equation*}
        \mu\bigl( \langle e^{s\mathcal L}\psi\rangle^2 \bigr) = \sum_{k=1}^\infty \mathfrak a_k(s)^2 .
    \end{equation*}
    Hence \eqref{eq: average split} and \eqref{eq: gaussian poincare} for $e^{s\mathcal L}\psi$ give
    \begin{equation}
        \mu\bigl( ( e^{s\mathcal L}\psi )^2 \bigr) \Le \sum_{k=1}^\infty \mathfrak a_k(s)^2 + \mu\bigl( |\nabla_v ( e^{s\mathcal L}\psi )|^2 \bigr) .
        \label{eq: split}
    \end{equation}
    By \eqref{eq: energy derivative}, for every $T > 0$,
    \begin{equation}
        2\gamma \int_0^T \mu\bigl( |\nabla_v ( e^{s\mathcal L}\psi )|^2 \bigr) \D s = \mu(\psi^2) - \mu\bigl( ( e^{T\mathcal L}\psi )^2 \bigr) \Le \mu(\psi^2) .
        \label{eq: dissipation integral}
    \end{equation}

    To bound $\int_0^\infty \sum_{k=1}^\infty \mathfrak a_k(s)^2 \,\D s$, we fix $n \Ge 1$ and $\varepsilon \in \bigl( 0, \beta(\gamma,m,K)^2 \bigr)$, with the root $\beta(\gamma,m,K)$ of \eqref{eq: root}. By \eqref{eq: root monotone}, $\varepsilon < \beta_{\omega_k}^2$ for every $k \Ge 1$. Hence Lemma~\ref{lem: mode} for $\omega = \omega_1, \dots, \omega_n$ gives a common $T_0 > 0$ such that, for every $T \Ge T_0$ and $k = 1, \dots, n$, there are $\mathfrak h_k, \mathfrak g_k \in H_0^1(0,T)$ that satisfy \eqref{eq: mode equation} and \eqref{eq: mode bounds} for $\omega = \omega_k$, $\mathfrak a = \mathfrak a_k$, $\mathfrak h = \mathfrak h_k$ and $\mathfrak g = \mathfrak g_k$. We fix $T \Ge T_0$ and define the functions of the time and the position
    \begin{equation*}
        \mathfrak h(s,x) := \sum_{k=1}^n \mathfrak h_k(s) \mathfrak b_k(x) , \qquad
        \mathfrak g(s,x) := \sum_{k=1}^n \mathfrak g_k(s) \mathfrak b_k(x) .
    \end{equation*}
    By \eqref{eq: eigenfunctions} and \eqref{eq: mode equation}, $\mathfrak h$ and $\mathfrak g$ solve the divergence equation
    \begin{equation}
        \partial_s \mathfrak h - \bar{\mathcal L}\mathfrak g = \sum_{k=1}^n \mathfrak a_k(s) \mathfrak b_k .
        \label{eq: divergence}
    \end{equation}
    Since $\nabla_x ( v \cdot \nabla\mathfrak g ) = \nabla^2\mathfrak g\, v$ by \eqref{eq: position gradient} and $\nabla_v ( v \cdot \nabla\mathfrak g ) = \nabla\mathfrak g$, \eqref{eq: generator} gives
    \begin{equation}
        \mathcal A ( \mathfrak h + v \cdot \nabla\mathfrak g ) = v \cdot \nabla\mathfrak h + v \cdot \nabla^2\mathfrak g\, v - \nabla U \cdot \nabla\mathfrak g .
        \label{eq: test transport}
    \end{equation}
    Let $W$ be the field
    \begin{equation}
        W := \partial_s \nabla\mathfrak g - \nabla\mathfrak h - \nabla^2\mathfrak g\, v ,
        \label{eq: field}
    \end{equation}
    which satisfies
    \begin{equation}
        \nabla_v^* \cdot W = \Delta\mathfrak g + v \cdot ( \partial_s \nabla\mathfrak g - \nabla\mathfrak h ) - v \cdot \nabla^2\mathfrak g\, v .
        \label{eq: field divergence}
    \end{equation}
    Then \eqref{eq: overdamped} and \eqref{eq: divergence}--\eqref{eq: field divergence} give
    \begin{equation}
        ( \partial_s - \mathcal A ) ( \mathfrak h + v \cdot \nabla\mathfrak g ) = \sum_{k=1}^n \mathfrak a_k(s) \mathfrak b_k + \nabla_v^* \cdot W .
        \label{eq: test identity}
    \end{equation}

    We multiply \eqref{eq: test identity} by $e^{s\mathcal L}\psi$ and integrate against $\mu$ and over $s \in [0,T]$. By Lemma~\ref{lem: eigenfunctions}, for almost every $s$, the function $\mathfrak h + v \cdot \nabla\mathfrak g$ and the field $W$ satisfy the conditions of \eqref{eq: domain form} and \eqref{eq: domain adjoints} in Lemma~\ref{lem: domain}. Since $\mathfrak h_k(0) = \mathfrak h_k(T) = 0$ and $\mathfrak g_k(0) = \mathfrak g_k(T) = 0$ for $k = 1, \dots, n$, integration by parts in $s$ with \eqref{eq: semigroup equation}, and \eqref{eq: domain form} for $\varrho = e^{s\mathcal L}\psi$ and the test function $\mathfrak h + v \cdot \nabla\mathfrak g$, turn the left side of \eqref{eq: test identity} into
    \begin{equation*}
    \begin{aligned}
        & \int_0^T \mu\bigl( e^{s\mathcal L}\psi\, ( \partial_s - \mathcal A ) ( \mathfrak h + v \cdot \nabla\mathfrak g ) \bigr) \D s \\
        &  \quad = \int_0^T \Bigl( \mu\bigl( ( -\mathcal L ( e^{s\mathcal L}\psi ) ) ( \mathfrak h + v \cdot \nabla\mathfrak g ) \bigr) - \mu\bigl( e^{s\mathcal L}\psi\, \mathcal A ( \mathfrak h + v \cdot \nabla\mathfrak g ) \bigr) \Bigr) \D s = \gamma \int_0^T \mu( \nabla_v ( e^{s\mathcal L}\psi ) \cdot \nabla\mathfrak g ) \,\D s .
    \end{aligned}
    \end{equation*}
    By \eqref{eq: coefficients} and \eqref{eq: domain adjoints} for $\varrho = e^{s\mathcal L}\psi$, the right side of \eqref{eq: test identity} gives
    \begin{equation*}
        \int_0^T \mu\Bigl( e^{s\mathcal L}\psi\, \Bigl( \sum_{k=1}^n \mathfrak a_k(s) \mathfrak b_k + \nabla_v^* \cdot W \Bigr) \Bigr) \D s = \sum_{k=1}^n \norm{\mathfrak a_k}_{L^2(0,T)}^2 + \int_0^T \mu( \nabla_v ( e^{s\mathcal L}\psi ) \cdot W ) \,\D s .
    \end{equation*}
    Hence
    \begin{equation}
        \sum_{k=1}^n \norm{\mathfrak a_k}_{L^2(0,T)}^2 = \int_0^T \mu\bigl( \nabla_v ( e^{s\mathcal L}\psi ) \cdot ( \gamma \nabla\mathfrak g - W ) \bigr) \D s .
        \label{eq: window duality}
    \end{equation}

    To bound the right side of \eqref{eq: window duality}, we write, by \eqref{eq: field},
    \begin{equation*}
        \gamma \nabla\mathfrak g - W = \nabla^2\mathfrak g\, v - ( \partial_s \nabla\mathfrak g - \nabla\mathfrak h - \gamma \nabla\mathfrak g ) .
    \end{equation*}
    Then \eqref{eq: gaussian moments} and \eqref{eq: hessian bound} for $g = \mathfrak g(s,\cdot)$ give
    \begin{equation*}
    \begin{aligned}
        \mu\bigl( |\gamma \nabla\mathfrak g - W|^2 \bigr) & = \mu\bigl( |\partial_s \nabla\mathfrak g - \nabla\mathfrak h - \gamma \nabla\mathfrak g|^2 \bigr) + \mu\bigl( |\nabla^2\mathfrak g|^2 \bigr) \\
        & \Le \mu\bigl( |\partial_s \nabla\mathfrak g - \nabla\mathfrak h - \gamma \nabla\mathfrak g|^2 \bigr) + \mu\bigl( ( \bar{\mathcal L}\mathfrak g )^2 \bigr) + K \mu\bigl( |\nabla\mathfrak g|^2 \bigr) .
    \end{aligned}
    \end{equation*}
    Since the functions $\mathfrak b_k$ are orthonormal and satisfy \eqref{eq: eigenfunctions}, the bounds \eqref{eq: mode bounds} and \eqref{eq: root monotone} give
    \begin{align}
        \int_0^T \mu\bigl( |\gamma \nabla\mathfrak g - W|^2 \bigr) \D s
        & \Le \sum_{k=1}^n \Bigl( \omega_k^2 \norm{\mathfrak g_k' - \mathfrak h_k - \gamma \mathfrak g_k}_{L^2(0,T)}^2 + ( \omega_k^4 + K \omega_k^2 ) \norm{\mathfrak g_k}_{L^2(0,T)}^2 \Bigr) \notag \\
        & \Le \frac{1}{\beta(\gamma,m,K)^2 - \varepsilon} \sum_{k=1}^n \norm{\mathfrak a_k}_{L^2(0,T)}^2 .
        \label{eq: field norm}
    \end{align}
    By the Cauchy--Schwarz inequality, \eqref{eq: window duality}, \eqref{eq: field norm} and \eqref{eq: dissipation integral} give
    \begin{equation}
        \sum_{k=1}^n \norm{\mathfrak a_k}_{L^2(0,T)}^2 \Le \frac{1}{\beta(\gamma,m,K)^2 - \varepsilon} \int_0^T \mu\bigl( |\nabla_v ( e^{s\mathcal L}\psi )|^2 \bigr) \D s \Le \frac{\mu(\psi^2)}{2\gamma \bigl( \beta(\gamma,m,K)^2 - \varepsilon \bigr)} .
        \label{eq: coefficient bound}
    \end{equation}
    The limits $T \to \infty$, then $\varepsilon \to 0$ and then $n \to \infty$ of \eqref{eq: coefficient bound} give
    \begin{equation*}
        \int_0^\infty \sum_{k=1}^\infty \mathfrak a_k(s)^2 \,\D s \Le \frac{\mu(\psi^2)}{2\gamma\, \beta(\gamma,m,K)^2} .
    \end{equation*}
    Hence the integral of \eqref{eq: split} over $s \Ge 0$, \eqref{eq: dissipation integral} and the definition of $\lambda$ in \eqref{eq: l2 decay} give
    \begin{equation*}
        \int_0^\infty \mu\bigl( ( e^{s\mathcal L}\psi )^2 \bigr) \D s \Le \Bigl( \frac{1}{\beta(\gamma,m,K)^2} + 1 \Bigr) \frac{\mu(\psi^2)}{2\gamma} = \frac{\mu(\psi^2)}{2\lambda} ,
    \end{equation*}
    which is \eqref{eq: window}.
\end{proof}

\subsection{Proof of Theorem~\ref{thm: decay}}
\label{app: decay proof}

Finally, we prove Theorem~\ref{thm: decay} by Lemma~\ref{lem: window}.

\begin{proof}[Proof of Theorem~\ref{thm: decay}]
    To prove \eqref{eq: l2 decay}, let $\psi \in \mathrm{Dom}(\mathcal L)$ with $\mu(\psi) = 0$. For $t \Ge 0$, the bound \eqref{eq: window} for $e^{t\mathcal L}\psi$ in place of $\psi$ gives
    \begin{equation*}
        \frac{\D}{\D t} \int_t^\infty \mu\bigl( ( e^{s\mathcal L}\psi )^2 \bigr) \D s = -\mu\bigl( ( e^{t\mathcal L}\psi )^2 \bigr) \Le -2\lambda \int_t^\infty \mu\bigl( ( e^{s\mathcal L}\psi )^2 \bigr) \D s ,
    \end{equation*}
    and thus, with \eqref{eq: window} for $\psi$,
    \begin{equation}
        \int_t^\infty \mu\bigl( ( e^{s\mathcal L}\psi )^2 \bigr) \D s \Le e^{-2\lambda t} \int_0^\infty \mu\bigl( ( e^{s\mathcal L}\psi )^2 \bigr) \D s \Le \frac{e^{-2\lambda t}}{2\lambda}\, \mu(\psi^2) .
        \label{eq: window decay}
    \end{equation}
    By \eqref{eq: energy derivative}, $s \mapsto \mu\bigl( ( e^{s\mathcal L}\psi )^2 \bigr)$ is nonincreasing. Hence, for $t \Ge \frac{1}{2\lambda}$, the bound \eqref{eq: window decay} at the time $t - \frac{1}{2\lambda}$ gives
    \begin{equation*}
        \mu\bigl( ( e^{t\mathcal L}\psi )^2 \bigr) \Le 2\lambda \int_{t - \frac{1}{2\lambda}}^t \mu\bigl( ( e^{s\mathcal L}\psi )^2 \bigr) \D s \Le e^{1 - 2\lambda t}\, \mu(\psi^2) .
    \end{equation*}
    For $t < \frac{1}{2\lambda}$, \eqref{eq: l2 decay} holds, since the semigroup contracts $L^2(\mu)$ and $e^{1 - 2\lambda t} \Ge 1$. By density, \eqref{eq: l2 decay} extends from $\mathrm{Dom}(\mathcal L)$ to every $\psi \in L^2(\mu)$ with $\mu(\psi) = 0$.

    Finally, to prove \eqref{eq: rate lower}, note that \eqref{eq: root} gives
    \begin{equation*}
        \frac{m}{\beta^2} = m + K + ( \gamma + \sqrt m \beta )^2 .
    \end{equation*}
    Hence the definition of $\lambda$ in \eqref{eq: l2 decay} gives
    \begin{equation}
        \lambda = \frac{\gamma}{1 + \frac{1}{\beta^2}} = \frac{m\gamma}{( \gamma + \sqrt m \beta )^2 + 2m + K} .
        \label{eq: rate form}
    \end{equation}
    Since \eqref{eq: root} gives $\beta^2 ( 1 + \beta^2 ) \Le 1$, that is, $\beta^2 \Le \frac{\sqrt5 - 1}{2}$, \eqref{eq: rate form} gives \eqref{eq: rate lower}.
\end{proof}

\section{Proofs for the regularity of the Poisson solution}
\label{app: poisson}

This appendix proves Theorem~\ref{thm: poisson}, Corollaries~\ref{cor: variance} and~\ref{cor: second moment} and the bounds \eqref{eq: variance overdamped} of Remark~\ref{rem: overdamped}. The Poisson solution \eqref{eq: poisson solution} and these statements are unchanged when a constant is added to $f$, so the proofs of this appendix assume $\mu(f) = 0$.

\subsection{Decay of the semigroup in $L^{\frac43}(\mu)$ and $L^4(\mu)$}
\label{app: lp decay}

In Lemma~\ref{lem: decay}, we extend the decay \eqref{eq: l2 decay} to $L^{\frac43}(\mu)$ and $L^4(\mu)$.

\begin{lemma}
\label{lem: decay}
    For $p = \frac43$ and $p = 4$, every $\psi \in L^p(\mu)$ satisfies, with the rate $\lambda$ of \eqref{eq: l2 decay},
    \begin{equation}
        \norm{e^{t\mathcal L}\psi - \mu(\psi)}_{L^p(\mu)} \Le \sqrt2 \min\bigl( 1, e^{\frac14 - \frac{\lambda t}2} \bigr) \norm{\psi}_{L^p(\mu)} , \qquad \forall t \Ge 0 .
        \label{eq: lp decay}
    \end{equation}
\end{lemma}

\begin{proof}[Proof of Lemma~\ref{lem: decay}]
    By Lemma~\ref{lem: dynamics}, the semigroup contracts $L^1(\mu)$ and $L^2(\mu)$, and $\sup | e^{t\mathcal L}\psi | \Le \sup |\psi|$ for every bounded function $\psi$, since $e^{t\mathcal L}\psi$ is an expectation. Hence every bounded function $\psi$ satisfies
    \begin{equation*}
        \norm{e^{t\mathcal L}\psi - \mu(\psi)}_{L^1(\mu)} \Le 2 \norm{\psi}_{L^1(\mu)} , \qquad
        \sup \bigl| e^{t\mathcal L}\psi - \mu(\psi) \bigr| \Le 2 \sup |\psi| ,
    \end{equation*}
    and the contraction of $L^2(\mu)$ with \eqref{eq: l2 decay} for $\psi - \mu(\psi)$ in place of $\psi$ gives
    \begin{equation}
        \norm{e^{t\mathcal L}\psi - \mu(\psi)}_{L^2(\mu)} \Le \min\bigl( 1, e^{\frac12 - \lambda t} \bigr) \norm{\psi}_{L^2(\mu)} .
        \label{eq: l2 contraction}
    \end{equation}
    Since $e^{t\mathcal L}$ is a real Markov operator, these bounds in $L^1(\mu)$, $L^2(\mu)$ and the supremum norm remain valid for complex-valued functions. Applying the Riesz--Thorin interpolation theorem \citep[Theorem~1.1.1]{bergh1976interpolation} with equal weights between $L^1(\mu)$ and $L^2(\mu)$, and between $L^2(\mu)$ and the supremum norm, gives \eqref{eq: lp decay}.
\end{proof}

\subsection{Poisson solution of a function in $L^2(\mu)$}
\label{app: solutions}

We define the Poisson solution of a function of the position and the velocity. For $F \in L^2(\mu)$ with $\mu(F) = 0$, the Poisson solution of $F$ is
\begin{equation}
    \Phi := \int_0^\infty e^{t\mathcal L} F \,\D t ,
    \label{eq: resolvent}
\end{equation}
so that, for $F = f$, $\Phi$ is the Poisson solution $\phi$ of \eqref{eq: poisson solution}. In Lemma~\ref{lem: poisson}, we show that $\Phi \in \mathrm{Dom}(\mathcal L)$ and bound $\norm{\Phi}_{L^2(\mu)}$, $\norm{\nabla_v\Phi}_{L^2(\mu)}$ and $\norm{\Phi}_{L^4(\mu)}$.

\begin{lemma}
\label{lem: poisson}
    Let $F \in L^2(\mu)$ with $\mu(F) = 0$. Then the integral \eqref{eq: resolvent} converges in $L^2(\mu)$, and the Poisson solution $\Phi$ of $F$ satisfies $\mu(\Phi) = 0$, belongs to $\mathrm{Dom}(\mathcal L)$ and solves $\mathcal L\Phi = -F$. With the rate $\lambda$ of \eqref{eq: l2 decay},
    \begin{equation}
        \norm{\Phi}_{L^2(\mu)} \Le \frac2\lambda \norm{F}_{L^2(\mu)} , \qquad
        \gamma \mu\bigl( |\nabla_v\Phi|^2 \bigr) = \mu(\Phi\, F) \Le \frac2\lambda \norm{F}_{L^2(\mu)}^2 .
        \label{eq: resolvent bound}
    \end{equation}
    If $F \in L^4(\mu)$, then
    \begin{equation}
        \norm{\Phi}_{L^4(\mu)} \Le \frac1\lambda \Bigl( \frac52 + \log 2 \Bigr) \norm{F}_{L^4(\mu)} .
        \label{eq: resolvent fourth}
    \end{equation}
\end{lemma}

\begin{proof}[Proof of Lemma~\ref{lem: poisson}]
    Minkowski's integral inequality with \eqref{eq: l2 decay} gives the convergence of \eqref{eq: resolvent} in $L^2(\mu)$ and the bound of $\norm{\Phi}_{L^2(\mu)}$ in \eqref{eq: resolvent bound}. The invariance of $\mu$ gives $\mu(\Phi) = 0$. The semigroup property gives
    \begin{equation}
        e^{s\mathcal L}\Phi - \Phi = -\int_0^s e^{r\mathcal L}F \,\D r , \qquad \forall s \Ge 0 .
        \label{eq: semigroup integral}
    \end{equation}
    Dividing \eqref{eq: semigroup integral} by $s$ and letting $s \to 0$ gives $\Phi \in \mathrm{Dom}(\mathcal L)$ and $\mathcal L\Phi = -F$. Then \eqref{eq: dissipation} for $\varrho = \Phi$ is the identity of \eqref{eq: resolvent bound}, and the Cauchy--Schwarz inequality with the bound of $\norm{\Phi}_{L^2(\mu)}$ gives the bound of $\mu(\Phi\, F)$.

    To prove \eqref{eq: resolvent fourth}, let $F \in L^4(\mu)$. Minkowski's integral inequality, the contraction of $L^4(\mu)$ and \eqref{eq: lp decay} for $p = 4$ give
    \begin{equation*}
        \norm{\Phi}_{L^4(\mu)} \Le \int_0^\infty \min\bigl( 1, \sqrt2 e^{\frac14 - \frac{\lambda t}2} \bigr) \D t\, \norm{F}_{L^4(\mu)} = \frac1\lambda \Bigl( \frac52 + \log 2 \Bigr) \norm{F}_{L^4(\mu)} .
    \end{equation*}
\end{proof}

\subsection{Poisson equation in the position}
\label{app: position}

We consider the Poisson equation of the diffusion operator $\bar{\mathcal L}$ for a function $g$ of the position with $\mu(g) = 0$,
\begin{equation}
    -\bar{\mathcal L} u = g , \qquad \mu(u) = 0 .
    \label{eq: position equation}
\end{equation}
In Lemma~\ref{lem: position}, we solve \eqref{eq: position equation} and bound $\nabla u$ and $\nabla^2 u$ in $L^2(\bar\mu)$.

\begin{lemma}
\label{lem: position}
    Let $g \in C^\infty(\mathcal X) \cap L^2(\bar\mu)$ with $\mu(g) = 0$. Then \eqref{eq: position equation} has a smooth solution $u$ on $\mathcal X$ with
    \begin{equation}
        \mu\bigl( |\nabla u|^2 \bigr) \Le \frac1m \mu(g^2) , \qquad
        \mu\bigl( |\nabla^2 u|^2 \bigr) \Le \frac{m + K}{m} \mu(g^2) .
        \label{eq: position bounds}
    \end{equation}
\end{lemma}

\begin{proof}[Proof of Lemma~\ref{lem: position}]
    By the Poincaré inequality \eqref{eq: poincare} (extended by density to functions with a weak gradient in $L^2(\bar\mu)$) and the Riesz representation theorem, there is a function $u \in L^2(\bar\mu)$ with weak gradient $\nabla u \in L^2(\bar\mu)$ and $\mu(u) = 0$ such that, for every function $h \in L^2(\bar\mu)$ with weak gradient $\nabla h \in L^2(\bar\mu)$,
    \begin{equation}
        \mu( \nabla u \cdot \nabla h ) = \mu( g h ) .
        \label{eq: position weak}
    \end{equation}
    For $h \in C_c^\infty(\mathcal X)$, \eqref{eq: position weak} is the weak form of $-\bar{\mathcal L} u = g$. Since $\bar{\mathcal L}$ is an elliptic operator with smooth coefficients and $g$ is smooth, $u$ is smooth with $-\bar{\mathcal L} u = g$ \citep[Corollary~8.11]{gilbarg2001elliptic}. Then \eqref{eq: position weak} for $h = u$ and \eqref{eq: poincare} for $u$ give the bound of $\mu\bigl( |\nabla u|^2 \bigr)$ in \eqref{eq: position bounds}. Since $u, \bar{\mathcal L} u \in L^2(\bar\mu)$, Lemma~\ref{lem: eigenfunctions} gives $\nabla^2 u \in L^2(\bar\mu)$ and the bound \eqref{eq: hessian bound} for the function $u$,
    \begin{equation*}
        \mu\bigl( |\nabla^2 u|^2 \bigr) \Le \mu(g^2) + K \mu\bigl( |\nabla u|^2 \bigr) ,
    \end{equation*}
    and the bound of $\mu\bigl( |\nabla u|^2 \bigr)$ in \eqref{eq: position bounds} then gives the bound of $\mu\bigl( |\nabla^2 u|^2 \bigr)$.
\end{proof}

\subsection{Proof of Theorem~\ref{thm: poisson}}
\label{app: conclusion}

We prove Theorem~\ref{thm: poisson} by Lemmas~\ref{lem: poisson} and~\ref{lem: position}.

\begin{proof}[Proof of Theorem~\ref{thm: poisson}]
    By Lemma~\ref{lem: poisson} for $F = f$, the integral \eqref{eq: poisson solution} converges in $L^2(\mu)$, and $\phi$ satisfies $\mu(\phi) = 0$, belongs to $\mathrm{Dom}(\mathcal L)$ and solves $\mathcal L\phi = -f$. By Lemma~\ref{lem: domain}, $\phi$ has a velocity gradient $\nabla_v\phi \in L^2(\mu)$, and \eqref{eq: domain form} for $\varrho = \phi$ is \eqref{eq: weak poisson}. Hence $\phi$ is a weak solution of \eqref{eq: poisson equation}. The identity of \eqref{eq: resolvent bound} for $F = f$ is \eqref{eq: energy}. If $f \in C^\infty(\mathcal X)$, Lemma~\ref{lem: domain} shows that $\phi$ is smooth and satisfies \eqref{eq: poisson equation} on $T\mathcal X$, so $\phi$ is a strong solution.

    To prove \eqref{eq: variance upper}, let $g \in C^\infty(\mathcal X) \cap L^2(\bar\mu)$ with $\mu(g) = 0$, and let $u$ be the solution of \eqref{eq: position equation} given by Lemma~\ref{lem: position}. Then $v \cdot \nabla u \in C^1(T\mathcal X)$ and $\nabla_v ( v \cdot \nabla u ) = \nabla u$. By the definitions of the position gradient \eqref{eq: position gradient} and of $\nabla_v^*$ and by \eqref{eq: overdamped} and \eqref{eq: position equation},
    \begin{equation}
        \mathcal A ( v \cdot \nabla u ) = v \cdot \nabla^2 u\, v - \nabla U \cdot \nabla u = \nabla_v^* \cdot ( \nabla^2 u\, v ) - g .
        \label{eq: transport}
    \end{equation}
    By \eqref{eq: transport} and \eqref{eq: position bounds}, $\psi = v \cdot \nabla u$ satisfies $\psi, \mathcal A\psi, \nabla_v\psi \in L^2(\mu)$ as in Definition~\ref{def: weak solution}, and $W = \nabla^2 u\, v$ satisfies $W, \nabla_v^* \cdot W \in L^2(\mu)$ as in Lemma~\ref{lem: domain}. Since $f$ is a function of the position, $\mu( f\, v \cdot \nabla u ) = 0$. Hence \eqref{eq: weak poisson} for $\psi = v \cdot \nabla u$, \eqref{eq: transport} and \eqref{eq: domain adjoints} for $\varrho = \phi$ and $W = \nabla^2 u\, v$ give
    \begin{equation}
        \mu( \phi\, g ) = \mu\bigl( \nabla_v\phi \cdot ( \nabla^2 u\, v + \gamma \nabla u ) \bigr) .
        \label{eq: variance identity}
    \end{equation}
    By \eqref{eq: gaussian moments} for $a = \gamma \nabla u$ and $A = -\nabla^2 u$ and by \eqref{eq: position bounds},
    \begin{equation}
        \mu\bigl( | \nabla^2 u\, v + \gamma \nabla u |^2 \bigr) = \mu\bigl( |\nabla^2 u|^2 \bigr) + \gamma^2 \mu\bigl( |\nabla u|^2 \bigr) \Le \frac{\gamma^2 + m + K}{m} \mu(g^2) .
        \label{eq: variance position}
    \end{equation}
    By \eqref{eq: average identity} for $\psi = \phi$, the velocity average $\langle\phi\rangle$ of \eqref{eq: velocity average} satisfies $\mu( \phi\, g ) = \mu( \langle\phi\rangle g )$, so the Cauchy--Schwarz inequality in \eqref{eq: variance identity} with \eqref{eq: variance position} gives
    \begin{equation*}
        \mu( \langle\phi\rangle g )^2 \Le \frac{\gamma^2 + m + K}{m} \norm{\nabla_v\phi}_{L^2(\mu)}^2 \mu(g^2) .
    \end{equation*}
    Since $\mu( \langle\phi\rangle ) = \mu(\phi) = 0$ and $\{ g \in C^\infty(\mathcal X) \cap L^2(\bar\mu) : \mu(g) = 0 \}$ is dense in $\{ g \in L^2(\bar\mu) : \mu(g) = 0 \}$,
    \begin{equation}
        \norm{\langle\phi\rangle}_{L^2(\mu)}^2 \Le \frac{\gamma^2 + m + K}{m} \norm{\nabla_v\phi}_{L^2(\mu)}^2 .
        \label{eq: average bound}
    \end{equation}
    By \eqref{eq: energy} and by \eqref{eq: average identity} for $\psi = \phi$, with $f$ in place of $g$, the Cauchy--Schwarz inequality with \eqref{eq: average bound} gives
    \begin{equation*}
        \gamma \norm{\nabla_v\phi}_{L^2(\mu)}^2 = \mu( \langle\phi\rangle f ) \Le \Bigl( \frac{\gamma^2 + m + K}{m} \Bigr)^{\frac12} \norm{\nabla_v\phi}_{L^2(\mu)} \norm{f}_{L^2(\mu)} ,
    \end{equation*}
    which gives \eqref{eq: variance upper} by \eqref{eq: asymptotic variance}.

    To prove \eqref{eq: poisson moments}, note that \eqref{eq: average split} for $\psi = \phi$, \eqref{eq: average bound} and \eqref{eq: gaussian poincare} for $\varrho = \phi$ give
    \begin{equation*}
        \norm{\phi}_{L^2(\mu)}^2 = \norm{\langle\phi\rangle}_{L^2(\mu)}^2 + \norm{\phi - \langle\phi\rangle}_{L^2(\mu)}^2 \Le \frac{\gamma^2 + 2m + K}{m} \norm{\nabla_v\phi}_{L^2(\mu)}^2 ,
    \end{equation*}
    which is \eqref{eq: poisson moments} by \eqref{eq: asymptotic variance}.

    Finally, It\^o's formula \eqref{eq: domain ito} for $\varrho = \phi$ gives \eqref{eq: martingale decomposition} along the stationary dynamics, and hence from $\mu$-almost every initial state, simultaneously for all $t \Ge 0$ after passing to the continuous modification of $t \mapsto \phi(Z_t)$.
\end{proof}

\subsection{Proofs of Corollary~\ref{cor: variance} and of the bounds of Remark~\ref{rem: overdamped}}
\label{app: variance proofs}

We first prove that, for $f$ and $\phi$ of Theorem~\ref{thm: poisson} with $\mu(f) = 0$, every function $g$ with $g, \nabla g \in L^2(\bar\mu)$ satisfies
\begin{equation}
    \mu( f g )^2 \Le \frac{\sigma_f^2}{2\gamma} \mu\bigl( |\nabla g|^2 \bigr) .
    \label{eq: variance lower}
\end{equation}
To prove \eqref{eq: variance lower}, first let $g \in C^1(\mathcal X)$, so that $\nabla_v g = 0$ and
\begin{equation*}
    \mathcal A g = v \cdot \nabla g = \nabla_v^* \cdot \nabla g \in L^2(\mu) .
\end{equation*}
Then \eqref{eq: weak poisson} for $\psi = g$ and \eqref{eq: domain adjoints} for $\varrho = \phi$ and $W = \nabla g$ give
\begin{equation*}
    \mu( f g ) = \mu( \nabla_v\phi \cdot \nabla g ) ,
\end{equation*}
and the Cauchy--Schwarz inequality with \eqref{eq: asymptotic variance} gives \eqref{eq: variance lower}, which extends by density to $g, \nabla g \in L^2(\bar\mu)$.

\begin{proof}[Proof of Corollary~\ref{cor: variance}]
    To prove (i), let $f \in L^4(\bar\mu)$. By \eqref{eq: rate lower},
    \begin{equation*}
        \frac1\lambda \Le \frac{3 ( \gamma + \sqrt m + \sqrt K )^2}{m\gamma} ,
    \end{equation*}
    so \eqref{eq: resolvent fourth} for $F = f$ gives \eqref{eq: poisson fourth}.
    To prove (ii), let $\nabla f \in L^2(\bar\mu)$ and $\mu\bigl( |\nabla f|^2 \bigr) > 0$. Then \eqref{eq: variance lower} for $g = f$, divided by $\mu\bigl( |\nabla f|^2 \bigr)$, gives \eqref{eq: variance below}.
\end{proof}

\begin{proof}[Proof of \eqref{eq: variance overdamped}]
    First let $f \in C^\infty(\mathcal X)$, so that $\bar\phi$ is the solution $u$ of \eqref{eq: position equation} for $g = f$. By \eqref{eq: position weak} for $h = u$, $\mu( f \bar\phi ) = \mu\bigl( |\nabla\bar\phi|^2 \bigr)$, so \eqref{eq: variance lower} for $g = \bar\phi$ gives the lower bound of $\sigma_f^2$ in \eqref{eq: variance overdamped}. By \eqref{eq: energy}, \eqref{eq: variance identity} for $g = f$, the Cauchy--Schwarz inequality and the identity of \eqref{eq: variance position},
    \begin{equation*}
        \gamma \norm{\nabla_v\phi}_{L^2(\mu)}^2 = \mu( \phi f ) \Le \norm{\nabla_v\phi}_{L^2(\mu)} \Bigl( \mu\bigl( |\nabla^2\bar\phi|^2 \bigr) + \gamma^2 \mu\bigl( |\nabla\bar\phi|^2 \bigr) \Bigr)^{\frac12} ,
    \end{equation*}
    which gives the upper bound of $\sigma_f^2$ in \eqref{eq: variance overdamped} by \eqref{eq: asymptotic variance}. The bounds \eqref{eq: resolvent bound} and \eqref{eq: position bounds} extend \eqref{eq: variance overdamped} to $f \in L^2(\bar\mu)$ by density.
\end{proof}

\subsection{Flip of the velocity}
\label{app: flip}

We define the flip of the velocity as the map
\begin{equation}
    \iota(x,v) := (x,-v) , \qquad (x,v) \in T\mathcal X .
    \label{eq: flip}
\end{equation}
In Lemma~\ref{lem: reversal}, we prove the symmetry of the generator and of the semigroup under $\iota$.

\begin{lemma}
\label{lem: reversal}
    For the flip $\iota$ of \eqref{eq: flip}, all $\varrho_1, \varrho_2 \in \mathrm{Dom}(\mathcal L)$ satisfy
    \begin{equation}
        \mu\bigl( ( \mathcal L\varrho_1 )\, ( \varrho_2 \circ \iota ) \bigr) = \mu\bigl( ( \varrho_1 \circ \iota )\, \mathcal L\varrho_2 \bigr) ,
        \label{eq: generator reversal}
    \end{equation}
    and all $\psi_1, \psi_2 \in L^2(\mu)$ and every $t \Ge 0$ satisfy
    \begin{equation}
        \mu\bigl( ( e^{t\mathcal L}\psi_1 )\, ( \psi_2 \circ \iota ) \bigr) = \mu\bigl( ( \psi_1 \circ \iota )\, e^{t\mathcal L}\psi_2 \bigr) .
        \label{eq: reversal}
    \end{equation}
\end{lemma}

\begin{proof}[Proof of Lemma~\ref{lem: reversal}]
    To prove \eqref{eq: generator reversal}, first let $\varrho_1, \varrho_2 \in C_c^\infty(T\mathcal X)$. Every $\varrho \in C_c^\infty(T\mathcal X)$ satisfies
    \begin{equation}
        \mathcal A ( \varrho \circ \iota ) = -( \mathcal A\varrho ) \circ \iota , \qquad
        \nabla_v ( \varrho \circ \iota ) = -( \nabla_v\varrho ) \circ \iota ,
        \label{eq: flip derivatives}
    \end{equation}
    so \eqref{eq: generator form} for $\psi_1 = \varrho_1$ and $\psi_2 = \varrho_2 \circ \iota$ reads
    \begin{equation}
        \mu\bigl( ( \mathcal L\varrho_1 )\, ( \varrho_2 \circ \iota ) \bigr) = \mu\bigl( \varrho_1\, ( ( \mathcal A\varrho_2 ) \circ \iota ) \bigr) + \gamma \mu\bigl( \nabla_v\varrho_1 \cdot ( ( \nabla_v\varrho_2 ) \circ \iota ) \bigr) .
        \label{eq: flip form}
    \end{equation}
    By \eqref{eq: flip derivatives}, the antisymmetry of $\mathcal A$ and the invariance of $\mu$ under $\iota$, the right side of \eqref{eq: flip form} is unchanged when $\varrho_1$ and $\varrho_2$ are exchanged. Hence the left side of \eqref{eq: flip form} is unchanged by this exchange, which is \eqref{eq: generator reversal}. Since $C_c^\infty(T\mathcal X)$ is a core of the generator by Lemma~\ref{lem: dynamics}, \eqref{eq: generator reversal} extends to all $\varrho_1, \varrho_2 \in \mathrm{Dom}(\mathcal L)$.

    To prove \eqref{eq: reversal}, first let $\psi_1, \psi_2 \in \mathrm{Dom}(\mathcal L)$ and $t \Ge 0$. Since $\iota$ preserves $\mu$, \eqref{eq: generator reversal} for $\varrho_1 = e^{s\mathcal L}\psi_1$ and $\varrho_2 = e^{(t-s)\mathcal L}\psi_2$ gives
    \begin{equation}
        \frac{\D}{\D s}\, \mu\bigl( ( e^{s\mathcal L}\psi_1 )\, ( ( e^{(t-s)\mathcal L}\psi_2 ) \circ \iota ) \bigr) = 0 , \qquad 0 \Le s \Le t ,
        \label{eq: reversal derivative}
    \end{equation}
    and the values of the mean in \eqref{eq: reversal derivative} at $s = t$ and $s = 0$ are the two sides of \eqref{eq: reversal}. By density, \eqref{eq: reversal} extends to all $\psi_1, \psi_2 \in L^2(\mu)$.
\end{proof}

\subsection{Proof of Corollary~\ref{cor: second moment}}
\label{app: second moment}

Finally, we prove Corollary~\ref{cor: second moment} by the identity \eqref{eq: semigroup integral} and the symmetry \eqref{eq: reversal}.

\begin{proof}[Proof of Corollary~\ref{cor: second moment}]
    Let $\phi$ be the Poisson solution \eqref{eq: poisson solution} of $f$. Since $Z_0 \sim \mu$, the Markov property gives
    \begin{equation*}
        \mathbb E\bigl[ f(X_{s'})\, f(X_s) \bigr] = \mu\bigl( f\, e^{(s-s')\mathcal L}f \bigr) , \qquad 0 \Le s' \Le s ,
    \end{equation*}
    so \eqref{eq: semigroup integral} for $F = f$ gives
    \begin{equation}
        \mathbb E\biggl[ \Bigl( \int_0^t f(X_s) \,\D s \Bigr)^2 \biggr] = 2 \int_0^t \mu\Bigl( f \int_0^s e^{r\mathcal L}f \,\D r \Bigr) \D s = 2 \int_0^t \mu\bigl( f\, ( \phi - e^{s\mathcal L}\phi ) \bigr) \D s .
        \label{eq: second moment double}
    \end{equation}
    By \eqref{eq: energy} and \eqref{eq: asymptotic variance}, $2 \mu( f \phi ) = \sigma_f^2$. Since $f \circ \iota = f$ for the flip $\iota$ of \eqref{eq: flip}, \eqref{eq: reversal} at the time $s$ for $\psi_1 = f$ and $\psi_2 = \phi$ gives
    \begin{equation}
        \mu\bigl( f\, e^{s\mathcal L}\phi \bigr) = \mu\bigl( ( e^{s\mathcal L}f )\, ( \phi \circ \iota ) \bigr) , \qquad \forall s \Ge 0 .
        \label{eq: second moment flip}
    \end{equation}
    Hence \eqref{eq: second moment double}, $2 \mu( f \phi ) = \sigma_f^2$, \eqref{eq: second moment flip} and \eqref{eq: semigroup integral} for $F = f$ and $s = t$ give
    \begin{equation}
        \mathbb E\biggl[ \Bigl( \int_0^t f(X_s) \,\D s \Bigr)^2 \biggr] = \sigma_f^2 t - 2 \mu\bigl( ( \phi - e^{t\mathcal L}\phi )\, ( \phi \circ \iota ) \bigr) .
        \label{eq: second moment integral}
    \end{equation}
    Since $\iota$ preserves $\mu$ and the semigroup contracts $L^2(\mu)$,
    \begin{equation*}
        \mu\bigl( \phi ( \phi \circ \iota ) \bigr) = \norm{\phi}_{L^2(\mu)}^2 - \frac12 \norm{\phi - \phi \circ \iota}_{L^2(\mu)}^2 \Le \norm{\phi}_{L^2(\mu)}^2 , \qquad
        \bigl| \mu\bigl( ( e^{t\mathcal L}\phi ) ( \phi \circ \iota ) \bigr) \bigr| \Le \norm{\phi}_{L^2(\mu)}^2 ,
    \end{equation*}
    so \eqref{eq: second moment integral} gives
    \begin{equation}
        -4 \norm{\phi}_{L^2(\mu)}^2 \Le \mathbb E\biggl[ \Bigl( \int_0^t f(X_s) \,\D s \Bigr)^2 \biggr] - \sigma_f^2 t \Le \norm{\phi - \phi \circ \iota}_{L^2(\mu)}^2 .
        \label{eq: second moment gap}
    \end{equation}
    After division by $t^2$, \eqref{eq: second moment gap} and \eqref{eq: poisson moments} give the lower bound of \eqref{eq: second moment}. Since the velocity average $\langle\phi\rangle$ satisfies $\langle\phi\rangle \circ \iota = \langle\phi\rangle$, \eqref{eq: gaussian poincare} for $\varrho = \phi$ and \eqref{eq: asymptotic variance} give
    \begin{equation*}
        \norm{\phi - \phi \circ \iota}_{L^2(\mu)}^2 \Le 4 \norm{\phi - \langle\phi\rangle}_{L^2(\mu)}^2 \Le 4 \norm{\nabla_v\phi}_{L^2(\mu)}^2 = \frac2\gamma \sigma_f^2 ,
    \end{equation*}
    which gives the upper bound of \eqref{eq: second moment} by \eqref{eq: second moment gap}.
\end{proof}

\section{Proofs for the properties of the normalized variance}
\label{app: variance}

This appendix proves Theorem~\ref{thm: properties} and Lemmas~\ref{lem: quadratic} and~\ref{lem: torus}. The ratio of \eqref{eq: normalized variance} is unchanged when a constant is added to $f$, so the proofs of this appendix assume $\mu(f) = 0$. By Lemma~\ref{lem: poisson} for $F = f$, the Poisson solution $\phi$ of $f \in L^2(\bar\mu)$ with $\mu(f) = 0$ belongs to $\mathrm{Dom}(\mathcal L)$ and satisfies the identity \eqref{eq: energy}. Hence \eqref{eq: asymptotic variance}, \eqref{eq: poisson solution} and the Markov property with $Z_0 \sim \mu$ give
\begin{equation}
    \sigma_f^2 = 2 \mu( \phi f ) = 2 \int_0^\infty \mu\bigl( f e^{t\mathcal L}f \bigr) \D t = 2 \int_0^\infty \mathbb E\bigl[ f(X_0) f(X_t) \bigr] \D t .
    \label{eq: variance correlation}
\end{equation}
The identity matrix of $\mathbb R^d$ is $I_d$.

\subsection{Proof of Theorem~\ref{thm: properties}}
\label{app: properties}

We prove (i) of Theorem~\ref{thm: properties} by \eqref{eq: generator form} for two damping rates, and (ii) by \eqref{eq: weak poisson} for a function of $U(x) + \frac12 |v|^2$.

\begin{proof}[Proof of Theorem~\ref{thm: properties}]
    To prove (i), let $0 < \gamma_1 \Le \gamma_2$, and for $i = 1, 2$, let $\mathcal L_i$ be the generator of the dynamics with the damping rate $\gamma_i$. Let $\varrho_1, \varrho_2 \in C_c^\infty(T\mathcal X)$. Since $\mathcal A$ is antisymmetric, the sum of \eqref{eq: generator form} for the damping rate $\gamma_1$ with $( \psi_1, \psi_2 ) = ( \varrho_1, \varrho_2 )$ and for the damping rate $\gamma_2$ with $( \psi_1, \psi_2 ) = ( \varrho_2, \varrho_1 )$ is
    \begin{equation}
        \mu\bigl( ( -\mathcal L_1\varrho_1 ) \varrho_2 \bigr) + \mu\bigl( ( -\mathcal L_2\varrho_2 ) \varrho_1 \bigr) = ( \gamma_1 + \gamma_2 ) \mu( \nabla_v\varrho_1 \cdot \nabla_v\varrho_2 ) .
        \label{eq: two rates}
    \end{equation}
    By Lemma~\ref{lem: dynamics} and \eqref{eq: dissipation}, \eqref{eq: two rates} extends to $\varrho_1 \in \mathrm{Dom}(\mathcal L_1)$ and $\varrho_2 \in \mathrm{Dom}(\mathcal L_2)$. Let $f \in L^2(\bar\mu)$ be nonconstant, and for $i = 1, 2$, let $\phi_i$ and $\sigma_i^2$, with $\sigma_i \Ge 0$, be the Poisson solution and the asymptotic variance of $f$ for the damping rate $\gamma_i$. By Lemma~\ref{lem: poisson}, \eqref{eq: asymptotic variance} and \eqref{eq: variance correlation}, $\phi_i \in \mathrm{Dom}(\mathcal L_i)$ and
    \begin{equation*}
        \mathcal L_i\phi_i = -f , \qquad
        \sigma_i^2 = 2 \gamma_i \norm{\nabla_v\phi_i}_{L^2(\mu)}^2 = 2 \mu( f \phi_i ) .
    \end{equation*}
    Hence \eqref{eq: two rates} for $( \varrho_1, \varrho_2 ) = ( \phi_1, \phi_2 )$ and the Cauchy--Schwarz inequality give
    \begin{equation*}
        \frac{\sigma_1^2 + \sigma_2^2}{2} \Le ( \gamma_1 + \gamma_2 ) \norm{\nabla_v\phi_1}_{L^2(\mu)} \norm{\nabla_v\phi_2}_{L^2(\mu)} = \frac{\gamma_1 + \gamma_2}{2 \sqrt{\gamma_1\gamma_2}} \sigma_1 \sigma_2 ,
    \end{equation*}
    that is,
    \begin{equation*}
        \bigl( \sqrt{\gamma_2} \sigma_2 - \sqrt{\gamma_1} \sigma_1 \bigr) \bigl( \sqrt{\gamma_2} \sigma_1 - \sqrt{\gamma_1} \sigma_2 \bigr) \Ge 0 .
    \end{equation*}
    Since the sum $( \sqrt{\gamma_2} - \sqrt{\gamma_1} ) ( \sigma_1 + \sigma_2 )$ of the two factors is nonnegative, both factors are nonnegative, that is,
    \begin{equation*}
        \gamma_1 \sigma_1^2 \Le \gamma_2 \sigma_2^2 , \qquad
        \frac{\sigma_2^2}{\gamma_2} \Le \frac{\sigma_1^2}{\gamma_1} .
    \end{equation*}
    Dividing these two inequalities by $\mu(f^2)$ and taking the supremum over $f$ gives (i).

    To prove (ii), let $U$ be nonconstant, and let $\psi := \exp\bigl( -U - \frac12 |v|^2 \bigr)$, which is bounded since $U$ is bounded from below. The function $\psi \in C^1(T\mathcal X)$ satisfies $\mathcal A\psi = 0$ and $\nabla_v\psi = -\psi v$. Its velocity average \eqref{eq: velocity average} is $\langle\psi\rangle = 2^{-\frac d2} e^{-U}$, which is not constant. Let $f := \langle\psi\rangle - \mu(\psi)$, and let $\phi$ be its Poisson solution. By \eqref{eq: average identity} for $g = f$ and \eqref{eq: weak poisson} for $\psi$,
    \begin{equation*}
        \mu(f^2) = \mu( f \psi ) = -\gamma \mu( \psi\, \nabla_v\phi \cdot v ) ,
    \end{equation*}
    and the Cauchy--Schwarz inequality with \eqref{eq: asymptotic variance} gives
    \begin{equation*}
        \mu(f^2)^2 \Le \frac{\gamma \sigma_f^2}{2} \mu\bigl( |v|^2 \psi^2 \bigr) .
    \end{equation*}
    Since $\psi$ and $f$ do not depend on $\gamma$, every $\gamma > 0$ satisfies
    \begin{equation*}
        \mathbb V(\gamma) \Ge \frac{\sigma_f^2}{\mu(f^2)} \Ge \frac c\gamma , \qquad
        c := \frac{2 \mu(f^2)}{\mu\bigl( |v|^2 \psi^2 \bigr)} > 0 .
    \end{equation*}
\end{proof}

\subsection{Normalized variance for the quadratic potential}
\label{app: quadratic}

We prove Lemma~\ref{lem: quadratic}. The lower bound of $\mathbb V(\gamma)$ follows from \eqref{eq: normalized bounds} and \eqref{eq: weak poisson} for a quadratic function of the position and the velocity, and the upper bound from \eqref{eq: variance correlation} and the expansion of $f$ in the Hermite polynomials.

\begin{proof}[Proof of Lemma~\ref{lem: quadratic}]
    For $\mathcal X = \mathbb R^d$ and $U(x) = \frac m2 |x|^2$, Assumptions~\ref{ass: manifold}, \ref{ass: curvature} and~\ref{ass: poincare} hold with $K = 0$ and the spectral gap $m$ of $\bar{\mathcal L}$, and $\bar\mu = \mathcal N(0, m^{-1} I_d)$. To prove
    \begin{equation}
        \mathbb V(\gamma) \Ge \max\Bigl( \frac{2\gamma}{m} , \frac\gamma m + \frac1\gamma \Bigr) ,
        \label{eq: quadratic below}
    \end{equation}
    let $f \in L^2(\bar\mu)$, and let $\phi$ be its Poisson solution \eqref{eq: poisson solution}. Let $x_1$ and $v_1$ be the first coordinates of the position and of the velocity, and let $\partial_{v_1}\phi$ be the first coordinate of $\nabla_v\phi$. Since $\nabla U(x) = m x$,
    \begin{equation*}
        \mathcal A \Bigl( x_1^2 + \frac{v_1^2}{m + \gamma^2} \Bigr) = \frac{2\gamma^2}{m + \gamma^2} x_1 v_1 .
    \end{equation*}
    The field $W = ( x_1, 0, \dots, 0 )$ has $\nabla_v^* \cdot W = x_1 v_1$, so \eqref{eq: domain adjoints} for $\varrho = \phi$ gives $\mu( \phi\, x_1 v_1 ) = \mu( x_1\, \partial_{v_1}\phi )$. Hence \eqref{eq: weak poisson} for $\psi = x_1^2 + \frac{v_1^2}{m + \gamma^2}$ gives
    \begin{equation*}
        \mu\Bigl( f \Bigl( x_1^2 + \frac{v_1^2}{m + \gamma^2} \Bigr) \Bigr) = \frac{2\gamma}{m + \gamma^2} \mu\bigl( ( \gamma x_1 + v_1 )\, \partial_{v_1}\phi \bigr) .
    \end{equation*}
    Since $\mu\bigl( ( \gamma x_1 + v_1 )^2 \bigr) = \frac{m + \gamma^2}{m}$, the Cauchy--Schwarz inequality with \eqref{eq: asymptotic variance} gives
    \begin{equation}
        \mu\Bigl( f \Bigl( x_1^2 + \frac{v_1^2}{m + \gamma^2} \Bigr) \Bigr)^2 \Le \frac{2\gamma}{m ( m + \gamma^2 )} \sigma_f^2 .
        \label{eq: quadratic lower}
    \end{equation}

    The function $f(x) = x_1^2 - m^{-1}$ satisfies
    \begin{equation*}
        \mu( f v_1^2 ) = 0 , \qquad
        \mu( f x_1^2 ) = \mu(f^2) = \frac{2}{m^2} ,
    \end{equation*}
    so \eqref{eq: quadratic lower} gives
    \begin{equation*}
        \mathbb V(\gamma) \Ge \frac{\sigma_f^2}{\mu(f^2)} \Ge \frac{m ( m + \gamma^2 )}{2\gamma} \mu(f^2) = \frac\gamma m + \frac1\gamma ,
    \end{equation*}
    and the lower bound of \eqref{eq: normalized bounds} gives $\mathbb V(\gamma) \Ge \frac{2\gamma}{m}$, which proves \eqref{eq: quadratic below}.

    To prove the reverse inequality of \eqref{eq: quadratic below}, let $\mathfrak s$ be the solution of
    \begin{equation}
        \mathfrak s'' + \gamma \mathfrak s' + m \mathfrak s = 0 , \qquad \mathfrak s(0) = 1 , \quad \mathfrak s'(0) = 0 .
        \label{eq: correlation}
    \end{equation}
    Since \eqref{eq: uld} is linear, the dynamics started at $z = (x,v)$ satisfies
    \begin{equation*}
        \mathbb E^z[X_t] = \mathfrak s(t) x - \frac{\mathfrak s'(t)}{m} v .
    \end{equation*}
    Since $Z_0 \sim \mu$, the vector $(X_0, X_t)$ is centered Gaussian in $\mathbb R^{2d}$ with $X_0, X_t \sim \mathcal N(0, m^{-1} I_d)$, and $X_0$ is independent of $V_0$. Hence, for $i, j = 1, \dots, d$ and the coordinates $X_0^i$ and $X_t^j$ of $X_0$ and $X_t$,
    \begin{equation}
        \mathbb E\bigl[ X_0^i X_t^j \bigr] = \frac{\mathfrak s(t)}{m} \mathbf 1_{i = j} , \qquad |\mathfrak s(t)| \Le 1 .
        \label{eq: position covariance}
    \end{equation}
    Since $\mathfrak s$ and $\mathfrak s'$ decay exponentially, multiplying \eqref{eq: correlation} by $1$, $\mathfrak s'$ and $\mathfrak s$ and integrating over $(0,\infty)$ gives
    \begin{equation}
        \int_0^\infty \mathfrak s(t) \,\D t = \frac\gamma m , \qquad
        \int_0^\infty \mathfrak s'(t)^2 \,\D t = \frac{m}{2\gamma} , \qquad
        \int_0^\infty \mathfrak s(t)^2 \,\D t = \frac{\gamma}{2m} + \frac{1}{2\gamma} .
        \label{eq: correlation integrals}
    \end{equation}

    Let $f \in L^2(\bar\mu)$ be nonconstant, and for $n \Ge 1$, let $f_n$ be the orthogonal projection of $f$ in $L^2(\bar\mu)$ onto the span of the Hermite polynomials of $\sqrt m\, x$ of total degree $n$ \citep[Section~2.7.1]{bakry2014analysis}. Since $\mu(f) = 0$,
    \begin{equation*}
        f = \sum_{n \Ge 1} f_n , \qquad
        \mu(f^2) = \sum_{n \Ge 1} \mu(f_n^2) .
    \end{equation*}
    The formula for the covariance of the Hermite polynomials of two jointly Gaussian variables \citep[Equation~(1.3)]{nourdin2012lectures}, with \eqref{eq: position covariance}, gives
    \begin{equation}
        \mathbb E\bigl[ f(X_0) f(X_t) \bigr] = \sum_{n \Ge 1} \mathfrak s(t)^n \mu(f_n^2) \Le \mathfrak s(t) \mu(f_1^2) + \mathfrak s(t)^2 \sum_{n \Ge 2} \mu(f_n^2) .
        \label{eq: quadratic correlation}
    \end{equation}
    Then \eqref{eq: variance correlation}, \eqref{eq: quadratic correlation} and \eqref{eq: correlation integrals} give
    \begin{equation*}
        \sigma_f^2 \Le \frac{2\gamma}{m} \mu(f_1^2) + \Bigl( \frac\gamma m + \frac1\gamma \Bigr) \sum_{n \Ge 2} \mu(f_n^2) \Le \max\Bigl( \frac{2\gamma}{m} , \frac\gamma m + \frac1\gamma \Bigr) \mu(f^2) .
    \end{equation*}
    The supremum over $f$ and \eqref{eq: quadratic below} give \eqref{eq: variance quadratic}.
    To prove \eqref{eq: optimal damping}, note that the right side of \eqref{eq: variance quadratic} equals $\frac{2}{\sqrt m}$ at $\gamma = \sqrt m$ and is larger for every other $\gamma > 0$.
\end{proof}

\subsection{Normalized variance for the flat torus}
\label{app: torus}

We prove Lemma~\ref{lem: torus}. The Fourier series of $f$ and \eqref{eq: variance correlation} give \eqref{eq: torus formula}. Two bounds of $s - 1 + e^{-s}$ from above give the lower bound of \eqref{eq: variance torus}, and a change of variable in the integral of \eqref{eq: torus formula} gives the upper bound.

\begin{proof}[Proof of Lemma~\ref{lem: torus}]
    For $\mathcal X = \bigl( \mathbb R / \frac{2\pi}{\sqrt m} \mathbb Z \bigr)^d$ and $U = 0$, Assumptions~\ref{ass: manifold}, \ref{ass: curvature} and~\ref{ass: poincare} hold with $K = 0$ and the spectral gap $m$ of $\bar{\mathcal L} = \Delta$, and $\bar\mu$ is the uniform distribution of $\mathcal X$. In the coordinates of $\mathbb R^d$, \eqref{eq: uld} reads
    \begin{equation*}
        \D X_t = V_t \,\D t , \qquad
        \D V_t = -\gamma V_t \,\D t + \sqrt{2\gamma}\,\D B_t .
    \end{equation*}
    Since $Z_0 \sim \mu$, the velocity $(V_t)_{t \Ge 0}$ is a stationary Ornstein--Uhlenbeck process independent of $X_0$, and the covariance of $V_s$ and $V_{s'}$ is $e^{-\gamma |s-s'|} I_d$. Hence the increment of the position is independent of $X_0$, with
    \begin{equation}
        X_t - X_0 = \int_0^t V_s \,\D s \sim \mathcal N\Bigl( 0 , \frac{2}{\gamma^2} \bigl( \gamma t - 1 + e^{-\gamma t} \bigr) I_d \Bigr) .
        \label{eq: torus increment}
    \end{equation}

    To prove \eqref{eq: torus formula}, let $f \in L^2(\bar\mu)$ be nonconstant, with the Fourier coefficients
    \begin{equation*}
        \hat f(k) := \mu\bigl( f\, e^{-\mathrm i \sqrt m\, k \cdot x} \bigr) , \qquad k \in \mathbb Z^d ,
    \end{equation*}
    where $\mathrm i$ is the imaginary unit, and $\hat f(0) = \mu(f) = 0$. Parseval's identity gives, for every $y \in \mathbb R^d$,
    \begin{equation}
        \mu\bigl( f\, f(\cdot + y) \bigr) = \int_{\mathcal X} f(x)\, f(x + y)\, \bar\mu(\D x) = \sum_{k \neq 0} |\hat f(k)|^2 e^{\mathrm i \sqrt m\, k \cdot y} .
        \label{eq: torus parseval}
    \end{equation}
    Since $X_t - X_0$ is independent of $X_0 \sim \bar\mu$, the expectation of \eqref{eq: torus parseval} at $y = X_t - X_0$ and the characteristic function of the distribution \eqref{eq: torus increment} give
    \begin{equation}
        \mathbb E\bigl[ f(X_0) f(X_t) \bigr] = \sum_{k \neq 0} |\hat f(k)|^2 \exp\Bigl( -\frac{m |k|^2}{\gamma^2} \bigl( \gamma t - 1 + e^{-\gamma t} \bigr) \Bigr) .
        \label{eq: torus correlation}
    \end{equation}
    Since $|k|^2 \Ge 1$ for $k \neq 0$, \eqref{eq: variance correlation}, \eqref{eq: torus correlation} and \eqref{eq: torus parseval} at $y = 0$ give
    \begin{equation*}
        \sigma_f^2 \Le 2 \mu(f^2) \int_0^\infty \exp\Bigl( -\frac{m}{\gamma^2} \bigl( \gamma t - 1 + e^{-\gamma t} \bigr) \Bigr) \D t .
    \end{equation*}
    Equality holds for $f(x) = \cos( \sqrt m\, x_1 )$ with the first coordinate $x_1$ of $x$, since $\hat f(k) = 0$ unless $|k| = 1$. The supremum over $f$ and the change of variable $s = \gamma t$ give \eqref{eq: torus formula}.

    To prove the lower bound of \eqref{eq: variance torus}, note that
    \begin{equation*}
        s - 1 + e^{-s} \Le \min\Bigl( \frac{s^2}{2} , s \Bigr) , \qquad \forall s \Ge 0 ,
    \end{equation*}
    so
    \begin{equation*}
        \int_0^\infty \exp\Bigl( -\frac{m}{\gamma^2} \bigl( s - 1 + e^{-s} \bigr) \Bigr) \D s \Ge \max\Bigl( \frac{\gamma}{\sqrt m} \sqrt{\frac\pi2} , \frac{\gamma^2}{m} \Bigr) ,
    \end{equation*}
    and \eqref{eq: torus formula} then gives the lower bound of \eqref{eq: variance torus}.

    To prove the upper bound of \eqref{eq: variance torus}, let $r := \sqrt{2 ( s - 1 + e^{-s} )}$ for $s \Ge 0$. Since $r \Le s \Le e^s - 1$,
    \begin{equation*}
        \frac{\D r}{\D s} = \frac{1 - e^{-s}}{r} \Ge \frac{1}{1 + r} , \qquad \forall s > 0 .
    \end{equation*}
    Hence the change of variable from $s$ to $r$ gives
    \begin{equation*}
        \int_0^\infty \exp\Bigl( -\frac{m}{\gamma^2} \bigl( s - 1 + e^{-s} \bigr) \Bigr) \D s \Le \int_0^\infty e^{-\frac{m r^2}{2\gamma^2}} ( 1 + r ) \,\D r = \frac{\gamma}{\sqrt m} \sqrt{\frac\pi2} + \frac{\gamma^2}{m} ,
    \end{equation*}
    and \eqref{eq: torus formula} then gives the upper bound of \eqref{eq: variance torus}.
\end{proof}

\section{Proofs for the relation to the mixing time}
\label{app: mixing}

This appendix proves Theorem~\ref{thm: mixing}.

\subsection{Resolvent of the generator on the imaginary axis}
\label{app: axis}

In Lemma~\ref{lem: axis}, we bound the resolvent of $\mathcal L$ on the imaginary axis and the Poisson solution \eqref{eq: resolvent}. In this subsection, the functions are complex valued, and $\mathrm i$ is the imaginary unit. The complex conjugate of such a function $\psi$ is $\overline\psi$, and $\norm{\psi}_{L^2(\mu)}^2 = \mu( \psi \overline\psi )$. The semigroup, the generator and the velocity average $\langle\psi\rangle$ of \eqref{eq: velocity average} act on the real and imaginary parts of $\psi$.

For a real number $\omega$ and $\zeta \in L^2(\mu)$ with $\mu(\zeta) = 0$, the resolvent of $\mathcal L$ at $\mathrm i\omega$ maps $\zeta$ to
\begin{equation}
    ( \mathrm i\omega - \mathcal L )^{-1} \zeta := \int_0^\infty e^{-\mathrm i\omega t} e^{t\mathcal L} \zeta \,\D t .
    \label{eq: axis resolvent}
\end{equation}
As in the proof of Lemma~\ref{lem: poisson}, the decay \eqref{eq: l2 decay} gives the convergence of this integral in $L^2(\mu)$, and the function $\varphi := ( \mathrm i\omega - \mathcal L )^{-1} \zeta$ satisfies $\mu(\varphi) = 0$ and $\mathrm{Re}\, \varphi, \mathrm{Im}\, \varphi \in \mathrm{Dom}(\mathcal L)$ and solves
\begin{equation}
    \mathrm i\omega \varphi - \mathcal L \varphi = \zeta .
    \label{eq: axis equation}
\end{equation}
The function $\varphi$ is the only solution of \eqref{eq: axis equation} with $\mu(\varphi) = 0$: a function $\varrho$ with $\mu(\varrho) = 0$ and $\mathcal L\varrho = \mathrm i\omega\varrho$ satisfies $e^{t\mathcal L}\varrho = e^{\mathrm i\omega t}\varrho$ for every $t \Ge 0$, so the decay \eqref{eq: l2 decay} gives $\varrho = 0$. For $\omega = 0$, the integral \eqref{eq: axis resolvent} is the Poisson solution \eqref{eq: resolvent} of $\zeta$, written $( -\mathcal L )^{-1} \zeta$.

\begin{lemma}
\label{lem: axis}
    For every $\zeta \in L^2(\mu)$ with $\mu(\zeta) = 0$, the Poisson solution $( -\mathcal L )^{-1} \zeta$ of \eqref{eq: resolvent} satisfies, with the normalized variance $\mathbb V(\gamma)$ of \eqref{eq: normalized variance},
    \begin{equation}
        \norm{( -\mathcal L )^{-1} \zeta}_{L^2(\mu)} \Le \Bigl( \frac12 \mathbb V(\gamma) + \frac1\gamma \Bigr) \norm{\zeta}_{L^2(\mu)} .
        \label{eq: axis poisson}
    \end{equation}
    Moreover, for every real number $\omega$, the function $( \mathrm i\omega - \mathcal L )^{-1} \zeta$ of \eqref{eq: axis resolvent} satisfies
    \begin{equation}
        \norm{( \mathrm i\omega - \mathcal L )^{-1} \zeta}_{L^2(\mu)} \Le \Bigl( \mathbb V(\gamma) + \frac3\gamma \Bigr) \norm{\zeta}_{L^2(\mu)} .
        \label{eq: axis bound}
    \end{equation}
\end{lemma}

\begin{proof}[Proof of Lemma~\ref{lem: axis}]
    Let $\omega$ be a real number, $\zeta \in L^2(\mu)$ with $\mu(\zeta) = 0$ and $\varphi = ( \mathrm i\omega - \mathcal L )^{-1} \zeta$. We multiply \eqref{eq: axis equation} by $\overline\varphi$ and take the real part of the mean. Then \eqref{eq: dissipation} and \eqref{eq: gaussian poincare} for $\mathrm{Re}\, \varphi$ and $\mathrm{Im}\, \varphi$ give
    \begin{equation}
        \mathrm{Re}\big\{ \mu( \zeta \overline\varphi ) \big\} \Ge \gamma \norm{\varphi - \langle \varphi\rangle}_{L^2(\mu)}^2 .
        \label{eq: axis dissipation}
    \end{equation}
    Let $\mathfrak k(\omega)$ be the norm of the map $f \mapsto \langle ( \mathrm i\omega - \mathcal L )^{-1} f \rangle$ on the functions of the position,
    \begin{equation}
        \mathfrak k(\omega) := \sup_{f} \frac{\norm{\langle ( \mathrm i\omega - \mathcal L )^{-1} f \rangle}_{L^2(\mu)}}{\norm{f}_{L^2(\mu)}} ,
        \label{eq: axis position}
    \end{equation}
    where the supremum is over the functions $f \neq 0$ in $L^2(\bar\mu)$ with $\mu(f) = 0$. For $\zeta = f$ in the supremum \eqref{eq: axis position}, \eqref{eq: average identity} gives $\mu( f \overline\varphi ) = \mu\bigl( f\, \overline{\langle \varphi\rangle} \bigr)$, so \eqref{eq: axis dissipation}, the Cauchy--Schwarz inequality and \eqref{eq: axis position} give
    \begin{equation}
        \norm{\varphi - \langle \varphi\rangle}_{L^2(\mu)}^2 \Le \frac1\gamma \mathrm{Re}\big\{ \mu\bigl( f\, \overline{\langle \varphi\rangle} \bigr) \big\} \Le \frac{\mathfrak k(\omega)}{\gamma} \norm{f}_{L^2(\mu)}^2 .
        \label{eq: axis cross}
    \end{equation}

    To prove the bound \eqref{eq: axis blocks} of $\norm{( \mathrm i\omega - \mathcal L )^{-1} \zeta}_{L^2(\mu)}$ by $\mathfrak k(\omega)$, let $\varphi_1$ and $\varphi_2$ be the functions
    \begin{equation*}
        \varphi_1 :=( \mathrm i\omega - \mathcal L )^{-1} \langle \zeta\rangle , \qquad
        \varphi_2 := ( \mathrm i\omega - \mathcal L )^{-1} ( \zeta - \langle \zeta\rangle ) .
    \end{equation*}
    By \eqref{eq: axis position} and \eqref{eq: axis cross} for $\langle \zeta\rangle$ and $\varphi_1$ in place of $f$ and $\varphi$,
    \begin{equation}
        \norm{\langle \varphi_1\rangle}_{L^2(\mu)} \Le \mathfrak k(\omega) \norm{\langle \zeta\rangle}_{L^2(\mu)} , \qquad
        \norm{\varphi_1 - \langle \varphi_1\rangle}_{L^2(\mu)} \Le \Bigl( \frac{\mathfrak k(\omega)}{\gamma} \Bigr)^{\frac12} \norm{\langle \zeta\rangle}_{L^2(\mu)} .
        \label{eq: axis position source}
    \end{equation}
    Since $\langle \zeta - \langle \zeta\rangle \rangle = 0$, \eqref{eq: axis dissipation} for $\zeta - \langle \zeta\rangle$ and $\varphi_2$ in place of $\zeta$ and $\varphi$ gives
    \begin{equation}
        \gamma \norm{\varphi_2 - \langle \varphi_2\rangle}_{L^2(\mu)}^2 \Le \mathrm{Re}\big\{ \mu\bigl( ( \zeta - \langle \zeta\rangle )\, \overline{( \varphi_2 - \langle \varphi_2\rangle )} \bigr) \big\} .
        \label{eq: axis velocity part}
    \end{equation}
    Let $f \in L^2(\bar\mu)$ with $\mu(f) = 0$. Since $f \circ \iota = f$ for the flip $\iota$ of \eqref{eq: flip} and $\langle ( \zeta - \langle \zeta\rangle ) \circ \iota \rangle = 0$, multiplying \eqref{eq: reversal} for $\psi_1 = f$ and $\psi_2 = \zeta - \langle \zeta\rangle$ by $e^{-\mathrm i\omega t}$ and integrating over $t \Ge 0$ gives
    \begin{equation*}
        \mu( f \varphi_2 ) = \mu\Bigl( \bigl( ( \mathrm i\omega - \mathcal L )^{-1} f - \langle ( \mathrm i\omega - \mathcal L )^{-1} f \rangle \bigr)\, ( ( \zeta - \langle \zeta\rangle ) \circ \iota ) \Bigr) .
    \end{equation*}
    Hence, for $f = \overline{\langle \varphi_2\rangle}$, the Cauchy--Schwarz inequality, \eqref{eq: axis cross} for $( \mathrm i\omega - \mathcal L )^{-1} f$ in place of $\varphi$ and \eqref{eq: axis velocity part} give
    \begin{equation}
        \norm{\langle \varphi_2\rangle}_{L^2(\mu)} \Le \Bigl( \frac{\mathfrak k(\omega)}{\gamma} \Bigr)^{\frac12} \norm{\zeta - \langle \zeta\rangle}_{L^2(\mu)} , \qquad
        \norm{\varphi_2 - \langle \varphi_2\rangle}_{L^2(\mu)} \Le \frac1\gamma \norm{\zeta - \langle \zeta\rangle}_{L^2(\mu)} .
        \label{eq: axis velocity source}
    \end{equation}
    By \eqref{eq: average split} for $\psi = \varphi_1 + \varphi_2 = ( \mathrm i\omega - \mathcal L )^{-1} \zeta$ and for $\psi = \zeta$, the bounds \eqref{eq: axis position source} and \eqref{eq: axis velocity source} and the Cauchy--Schwarz inequality give
    \begin{align}
        \norm{( \mathrm i\omega - \mathcal L )^{-1} \zeta}_{L^2(\mu)}^2 & \Le \Bigl( \mathfrak k(\omega) + \frac1\gamma \Bigr) \Bigl( \mathfrak k(\omega)^{\frac12} \norm{\langle \zeta\rangle}_{L^2(\mu)} + \gamma^{-\frac12} \norm{\zeta - \langle \zeta\rangle}_{L^2(\mu)} \Bigr)^2 \notag \\
        & \Le \Bigl( \mathfrak k(\omega) + \frac1\gamma \Bigr)^2 \norm{\zeta}_{L^2(\mu)}^2 .
        \label{eq: axis blocks}
    \end{align}

    To bound $\mathfrak k(0)$, let $f_1, f_2 \in L^2(\bar\mu)$ be real valued with $\mu(f_1) = \mu(f_2) = 0$, and let $\phi_1 := ( -\mathcal L )^{-1} f_1$ and $\phi_2 := ( -\mathcal L )^{-1} f_2$ be their Poisson solutions \eqref{eq: poisson solution}. Integrating \eqref{eq: reversal} for $\psi_1 = f_1$ and $\psi_2 = f_2$ over $t \Ge 0$ gives $\mu( \phi_1 f_2 ) = \mu( f_1 \phi_2 )$, and \eqref{eq: variance correlation} and \eqref{eq: normalized variance} give
    \begin{equation*}
        0 \Le \mu( \phi_1 f_1 ) \Le \frac12 \mathbb V(\gamma) \mu(f_1^2) .
    \end{equation*}
    Hence $(f_1,f_2) \mapsto \mu( \phi_1 f_2 )$ is a symmetric nonnegative bilinear form, and its Cauchy--Schwarz inequality for $f_2 = \langle\phi_1\rangle$, with the upper bound of $\mu( \phi_1 f_1 )$ and of $\mu( \phi_2 f_2 )$, gives
    \begin{equation*}
        \norm{\langle\phi_1\rangle}_{L^2(\mu)}^2 = \mu\bigl( \phi_1 \langle\phi_1\rangle \bigr) \Le \frac12 \mathbb V(\gamma) \norm{f_1}_{L^2(\mu)} \norm{\langle\phi_1\rangle}_{L^2(\mu)} .
    \end{equation*}
    For a function $f$ in the supremum \eqref{eq: axis position}, this bound of $\norm{\langle\phi_1\rangle}_{L^2(\mu)}$ for $f_1 = \mathrm{Re}\, f$ and for $f_1 = \mathrm{Im}\, f$ gives
    \begin{equation}
        \mathfrak k(0) \Le \frac12 \mathbb V(\gamma) .
        \label{eq: axis zero}
    \end{equation}
    For $\omega = 0$, the bounds \eqref{eq: axis blocks} and \eqref{eq: axis zero} give \eqref{eq: axis poisson}.

    To bound $\mathfrak k(\omega)$ for $\omega \neq 0$, let $f \in L^2(\bar\mu)$ with $\mu(f) = 0$, and let $\varphi = ( \mathrm i\omega - \mathcal L )^{-1} f$. The mean $\mu\bigl( ( \mathcal L \varphi )\, ( \overline\varphi \circ \iota ) \bigr)$ is real by \eqref{eq: generator reversal} for the real and imaginary parts of $\varphi$, and $\mu\bigl( \varphi\, ( \overline\varphi \circ \iota ) \bigr)$ is real since $\iota$ preserves $\mu$. Since $f$ is a function of the position and $\langle \varphi \circ \iota \rangle = \langle \varphi\rangle$, multiplying \eqref{eq: axis equation} for $\zeta = f$ by $\overline\varphi \circ \iota$ and taking the imaginary part of the mean gives
    \begin{equation*}
        \mathrm{Im}\, \mu\bigl( f\, \overline{\langle \varphi\rangle} \bigr) = \omega \mu\bigl( \varphi\, ( \overline\varphi \circ \iota ) \bigr) = \omega \Bigl( \norm{\langle \varphi\rangle}_{L^2(\mu)}^2 + \mu\bigl( ( \varphi - \langle \varphi\rangle )\, ( \overline{( \varphi - \langle \varphi\rangle )} \circ \iota ) \bigr) \Bigr) .
    \end{equation*}
    Since $\iota$ preserves $\mu$, this identity, the Cauchy--Schwarz inequality and the first inequality of \eqref{eq: axis cross} give
    \begin{equation*}
        |\omega| \norm{\langle \varphi\rangle}_{L^2(\mu)}^2 \Le \bigl| \mathrm{Im}\big\{ \mu\bigl( f\, \overline{\langle \varphi\rangle} \bigr) \big\} \bigr| + \frac{|\omega|}{\gamma} \mathrm{Re}\big\{ \mu\bigl( f\, \overline{\langle \varphi\rangle} \bigr) \big\} \Le \Bigl( 1 + \frac{\omega^2}{\gamma^2} \Bigr)^{\frac12} \norm{f}_{L^2(\mu)} \norm{\langle \varphi\rangle}_{L^2(\mu)} ,
    \end{equation*}
    and the supremum over $f$ gives
    \begin{equation}
        \mathfrak k(\omega) \Le \Bigl( \frac{1}{\omega^2} + \frac{1}{\gamma^2} \Bigr)^{\frac12} , \qquad \omega \neq 0 .
        \label{eq: axis frequency}
    \end{equation}

    To prove \eqref{eq: axis bound}, let $\varphi = ( \mathrm i\omega - \mathcal L )^{-1} \zeta$. By \eqref{eq: axis equation}, $-\mathcal L \varphi = \zeta - \mathrm i\omega \varphi$, so the uniqueness of the solution of \eqref{eq: axis equation}, with $0$ and $\zeta - \mathrm i\omega \varphi$ in place of $\omega$ and $\zeta$, gives $\varphi = ( -\mathcal L )^{-1} ( \zeta - \mathrm i\omega \varphi )$. Then the bound \eqref{eq: axis poisson} for $\zeta - \mathrm i\omega \varphi$ in place of $\zeta$ gives
    \begin{equation}
        \norm{\varphi}_{L^2(\mu)} \Le \Bigl( \frac12 \mathbb V(\gamma) + \frac1\gamma \Bigr) \bigl( \norm{\zeta}_{L^2(\mu)} + |\omega| \norm{\varphi}_{L^2(\mu)} \bigr) .
        \label{eq: axis implicit}
    \end{equation}
    If
    \begin{equation*}
        |\omega| \Bigl( \mathbb V(\gamma) + \frac{3}{2\gamma} \Bigr) \Le 1 ,
    \end{equation*}
    then the bound \eqref{eq: axis implicit} gives \eqref{eq: axis bound}. Otherwise, $\frac{\gamma}{|\omega|} < \gamma \mathbb V(\gamma) + \frac32$, and \eqref{eq: axis frequency} and \eqref{eq: axis blocks} give
    \begin{equation*}
        \gamma \norm{\varphi}_{L^2(\mu)} \Le \Bigl( 1 + \Bigl( 1 + \frac{\gamma^2}{\omega^2} \Bigr)^{\frac12} \Bigr) \norm{\zeta}_{L^2(\mu)} \Le ( \gamma \mathbb V(\gamma) + 3 ) \norm{\zeta}_{L^2(\mu)} ,
    \end{equation*}
    which is \eqref{eq: axis bound}.
\end{proof}

\subsection{Proof of Theorem~\ref{thm: mixing}}
\label{app: mixing proof}

We prove the lower bound of \eqref{eq: mixing bounds} from the definition \eqref{eq: mixing time} of the mixing time, and the upper bound from Lemma~\ref{lem: axis} and \citep[Theorem~1.3]{wei2020diffusion}.

\begin{proof}[Proof of Theorem~\ref{thm: mixing}]
    The definition \eqref{eq: mixing time}, the semigroup property and the contraction of $L^2(\mu)$ in Lemma~\ref{lem: dynamics} give, for every integer $k \Ge 0$, every $t \Ge k t_{\mathrm{mix}}(\gamma)$ and every $\psi \in L^2(\mu)$ with $\mu(\psi) = 0$,
    \begin{equation}
        \norm{e^{t\mathcal L}\psi}_{L^2(\mu)} \Le e^{-k} \norm{\psi}_{L^2(\mu)} .
        \label{eq: mixing powers}
    \end{equation}

    To prove $t_{\mathrm{mix}}(\gamma) \Ge \frac{1 - e^{-1}}{2} \mathbb V(\gamma)$, let $f \in L^2(\bar\mu)$ be nonconstant. The ratio of \eqref{eq: normalized variance} is unchanged when a constant is added to $f$, so we assume $\mu(f) = 0$. By \eqref{eq: mixing powers} for $\psi = f$, every integer $k \Ge 0$ and every $t \Ge k t_{\mathrm{mix}}(\gamma)$ satisfy $\mu\bigl( f e^{t\mathcal L}f \bigr) \Le e^{-k} \mu(f^2)$, so \eqref{eq: variance correlation} gives
    \begin{equation*}
        \frac{\sigma_f^2}{2} = \int_0^\infty \mu\bigl( f e^{t\mathcal L}f \bigr) \D t \Le \sum_{k=0}^\infty e^{-k} t_{\mathrm{mix}}(\gamma) \mu(f^2) ,
    \end{equation*}
    and the supremum over $f$ gives $t_{\mathrm{mix}}(\gamma) \Ge \frac{1 - e^{-1}}{2} \mathbb V(\gamma)$.

    To prove $t_{\mathrm{mix}}(\gamma) \Ge \frac{1 - e^{-1}}{2\gamma}$, we first prove that $U \in L^2(\bar\mu)$. Let $H(x,v) := U(x) - \inf U + \frac12 |v|^2 \Ge 0$. For an integer $n \Ge 1$, let $F_n := \exp\bigl( -\frac1n H \bigr)$, which satisfies $0 < F_n \Le 1$, and let $\Phi_n$ be the Poisson solution \eqref{eq: resolvent} of $F_n - \mu(F_n)$. By Lemma~\ref{lem: poisson}, $\Phi_n \in \mathrm{Dom}(\mathcal L)$ and $\mathcal L\Phi_n = \mu(F_n) - F_n$. Since $\mathcal A F_n = 0$ and $\nabla_v F_n = -\frac1n F_n v$, \eqref{eq: domain form} for $\varrho = \Phi_n$ and $\psi = F_n$ and \eqref{eq: domain adjoints} for $\varrho = \Phi_n$ and $W = F_n v$ give
    \begin{equation*}
        \mu\bigl( ( F_n - \mu(F_n) )^2 \bigr) = -\frac\gamma n \mu( \nabla_v\Phi_n \cdot F_n v ) = -\frac\gamma n \mu\Bigl( \Phi_n\, F_n \Bigl( \Bigl( 1 + \frac1n \Bigr) |v|^2 - d \Bigr) \Bigr) .
    \end{equation*}
    Hence the Cauchy--Schwarz inequality and Minkowski's integral inequality for \eqref{eq: resolvent} with \eqref{eq: mixing powers} give
    \begin{equation}
        n \norm{F_n - \mu(F_n)}_{L^2(\mu)} \Le \frac{\gamma\, t_{\mathrm{mix}}(\gamma)}{1 - e^{-1}} \Bigl\lVert F_n \Bigl( \Bigl( 1 + \frac1n \Bigr) |v|^2 - d \Bigr) \Bigr\rVert_{L^2(\mu)} .
        \label{eq: mixing truncation}
    \end{equation}
    The right side of \eqref{eq: mixing truncation} is bounded in $n$, since $t_{\mathrm{mix}}(\gamma) < \infty$ by \eqref{eq: mixing decay} and $0 < F_n \Le 1$. The left side of \eqref{eq: mixing truncation} is $\norm{G_n - \mu(G_n)}_{L^2(\mu)}$ for $G_n := n ( 1 - F_n )$, which increases pointwise to $H$ as $n \to \infty$. With an independent copy $Z' \sim \mu$ of $Z \sim \mu$, Fatou's lemma for $\mathbb E\bigl[ ( G_n(Z) - G_n(Z') )^2 \bigr] = 2 \norm{G_n - \mu(G_n)}_{L^2(\mu)}^2$ gives $\mathbb E\bigl[ ( H(Z) - H(Z') )^2 \bigr] < \infty$, hence $H \in L^2(\mu)$ and $U \in L^2(\bar\mu)$.

    Let $\Phi$ be the Poisson solution \eqref{eq: resolvent} of $F = U + \frac12 |v|^2 - \mu(U) - \frac d2$. By Lemma~\ref{lem: poisson}, $\Phi \in \mathrm{Dom}(\mathcal L)$ and $\mathcal L\Phi = -F$. Since $\mathcal A F = 0$ and $\nabla_v F = v$, \eqref{eq: domain form} for $\varrho = \Phi$ and $\psi = F$ and \eqref{eq: domain adjoints} for $\varrho = \Phi$ and $W = v$ give
    \begin{equation}
        \mu(F^2) = \gamma \mu( \nabla_v\Phi \cdot v ) = \gamma \mu\bigl( \Phi\, ( |v|^2 - d ) \bigr) .
        \label{eq: mixing energy}
    \end{equation}
    Since
    \begin{equation*}
        \mu\bigl( ( |v|^2 - d )^2 \bigr) = 2d , \qquad
        \mu(F^2) = \mu\bigl( ( U - \mu(U) )^2 \bigr) + \frac d2 \Ge \frac d2 ,
    \end{equation*}
    the Cauchy--Schwarz inequality in \eqref{eq: mixing energy} and Minkowski's integral inequality for \eqref{eq: resolvent} with \eqref{eq: mixing powers} give
    \begin{equation*}
        \sqrt{\frac d2} \norm{F}_{L^2(\mu)} \Le \mu(F^2) \Le \gamma \sqrt{2d} \norm{\Phi}_{L^2(\mu)} \Le \frac{\gamma \sqrt{2d}}{1 - e^{-1}} t_{\mathrm{mix}}(\gamma) \norm{F}_{L^2(\mu)} ,
    \end{equation*}
    and thus $t_{\mathrm{mix}}(\gamma) \Ge \frac{1 - e^{-1}}{2\gamma}$.

    To prove the upper bound of \eqref{eq: mixing bounds}, let the functions be complex valued, as in Appendix~\ref{app: axis}. By Lemma~\ref{lem: dynamics}, $(e^{t\mathcal L})_{t \Ge 0}$ is a strongly continuous semigroup of contractions on the space $\{ \psi \in L^2(\mu) : \mu(\psi) = 0 \}$, so the operator $-\mathcal L$ on this space is m-accretive. By the uniqueness of the solution of \eqref{eq: axis equation} and by \eqref{eq: axis bound} for $\zeta = \mathrm i\omega \varphi - \mathcal L \varphi$, every real number $\omega$ and every function $\varphi$ with $\mu(\varphi) = 0$ and $\mathrm{Re}\, \varphi, \mathrm{Im}\, \varphi \in \mathrm{Dom}(\mathcal L)$ satisfy
    \begin{equation*}
        \norm{\varphi}_{L^2(\mu)} \Le \Bigl( \mathbb V(\gamma) + \frac3\gamma \Bigr) \norm{\mathrm i\omega \varphi - \mathcal L \varphi}_{L^2(\mu)} .
    \end{equation*}
    Hence \citep[Theorem~1.3]{wei2020diffusion} gives, for every $\psi \in L^2(\mu)$ with $\mu(\psi) = 0$ and every $t \Ge 0$,
    \begin{equation}
        \norm{e^{t\mathcal L}\psi}_{L^2(\mu)} \Le \exp\Bigl( \frac\pi2 - \frac{t}{\mathbb V(\gamma) + \frac3\gamma} \Bigr) \norm{\psi}_{L^2(\mu)} .
        \label{eq: mixing rate}
    \end{equation}
    At $t = \bigl( 1 + \frac\pi2 \bigr) \bigl( \mathbb V(\gamma) + \frac3\gamma \bigr)$, the right side of \eqref{eq: mixing rate} is $e^{-1} \norm{\psi}_{L^2(\mu)}$, so the definition \eqref{eq: mixing time} gives the upper bound of \eqref{eq: mixing bounds}.
\end{proof}

\section{Proofs for the central limit theorem}
\label{app: clt proof}

This appendix proves Theorem~\ref{thm: clt}.

The distance $\mathbb W(f,t)$ of \eqref{eq: wasserstein} and $\norm{\nabla f}_\infty$ are unchanged when a constant is added to $f$, and for a constant $f$ both distributions of \eqref{eq: wasserstein} are the point mass at $0$. Hence, throughout this appendix, $f \in C^1(\mathcal X)$ is a nonconstant function with $\norm{\nabla f}_\infty < \infty$ and $\mu(f) = 0$, $\phi$ is the Poisson solution \eqref{eq: poisson solution} of $f$, and the asymptotic variance \eqref{eq: asymptotic variance} satisfies $\sigma_f^2 > 0$ by Corollary~\ref{cor: variance}, with $\sigma_f \Ge 0$. The time integral of $f$ is
\begin{equation}
    S_t := \int_0^t f(X_s) \,\D s .
    \label{eq: time integral}
\end{equation}
With the ratio $\alpha$ of \eqref{eq: ratio}, the reference rate and the scale of the distances are
\begin{equation}
    \lambda_0 := \frac{m\gamma}{( \gamma + \sqrt m + \sqrt K )^2} = \frac{m}{\gamma \alpha^2} , \qquad
    \mathfrak D_0 := \frac{\norm{\nabla f}_\infty}{\lambda_0 \sqrt m} = \frac{\gamma \alpha^2}{m^{\frac32}} \norm{\nabla f}_\infty ,
    \label{eq: scale}
\end{equation}
so that \eqref{eq: clt} reads $\sqrt t\, \mathbb W(f,t) \Le C ( 1 + \log\alpha )\, \mathfrak D_0$. By \eqref{eq: rate lower} for the rate $\lambda$ of \eqref{eq: l2 decay}, and by $\sqrt m \Le \gamma \alpha$, the reference rate $\lambda_0$ satisfies
\begin{equation}
    \lambda_0 \Le 3 \lambda , \qquad
    \lambda_0 \Le \gamma .
    \label{eq: rate bounds}
\end{equation}
Let $c_0$ be the number
\begin{equation}
    c_0 :=\frac{\sigma_f^2}{\lambda_0 \mathfrak D_0^2} .
    \label{eq: normalized sigma}
\end{equation}
The upper bound of \eqref{eq: variance} with $\gamma^2 + m + K \Le \gamma^2 \alpha^2$, and \eqref{eq: poisson moments} with $\gamma^2 + 2m + K \Le 2 \gamma^2 \alpha^2$, give
\begin{equation}
    c_0 \Le 2 , \qquad
    \norm{\phi}_{L^2(\mu)} \Le \frac{\sigma_f}{\sqrt{\lambda_0}} = \sqrt{c_0}\, \mathfrak D_0 .
    \label{eq: poisson scale}
\end{equation}

The parameter of this appendix is an even integer $p \Ge 6$, the exponent of the moment of $f$ in H\"older's inequality, and the moment factor is
\begin{equation}
    M := p \alpha^{\frac4p} .
    \label{eq: moment factor}
\end{equation}
No constant $C$ of this appendix depends on $p$.

\subsection{Moments of the test function and of the Poisson solutions}
\label{app: moments}

By \eqref{eq: variance correlation}, $\mu( 2\phi f - \sigma_f^2 ) = 0$, and the Poisson solution \eqref{eq: resolvent} of $2\phi f - \sigma_f^2$ is
\begin{equation}
    \Phi_2 := \int_0^\infty e^{r\mathcal L} ( 2\phi f - \sigma_f^2 ) \,\D r .
    \label{eq: second poisson}
\end{equation}
By Lemma~\ref{lem: poisson} for $F = f$ and for $F = 2\phi f - \sigma_f^2$,
\begin{equation}
    \mathcal L\phi = -f , \qquad
    \mathcal L\Phi_2 = \sigma_f^2 - 2\phi f , \qquad
    \mu(\phi) = \mu(\Phi_2) = 0 .
    \label{eq: poisson pair}
\end{equation}
Let $F_2$ be the function
\begin{equation}
    F_2 :=\Phi_2 f - \sigma_f^2\, \phi .
    \label{eq: third coefficient}
\end{equation}
In Lemma~\ref{lem: moments}, we bound $\mu(f^4)$, $\norm{f}_{L^p(\mu)}$, $\norm{\Phi_2}_{L^2(\mu)}$ and $|\mu(F_2)|$.

\begin{lemma}
\label{lem: moments}
    Let $p \Ge 6$ be an even integer. With the ratio $\alpha$ of \eqref{eq: ratio}, the reference rate $\lambda_0$ and the scale $\mathfrak D_0$ of \eqref{eq: scale} and the number $c_0$ of \eqref{eq: normalized sigma}, the test function $f$ satisfies
    \begin{equation}
        \mu(f^4) \Le C \alpha^2 c_0 ( \lambda_0 \mathfrak D_0 )^4 , \qquad
        \norm{f}_{L^p(\mu)} \Le C p \alpha^{\frac2p} c_0^{\frac1p} \lambda_0 \mathfrak D_0 .
        \label{eq: source moments}
    \end{equation}
    With the moment factor $M$ of \eqref{eq: moment factor}, the function $\Phi_2$ of \eqref{eq: second poisson} and the function $F_2$ of \eqref{eq: third coefficient} satisfy
    \begin{equation}
        \norm{\Phi_2}_{L^2(\mu)} \Le C M \sqrt{c_0}\, \mathfrak D_0^2 , \qquad
        |\mu(F_2)| \Le C M \sigma_f^2 \mathfrak D_0 .
        \label{eq: second poisson bounds}
    \end{equation}
\end{lemma}

\begin{proof}[Proof of Lemma~\ref{lem: moments}]
    We first verify that $f \in L^{2j}(\bar\mu)$ for every integer $j \Ge 1$. For $R > 0$, let $\chi_R \in C^\infty(\mathbb R)$ be odd with $|\chi_R(y)| \Le \min(|y|, R)$, $0 \Le \chi_R'(y) \Le 1$ and $\chi_R(y) = y$ for $|y| \Le R/2$, and let $g_R := \chi_R(f) \in C^1(\mathcal X) \cap L^\infty(\bar\mu)$, so $|\nabla g_R| \Le \norm{\nabla f}_\infty$. Applying \eqref{eq: poincare} to $g_R^j$ and the Cauchy--Schwarz inequality to $\mu(g_R^j)^2 \Le \mu(g_R^2)\mu(g_R^{2j-2})$ inductively in $j \Ge 1$ yields uniform bounds in $R$ on $\mu(g_R^{2j})$, and Fatou's lemma as $R \to \infty$ gives $f \in L^{2j}(\bar\mu)$ for every $j \Ge 1$.

    To prove \eqref{eq: source moments}, note that \eqref{eq: scale} gives $\frac{\norm{\nabla f}_\infty^2}{m} = ( \lambda_0 \mathfrak D_0 )^2$ and $\frac m\gamma = \alpha^2 \lambda_0$. Since $g = f^3$ satisfies $g, \nabla g \in L^2(\bar\mu)$ and $\mu\bigl( |\nabla g|^2 \bigr) \Le 9 \norm{\nabla f}_\infty^2\, \mu(f^4)$, \eqref{eq: variance lower} for $g$ and \eqref{eq: normalized sigma} give the bound of $\mu(f^4)$ in \eqref{eq: source moments},
    \begin{equation*}
        \mu(f^4) \Le \frac{9 \sigma_f^2}{2\gamma} \norm{\nabla f}_\infty^2 = \frac92 \alpha^2 c_0 ( \lambda_0 \mathfrak D_0 )^4 .
    \end{equation*}
    For an integer $j \Ge 3$, \eqref{eq: poincare} for $f$ and for $f^j$ and the Cauchy--Schwarz inequality give
    \begin{equation*}
        \mu(f^{2j}) \Le \mu(f^j)^2 + \frac{j^2}{m} \norm{\nabla f}_\infty^2\, \mu(f^{2j-2}) \Le ( 1 + j^2 ) ( \lambda_0 \mathfrak D_0 )^2 \mu(f^{2j-2}) .
    \end{equation*}
    Since $1 + j^2 \Le p^2$ for $j \Le \frac p2$, this bound for $j = 3, \dots, \frac p2$ gives $\mu(f^p) \Le ( p \lambda_0 \mathfrak D_0 )^{p-4} \mu(f^4)$, and the bound of $\mu(f^4)$ in \eqref{eq: source moments} then gives the bound of $\norm{f}_{L^p(\mu)}$.

    To prove the bound of $\norm{\Phi_2}_{L^2(\mu)}$ in \eqref{eq: second poisson bounds}, we first prove
    \begin{equation}
        \norm{\phi f}_{L^2(\mu)} \Le C M \sqrt{c_0}\, \lambda_0 \mathfrak D_0^2 .
        \label{eq: source bounds}
    \end{equation}
    By H\"older's inequality with the exponents $\frac{2p}{p-2}$ and $p$, and then between $L^2(\mu)$ and $L^4(\mu)$,
    \begin{equation*}
        \norm{\phi f}_{L^2(\mu)} \Le \norm{\phi}_{L^{\frac{2p}{p-2}}(\mu)} \norm{f}_{L^p(\mu)} \Le \norm{\phi}_{L^2(\mu)}^{1 - \frac4p} \norm{\phi}_{L^4(\mu)}^{\frac4p} \norm{f}_{L^p(\mu)} .
    \end{equation*}
    By \eqref{eq: poisson fourth}, \eqref{eq: scale} and the bound of $\mu(f^4)$ in \eqref{eq: source moments}, $\norm{\phi}_{L^4(\mu)} \Le C \sqrt\alpha\, c_0^{\frac14} \mathfrak D_0$. Hence the bound of $\norm{\phi}_{L^2(\mu)}$ in \eqref{eq: poisson scale} and the bound of $\norm{f}_{L^p(\mu)}$ in \eqref{eq: source moments} give
    \begin{equation*}
        \norm{\phi f}_{L^2(\mu)} \Le \bigl( \sqrt{c_0}\, \mathfrak D_0 \bigr)^{1 - \frac4p} \bigl( C \sqrt\alpha\, c_0^{\frac14} \mathfrak D_0 \bigr)^{\frac4p}\, C p \alpha^{\frac2p} c_0^{\frac1p} \lambda_0 \mathfrak D_0 \Le C p \alpha^{\frac4p} \sqrt{c_0}\, \lambda_0 \mathfrak D_0^2 ,
    \end{equation*}
    which is \eqref{eq: source bounds} by \eqref{eq: moment factor}.

    By \eqref{eq: variance correlation}, $2\mu(\phi f) = \sigma_f^2$, so $\norm{2\phi f - \sigma_f^2}_{L^2(\mu)} \Le 2 \norm{\phi f}_{L^2(\mu)}$. Hence the bound of $\norm{\Phi}_{L^2(\mu)}$ in \eqref{eq: resolvent bound} for $F = 2\phi f - \sigma_f^2$, \eqref{eq: rate bounds} and \eqref{eq: source bounds} give the bound of $\norm{\Phi_2}_{L^2(\mu)}$ in \eqref{eq: second poisson bounds}.

    To prove the bound of $|\mu(F_2)|$ in \eqref{eq: second poisson bounds}, note that $f \circ \iota = f$ for the flip $\iota$ of \eqref{eq: flip} and that $\iota$ preserves $\mu$. Hence \eqref{eq: poisson pair} and \eqref{eq: generator reversal} for $\varrho_1 = \Phi_2$ and $\varrho_2 = \phi$ give
    \begin{equation*}
        \mu(F_2) = \mu( \Phi_2 f ) = \mu\bigl( ( 2\phi f - \sigma_f^2 )\, ( \phi \circ \iota ) \bigr) = 2 \mu\bigl( \phi f\, ( \phi \circ \iota ) \bigr) .
    \end{equation*}
    Hence the Cauchy--Schwarz inequality with \eqref{eq: source bounds} and the bound of $\norm{\phi}_{L^2(\mu)}$ in \eqref{eq: poisson scale} gives $|\mu(F_2)| \Le C M c_0 \lambda_0 \mathfrak D_0^3$, which is the bound of $|\mu(F_2)|$ in \eqref{eq: second poisson bounds} by \eqref{eq: normalized sigma}.
\end{proof}

\subsection{Correlation of the time integral with the state}
\label{app: increments}

In Lemma~\ref{lem: decorrelation}, for a Lipschitz function $G$ and $\psi \in L^{\frac{2p}{p+2}}(\mu)$, we bound $\mathbb E\bigl[ G(S_s)\, \bigl( \psi(Z_s) - \mu(\psi) \bigr) \bigr]$ with the time integral $S_s$ of \eqref{eq: time integral}.

\begin{lemma}
\label{lem: decorrelation}
    Let $G \colon \mathbb R \to \mathbb R$ have Lipschitz constant at most $1$, let $s \Ge 0$ and let $p \Ge 6$ be an even integer. With the ratio $\alpha$ of \eqref{eq: ratio}, the scale $\mathfrak D_0$ of \eqref{eq: scale} and the number $c_0$ of \eqref{eq: normalized sigma}, every $\psi \in L^{\frac{2p}{p+2}}(\mu)$ satisfies
    \begin{equation}
        \bigl| \mathbb E\bigl[ G(S_s)\, \bigl( \psi(Z_s) - \mu(\psi) \bigr) \bigr] \bigr| \Le C \alpha^{\frac2p} c_0^{\frac12 - \frac1p} \mathfrak D_0 \norm{\psi}_{L^{\frac{2p}{p+2}}(\mu)} .
        \label{eq: decorrelation block}
    \end{equation}
\end{lemma}

\begin{proof}[Proof of Lemma~\ref{lem: decorrelation}]
    Let $n$ be the smallest integer with $n \Ge \lambda s$, for the rate $\lambda$ of \eqref{eq: l2 decay}, and let $r_j = \max\bigl( s - \frac{n-j}{\lambda}, 0 \bigr)$, $j = 0, \dots, n$, so that $r_0 = 0$, $r_n = s$ and $r_j - r_{j-1} \Le \frac1\lambda$. Since $S_0 = 0$ and $\mu( e^{s\mathcal L}\psi ) = \mu(\psi)$, the Markov property at the times $r_j$ gives
    \begin{equation}
        \mathbb E\bigl[ G(S_s)\, \bigl( \psi(Z_s) - \mu(\psi) \bigr) \bigr] = \sum_{j=1}^{n} \mathbb E\Bigl[ \bigl( G(S_{r_j}) - G(S_{r_{j-1}}) \bigr) \bigl( e^{(s-r_j)\mathcal L}\psi - \mu(\psi) \bigr)(Z_{r_j}) \Bigr] .
        \label{eq: decorrelation sum}
    \end{equation}
    By stationarity and $r_j - r_{j-1} \Le \frac1\lambda$, the upper bound of \eqref{eq: second moment}, Minkowski's integral inequality and \eqref{eq: rate bounds} give
    \begin{equation*}
        \mathbb E| S_{r_j} - S_{r_{j-1}} |^2 \Le \sigma_f^2 \Bigl( \frac1\lambda + \frac2\gamma \Bigr) \Le \frac{5 \sigma_f^2}{\lambda_0} , \qquad
        \bigl( \mathbb E| S_{r_j} - S_{r_{j-1}} |^4 \bigr)^{\frac14} \Le \frac1\lambda \norm{f}_{L^4(\mu)} \Le \frac{3}{\lambda_0} \norm{f}_{L^4(\mu)} .
    \end{equation*}
    Hence H\"older's inequality with \eqref{eq: normalized sigma} and the bound of $\mu(f^4)$ in \eqref{eq: source moments} gives
    \begin{equation}
        \bigl( \mathbb E| S_{r_j} - S_{r_{j-1}} |^{\frac{2p}{p-2}} \bigr)^{\frac{p-2}{2p}} \Le C \Bigl( \frac{\sigma_f^2}{\lambda_0} \Bigr)^{\frac12 - \frac2p} \Bigl( \frac{1}{\lambda_0} \norm{f}_{L^4(\mu)} \Bigr)^{\frac4p} \Le C \alpha^{\frac2p} c_0^{\frac12 - \frac1p} \mathfrak D_0 .
        \label{eq: block increment}
    \end{equation}
    The Riesz--Thorin interpolation theorem between \eqref{eq: lp decay} in $L^{\frac43}(\mu)$ and \eqref{eq: l2 contraction} in $L^2(\mu)$, with the weights $\frac4p$ and $1 - \frac4p$, gives, for every $t \Ge 0$,
    \begin{equation}
        \norm{e^{t\mathcal L}\psi - \mu(\psi)}_{L^{\frac{2p}{p+2}}(\mu)} \Le 2^{\frac2p} \min\bigl( 1, e^{\frac12 - \lambda t} \bigr)^{1 - \frac2p} \norm{\psi}_{L^{\frac{2p}{p+2}}(\mu)} .
        \label{eq: lu decay}
    \end{equation}
    Since $s - r_j = \frac{n-j}{\lambda}$ for $j \Ge 1$ and $1 - \frac2p \Ge \frac23$, H\"older's inequality with the exponents $\frac{2p}{p-2}$ and $\frac{2p}{p+2}$, \eqref{eq: block increment} and \eqref{eq: lu decay} at $t = s - r_j$ give
    \begin{equation*}
        \Bigl| \mathbb E\Bigl[ \bigl( G(S_{r_j}) - G(S_{r_{j-1}}) \bigr) \bigl( e^{(s-r_j)\mathcal L}\psi - \mu(\psi) \bigr)(Z_{r_j}) \Bigr] \Bigr| \Le C \alpha^{\frac2p} c_0^{\frac12 - \frac1p} \mathfrak D_0 e^{-\frac23 (n-j)} \norm{\psi}_{L^{\frac{2p}{p+2}}(\mu)} ,
    \end{equation*}
    and the sum over $j$ in \eqref{eq: decorrelation sum} gives \eqref{eq: decorrelation block}.
\end{proof}

\subsection{Kernels with a cubic phase}
\label{app: kernels}

In Lemma~\ref{lem: kernels}, we bound the integrals of the absolute values of the derivatives of the kernels whose Fourier transform is a Gaussian function with a cubic phase. For $\tau \in \mathbb R$ and $y \in \mathbb R$, the kernel is
\begin{equation}
    \mathfrak p_\tau(y) := \frac{1}{2\pi} \int_{\mathbb R} \exp\Bigl( \mathrm i y \xi - \frac{\xi^2}{2} - \mathrm i \tau \xi^3 \Bigr) \D\xi ,
    \label{eq: cubic kernel}
\end{equation}
where $\mathrm i$ is the imaginary unit. The function $\mathfrak p_\tau$ is real and smooth and decreases rapidly with all its derivatives, and $\mathfrak p_0$ is the density of $\mathcal N(0,1)$. Differentiation under the integral gives
\begin{equation}
    \partial_\tau \mathfrak p_\tau = \partial_y^3 \mathfrak p_\tau .
    \label{eq: phase derivative}
\end{equation}
Up to the sign of $y$, \eqref{eq: phase derivative} is the linearized Korteweg--de Vries (KdV) equation, and $\mathfrak p_\tau$ is its solution with the value $\mathfrak p_0$ at $\tau = 0$. For an integer $j \Ge 1$ and $\tau \in \mathbb R$, let $\mathfrak N_j(\tau)$ be the integral
\begin{equation}
    \mathfrak N_j(\tau) := \int_{\mathbb R} \bigl| \partial_y^j \mathfrak p_\tau(y) \bigr| \D y .
    \label{eq: kernel norms}
\end{equation}

\begin{lemma}
\label{lem: kernels}
    For $j = 1, \dots, 5$ and every $\tau \in \mathbb R$, the integral \eqref{eq: kernel norms} satisfies
    \begin{equation}
        \mathfrak N_j(\tau) \Le C ( 1 + |\tau| ) .
        \label{eq: kernel second}
    \end{equation}
\end{lemma}

\begin{proof}[Proof of Lemma~\ref{lem: kernels}]
    The Cauchy--Schwarz inequality with the weight $1 + y^2$ gives
    \begin{equation}
        \mathfrak N_j(\tau)^2 \Le \pi \int_{\mathbb R} ( 1 + y^2 )\, \bigl| \partial_y^j \mathfrak p_\tau \bigr|^2 \D y .
        \label{eq: kernel weights}
    \end{equation}
    The Fourier transform of $\partial_y^j \mathfrak p_\tau$ is $( \mathrm i \xi )^j \exp\bigl( -\frac{\xi^2}{2} - \mathrm i \tau \xi^3 \bigr)$, and multiplying a function by $y$ applies the operator $\mathrm i\, \partial_\xi$ to its Fourier transform. Hence Plancherel's identity gives
    \begin{equation*}
        \int_{\mathbb R} ( 1 + y^2 )\, \bigl| \partial_y^j \mathfrak p_\tau \bigr|^2 \D y = \frac{1}{2\pi} \int_{\mathbb R} \xi^{2j} e^{-\xi^2} \Bigl( 1 + \Bigl( \frac j\xi - \xi \Bigr)^2 + 9 \tau^2 \xi^4 \Bigr) \D\xi \Le C ( 1 + \tau^2 ) ,
    \end{equation*}
    and \eqref{eq: kernel weights} gives \eqref{eq: kernel second}.
\end{proof}

\subsection{Proof of Theorem~\ref{thm: clt}}
\label{app: comparison}

The proof of Theorem~\ref{thm: clt} compares the mean of a Lipschitz function of $S_t$ with its mean under $\mathcal N(0,\sigma_f^2 t)$, through a solution of a heat equation with a term of third order. The Poisson solution $\phi$ removes the term with the first derivative of the solution of the heat equation, the Poisson solution $\Phi_2$ of $2\phi f - \sigma_f^2$ removes the terms with the second derivative, and the term of third order of the equation removes the mean of the coefficient $\frac12 F_2$ of the third derivative. The remaining terms are correlations of $\phi$, of $\Phi_2$ and of $F_2$ with Lipschitz functions of the time integral, and Lemma~\ref{lem: decorrelation} bounds these correlations by time blocks. The proof takes $p$ with $M \Le C ( 1 + \log\alpha )$, which gives the factor $1 + \log\alpha$ of \eqref{eq: clt}.

Let $\chi \sim \mathcal N(0,1)$. By \eqref{eq: wasserstein} with $\mu(f) = 0$ and the Kantorovich--Rubinstein duality \citep[Theorem~1.14]{villani2021topics},
\begin{equation}
    \sqrt t\, \mathbb W(f,t) = \mathcal W_1\bigl( \mathrm{Law}(S_t), \mathcal N(0,\sigma_f^2 t) \bigr) = \sup_\ell \Bigl( \mathbb E[\ell(S_t)] - \mathbb E\bigl[ \ell(\sigma_f\sqrt t\, \chi) \bigr] \Bigr) ,
    \label{eq: distance}
\end{equation}
where the supremum is over the functions $\ell \colon \mathbb R \to \mathbb R$ with Lipschitz constant at most $1$.

\begin{proof}[Proof of Theorem~\ref{thm: clt}]
    Let $p \Ge 6$ be an even integer. With the moment factor $M$ of \eqref{eq: moment factor}, we first prove
    \begin{equation}
        \sqrt t\, \mathbb W(f,t) \Le C M \mathfrak D_0 , \qquad \forall t > 0 .
        \label{eq: distance scalar}
    \end{equation}
    Let $s \Ge 0$, and let $G \colon \mathbb R \to \mathbb R$ have Lipschitz constant at most $1$. By \eqref{eq: moment factor}, $\alpha^{\frac2p} \Le M$, and $\norm{\psi}_{L^{\frac{2p}{p+2}}(\mu)} \Le \norm{\psi}_{L^2(\mu)}$ for every $\psi \in L^2(\mu)$. Since $\mu(\phi) = \mu(\Phi_2) = 0$ by \eqref{eq: poisson pair}, the bound \eqref{eq: decorrelation block} for $\psi = \phi$ and for the function $\psi = \Phi_2$ of \eqref{eq: second poisson}, with \eqref{eq: poisson scale} and the bound of $\norm{\Phi_2}_{L^2(\mu)}$ in \eqref{eq: second poisson bounds}, gives
    \begin{equation}
        \bigl| \mathbb E\bigl[ G(S_s)\, \phi(Z_s) \bigr] \bigr| \Le C M \mathfrak D_0^2 , \qquad
        \bigl| \mathbb E\bigl[ G(S_s)\, \Phi_2(Z_s) \bigr] \bigr| \Le C M^2 \mathfrak D_0^3 .
        \label{eq: corrector correlation}
    \end{equation}
    By H\"older's inequality with the exponents $2$ and $p$, the bound of $\norm{\Phi_2}_{L^2(\mu)}$ in \eqref{eq: second poisson bounds} and the bound of $\norm{f}_{L^p(\mu)}$ in \eqref{eq: source moments},
    \begin{equation*}
        \norm{\Phi_2 f}_{L^{\frac{2p}{p+2}}(\mu)} \Le \norm{\Phi_2}_{L^2(\mu)} \norm{f}_{L^p(\mu)} \Le C M p \alpha^{\frac2p} c_0^{\frac12 + \frac1p} \lambda_0 \mathfrak D_0^3 .
    \end{equation*}
    Hence, with $p \alpha^{\frac4p} = M \Ge 1$ and $c_0 \lambda_0 \mathfrak D_0^2 = \sigma_f^2$, \eqref{eq: decorrelation block} for $\psi = \Phi_2 f$ and the bound of $\mathbb E\bigl[ G(S_s)\, \phi(Z_s) \bigr]$ in \eqref{eq: corrector correlation} give, for the function $F_2$ of \eqref{eq: third coefficient},
    \begin{equation}
        \bigl| \mathbb E\bigl[ G(S_s)\, \bigl( F_2(Z_s) - \mu(F_2) \bigr) \bigr] \bigr| \Le C M^2 \sigma_f^2 \mathfrak D_0^2 .
        \label{eq: third correlation}
    \end{equation}

    Fix $t > 0$ and a function $\ell \colon \mathbb R \to \mathbb R$ with Lipschitz constant at most $1$. For $s \in [0,t]$ and $y \in \mathbb R$, let $r_s$ and $\tau_s$ be the numbers and $\Psi$ the function
    \begin{equation}
        r_s := \bigl( \sigma_f^2 (t-s) + M^2 \mathfrak D_0^2 \bigr)^{\frac12} , \qquad
        \tau_s := \frac{\mu(F_2)\, (t-s)}{2 r_s^3} , \qquad
        \Psi(s,y) := \int_{\mathbb R} \ell( y - r_s y' )\, \mathfrak p_{\tau_s}(y') \,\D y' ,
        \label{eq: cubic function}
    \end{equation}
    with the kernels of \eqref{eq: cubic kernel}. The function $\Psi(s,\cdot)$ is the convolution of $\ell$ with the kernel whose Fourier transform is $\xi \mapsto \exp\bigl( -\frac12 r_s^2 \xi^2 - \frac{\mathrm i}{2} \mu(F_2)\, (t-s)\, \xi^3 \bigr)$. Hence $\Psi$ is smooth on $[0,t] \times \mathbb R$ and solves
    \begin{equation}
        \partial_s\Psi + \frac{\sigma_f^2}{2}\, \partial_y^2\Psi + \frac{\mu(F_2)}{2}\, \partial_y^3\Psi = 0 .
        \label{eq: cubic equation}
    \end{equation}
    Working with the continuous semimartingale modifications of $s \mapsto \phi(Z_s)$ and $s \mapsto \Phi_2(Z_s)$ from Lemma~\ref{lem: domain}, \eqref{eq: poisson pair}, \eqref{eq: cubic equation} and \eqref{eq: domain ito} for $\varrho = \phi$ and $\varrho = \Phi_2$, It\^o's formula for $\Psi(s,S_s)$, for $\phi(Z_s)\, \partial_y\Psi(s,S_s)$ and for $\Phi_2(Z_s)\, \partial_y^2\Psi(s,S_s)$ gives
    \begin{align}
        \mathbb E\bigl[ \Psi(t,S_t) \bigr] - \Psi(0,0) & = \mathbb E\int_0^t \Bigl( f\, \partial_y\Psi - \frac{\sigma_f^2}{2}\, \partial_y^2\Psi - \frac{\mu(F_2)}{2}\, \partial_y^3\Psi \Bigr) \D s ,
        \label{eq: ito heat} \\
        \mathbb E\bigl[ \phi(Z_t)\, \partial_y\Psi(t,S_t) \bigr] & = \mathbb E\int_0^t \Bigl( -f\, \partial_y\Psi + \phi f\, \partial_y^2\Psi - \frac{\sigma_f^2}{2}\, \phi\, \partial_y^3\Psi - \frac{\mu(F_2)}{2}\, \phi\, \partial_y^4\Psi \Bigr) \D s ,
        \label{eq: ito first} \\
        \mathbb E\bigl[ \Phi_2(Z_t)\, \partial_y^2\Psi(t,S_t) \bigr] & = \mathbb E\int_0^t \Bigl( ( \sigma_f^2 - 2\phi f )\, \partial_y^2\Psi + \Phi_2 f\, \partial_y^3\Psi - \frac{\sigma_f^2}{2}\, \Phi_2\, \partial_y^4\Psi - \frac{\mu(F_2)}{2}\, \Phi_2\, \partial_y^5\Psi \Bigr) \D s ,
        \label{eq: ito second}
    \end{align}
    where the functions on $T\mathcal X$ are evaluated at $Z_s$ and the derivatives of $\Psi$ at $(s,S_s)$. With the function $F_2$ of \eqref{eq: third coefficient}, the sum of \eqref{eq: ito heat}, \eqref{eq: ito first} and half of \eqref{eq: ito second} is
    \begin{equation}
    \begin{aligned}
        & \mathbb E\bigl[ \Psi(t,S_t) \bigr] - \Psi(0,0) = -\mathbb E\bigl[ \phi(Z_t)\, \partial_y\Psi(t,S_t) \bigr] - \frac12\, \mathbb E\bigl[ \Phi_2(Z_t)\, \partial_y^2\Psi(t,S_t) \bigr] \\
        & \quad + \frac12 \int_0^t \mathbb E\bigl[ \bigl( F_2(Z_s) - \mu(F_2) \bigr)\, \partial_y^3\Psi(s,S_s) \bigr] \D s - \int_0^t \mathbb E\Bigl[ \Bigl( \frac{\sigma_f^2}{4}\, \Phi_2(Z_s) + \frac{\mu(F_2)}{2}\, \phi(Z_s) \Bigr) \partial_y^4\Psi(s,S_s) \Bigr] \D s \\
        & \quad - \frac{\mu(F_2)}{4} \int_0^t \mathbb E\bigl[ \Phi_2(Z_s)\, \partial_y^5\Psi(s,S_s) \bigr] \D s .
    \end{aligned}
        \label{eq: heat identity}
    \end{equation}

    To bound the terms of the right side of \eqref{eq: heat identity}, note that \eqref{eq: cubic function} gives $r_s \Ge r_t = M \mathfrak D_0$ and $\sigma_f^2 (t-s) \Le r_s^2$. Hence \eqref{eq: cubic function} and the bound of $|\mu(F_2)|$ in \eqref{eq: second poisson bounds} give
    \begin{equation}
        |\tau_s| \Le \frac{|\mu(F_2)|}{2 \sigma_f^2\, r_t} \Le C , \qquad
        \int_0^t r_s^{-j} \,\D s \Le \frac{2}{\sigma_f^2\, r_t^{j-2}} , \quad j = 3, 4, 5 .
        \label{eq: cubic variable}
    \end{equation}
    Since $\ell$ has Lipschitz constant at most $1$ and $|\tau_s| \Le C$, \eqref{eq: cubic function}, \eqref{eq: kernel norms} and \eqref{eq: kernel second} give, for $j = 1, \dots, 5$, $0 \Le s \Le t$ and $y \in \mathbb R$,
    \begin{equation}
        \bigl| \partial_y^{j+1}\Psi(s,y) \bigr| \Le \mathfrak N_j( \tau_s )\, r_s^{-j} \Le C r_s^{-j} ,
        \label{eq: cubic lipschitz}
    \end{equation}
    which bounds the Lipschitz constant of $\partial_y^j\Psi(s,\cdot)$. By \eqref{eq: corrector correlation} and \eqref{eq: cubic lipschitz} at $s = t$,
    \begin{equation}
        \bigl| \mathbb E\bigl[ \phi(Z_t)\, \partial_y\Psi(t,S_t) \bigr] \bigr| + \frac12 \bigl| \mathbb E\bigl[ \Phi_2(Z_t)\, \partial_y^2\Psi(t,S_t) \bigr] \bigr| \Le C \mathfrak D_0 .
        \label{eq: terminal correctors}
    \end{equation}
    By the bound of $|\mu(F_2)|$ in \eqref{eq: second poisson bounds}, \eqref{eq: corrector correlation}, \eqref{eq: third correlation} and the bounds \eqref{eq: cubic lipschitz}, and then by the bounds of $\int_0^t r_s^{-j} \,\D s$ in \eqref{eq: cubic variable},
    \begin{equation}
    \begin{aligned}
        \frac12 \int_0^t \bigl| \mathbb E\bigl[ \bigl( F_2(Z_s) - \mu(F_2) \bigr)\, \partial_y^3\Psi(s,S_s) \bigr] \bigr| \D s & \Le C M^2 \sigma_f^2 \mathfrak D_0^2 \int_0^t r_s^{-3} \,\D s \Le C M \mathfrak D_0 , \\
        \int_0^t \Bigl| \mathbb E\Bigl[ \Bigl( \frac{\sigma_f^2}{4}\, \Phi_2(Z_s) + \frac{\mu(F_2)}{2}\, \phi(Z_s) \Bigr) \partial_y^4\Psi(s,S_s) \Bigr] \Bigr| \D s & \Le C M^2 \sigma_f^2 \mathfrak D_0^3 \int_0^t r_s^{-4} \,\D s \Le C \mathfrak D_0 , \\
        \frac{|\mu(F_2)|}{4} \int_0^t \bigl| \mathbb E\bigl[ \Phi_2(Z_s)\, \partial_y^5\Psi(s,S_s) \bigr] \bigr| \D s & \Le C M^3 \sigma_f^2 \mathfrak D_0^4 \int_0^t r_s^{-5} \,\D s \Le C \mathfrak D_0 .
    \end{aligned}
        \label{eq: cubic integrals}
    \end{equation}

    By \eqref{eq: cubic function}, $\tau_t = 0$, so $\Psi(t,y) = \mathbb E\bigl[ \ell( y - r_t \chi ) \bigr]$. Since $\ell$ has Lipschitz constant at most $1$ and $\mathbb E|\chi| \Le 1$,
    \begin{equation}
        \bigl| \mathbb E[\ell(S_t)] - \mathbb E\bigl[ \Psi(t,S_t) \bigr] \bigr| \Le r_t\, \mathbb E|\chi| \Le M \mathfrak D_0 .
        \label{eq: end comparison}
    \end{equation}
    By \eqref{eq: phase derivative} and an integration by parts, every $\tau \in \mathbb R$ satisfies
    \begin{equation*}
        \Bigl| \frac{\D}{\D\tau} \int_{\mathbb R} \ell( -r_0 y' )\, \mathfrak p_\tau(y') \,\D y' \Bigr| \Le r_0\, \mathfrak N_2(\tau) ,
    \end{equation*}
    and $\int_{\mathbb R} \ell( -r_0 y' )\, \mathfrak p_\tau(y') \,\D y'$ is $\mathbb E[ \ell( r_0 \chi ) ]$ at $\tau = 0$ and $\Psi(0,0)$ at $\tau = \tau_0$. By \eqref{eq: cubic variable}, $|\tau_0| \Le C$, so \eqref{eq: kernel second}, \eqref{eq: cubic function} and the bound of $|\mu(F_2)|$ in \eqref{eq: second poisson bounds} give
    \begin{equation}
        \bigl| \Psi(0,0) - \mathbb E[ \ell( r_0 \chi ) ] \bigr| \Le C r_0 |\tau_0| = \frac{C |\mu(F_2)|\, t}{2 r_0^2} \Le \frac{C |\mu(F_2)|}{\sigma_f^2} \Le C M \mathfrak D_0 .
        \label{eq: start comparison}
    \end{equation}
    Moreover, $r_0^2 = \sigma_f^2 t + r_t^2$ gives $r_0 - \sigma_f\sqrt t \Le r_t$, so
    \begin{equation}
        \bigl| \mathbb E[ \ell( r_0 \chi ) ] - \mathbb E\bigl[ \ell( \sigma_f\sqrt t\, \chi ) \bigr] \bigr| \Le \bigl( r_0 - \sigma_f\sqrt t \bigr)\, \mathbb E|\chi| \Le M \mathfrak D_0 .
        \label{eq: variance comparison}
    \end{equation}
    Since $M \Ge 1$, the identity \eqref{eq: heat identity} and the bounds \eqref{eq: terminal correctors}--\eqref{eq: variance comparison} give
    \begin{equation*}
        \mathbb E[\ell(S_t)] - \mathbb E\bigl[ \ell( \sigma_f\sqrt t\, \chi ) \bigr] \Le C M \mathfrak D_0 ,
    \end{equation*}
    and the supremum over $\ell$ in \eqref{eq: distance} gives \eqref{eq: distance scalar}.

    To prove \eqref{eq: clt}, let $p$ be the smallest even integer with $p \Ge 6$ and $p \Ge 4 \log\alpha$. Then $\alpha^{\frac4p} \Le e$ and $p \Le 8 ( 1 + \log\alpha )$, so \eqref{eq: moment factor} gives $M \Le C ( 1 + \log\alpha )$, and \eqref{eq: distance scalar} with \eqref{eq: scale} gives \eqref{eq: clt}.
\end{proof}


\end{document}